# EXTENDED 3-DIMENSIONAL BORDISM AS THE THEORY OF MODULAR OBJECTS

BRUCE BARTLETT, CHRISTOPHER L. DOUGLAS,
CHRISTOPHER J. SCHOMMER-PRIES, AND JAMIE VICARY

Abstract. A modular object in a symmetric monoidal bicategory is a Frobenius algebra object whose product and coproduct are biadjoint, equipped with a braided structure and a compatible twist, satisfying rigidity, ribbon, pivotality, and modularity conditions. We prove that the oriented 3-dimensional bordism bicategory of 1-, 2-, and 3-manifolds is the free symmetric monoidal bicategory on a single anomaly-free modular object.

Contents

## 1. Introduction

This is the second paper in a series [7, 8, 9] on a generators-and-relations approach to 3-dimensional topological quantum field theories extended to 1-manifolds. Compact oriented manifolds of dimensions 1, 2, and 3 can be organized into a symmetric monoidal bicategory $\mathbf{Bord}^{\mathrm{or}}_{123}$, as follows:

- objects are closed oriented 1-dimensional manifolds;
- 1-morphisms are 2-dimensional compact oriented bordisms;
- 2-morphisms are diffeomorphism classes of compact oriented 3-dimensional bordisms.

In our first paper [7] we use higher Cerf theory to give a presentation $\mathcal{F}$ of $\mathbf{Bord}^{\mathrm{or}}_{123}$ which is finite, but very large, containing 132 relations. For reference, we repeat this presentation in Appendix A. In this paper we algebraically simplify the presentation $\mathcal{F}$ to obtain an equivalent presentation $\mathcal{O}$ having only 33 relations, all of which have clear conceptual interpretations from a categorical and topological point of view (namely as the structure of an anomaly-free, modular, ribbon, pivotal, rigid, biadjoint Frobenius structure).

Given a presentation $\mathcal{H}$, we write $\mathbf{F}(\mathcal{H})$ for the free symmetric monoidal bicategory generated by $\mathcal{H}$. Our main result in this paper can then be stated as follows.

**Theorem 1.** *There is a symmetric monoidal equivalence $\mathbf{F}(\mathcal{F}) \simeq \mathbf{F}(\mathcal{O})$ between the bicategories generated by the 2-Morse presentation $\mathcal{F}$ and the anomaly-free modular presentation $\mathcal{O}$.*

We can compose Theorem 1 with the main result of [7] that there is a symmetric monoidal equivalence $\mathbf{Bord}_{123}^{\mathrm{or}} \simeq \mathbf{F}(\mathcal{F})$, to give the following result.

**Corollary 2.** *There is a symmetric monoidal equivalence $\mathbf{Bord}_{123}^{\mathrm{or}} \simeq \mathbf{F}(\mathcal{O})$ between the oriented 3-dimensional bordism bicategory and the bicategory generated by the anomaly-free modular presentation $\mathcal{O}$.*

The Kirby calculus [29] gives a simple set of moves relating surgery descriptions of closed oriented 3-manifolds by framed links in $S^3$, replacing the relations naturally coming from Cerf theory; analogously, Theorem 1 gives a simple presentation of extended bordisms, replacing the presentation naturally coming from higher Cerf theory. Also, the equivalence in Corollary 2 may be seen as a categorification of the classic result (see [31] and references therein) that the symmetric monoidal category of closed oriented 1-dimensional manifolds and diffeomorphism classes of compact oriented 2-dimensional bordisms is symmetric-monoidally equivalent to the free symmetric monoidal category on a single commutative Frobenius algebra object.

In our third paper [8] we show how to modify the presentation $\mathcal{O}$ to obtain presentations of certain extensions of the oriented bordism bicategory $\mathbf{Bord}_{123}^{\mathrm{or}}$, corresponding geometrically to manifolds with signature, componentwise signature, or $p_1$ structure. These extra structures on manifolds have played an important role in topological quantum field theory [1, 10, 16, 35, 39, 43].

Finally, in our fourth paper [9] we use these presentations to classify $k$-linear representations of $\mathbf{Bord}_{123}^{\mathrm{or}}$ and its extensions, where $k$ is an algebraically-closed field, in terms of modular tensor categories equipped with extra structure.

Our approach of using generators and relations to relate 3-dimensional topological quantum field theories and modular categories has much in common with other work [2, 3, 4, 10, 14, 16, 32, 35, 36, 37, 39, 43]. Since this present paper is about algebraically simplifying the presentation $\mathcal{F}$ from [7], and it does not deal with topological quantum field theories or modular categories directly, we leave discussion of the relationship with the literature to the papers [7, 8, 9].

1.1. **Outline.** In Section 2 we define the ribbon presentation $\mathcal{R}$, the modular presentation $\mathcal{M}$, and the anomaly-free presentation $\mathcal{O}$, which have the same generators but progressively more relations. In Section 3 we work with the ribbon presentation $\mathcal{R}$ and show how this set of relations implies that $\mathbf{F}(\mathcal{R})$ (and hence also $\mathbf{F}(\mathcal{M})$ and $\mathbf{F}(\mathcal{O})$) is a *pivotal bicategory*. In Section 4 we work with the modular presentation $\mathcal{M}$ and show how the Modularity relation (20) implies the algebraic versions of invertibility of the $s$-matrix, and charge conjugation. In Section 5 we define the anomaly 2-morphism $x$ and derive its basic properties, leading to a proof of the anomalous braid relation for Dehn twists (54). In Section 6 we show that the relations in $\mathcal{O}$ imply those in $\mathcal{F}$, and in Section 7 we show the converse, leading to a proof of Theorem 1 in Section 7.10.

In Appendix A we recall the 2-Morse presentation $\mathcal{F}$ from [7]. In Appendix B we give our conventions regarding *daggers* and *rotations* of relations. In Appendix C we review the notions of quasistrict symmetric monoidal bicategories and presentations of such from [40], as well as the wire diagrams graphical calculus from [5] specialized to the present context.

1.2. **Acknowledgements.** We are grateful to John Baez, who encouraged and supported this project from the beginning, and provided key insights and inspiration. We also thank Sylvain Gervais, Mark Lackenby, Justin Roberts, Daniel Ruberman and Kevin Walker for helpful discussions.

## 2. The modular presentation

In this section we define the *ribbon presentation* $\mathcal{R}$, the *modular presentation* $\mathcal{M}$ and the *anomaly-free modular presentation* $\mathcal{O}$. These presentations have the same generators but progressively more relations, and it is useful to investigate the consequences of each relation in turn. Note that equations which hold in a weaker presentation continue to hold in a strengthening of it, so that for instance every equation between 2-morphisms in $\mathbf{F}(\mathcal{R})$ also holds in $\mathbf{F}(\mathcal{M})$ and $\mathbf{F}(\mathcal{O})$.

In the paper [9] we show that a linear representation of $\mathbf{F}(\mathcal{M})$ corresponds to a modular category equipped with a square root of the anomaly in each factor, and a linear representation of $\mathbf{F}(\mathcal{O})$ corresponds to a modular category whose anomaly is the identity in each factor. (We allow the unit object of our modular tensor category to have proper subobjects, and by 'factor' we mean a summand of the tensor category with a simple unit object.)

In this paper, whenever we refer to a 'presentation' of a symmetric monoidal bicategory, we are referring more precisely to a 3-computad for a quasistrict symmetric monoidal bicategory. For a review of 3-computads and the quasistrict symmetric monoidal bicategories presented by them, see Appendix C, which is a summary of the relevant material from [40]. Appendix C also reviews the graphical notation for quasistrict symmetric monoidal bicategories from [5], adapted to our context, which we will use throughout this paper. For our conventions regarding *daggers* and *rotations* of relations, see Appendix B.

**Definition 3.** The *ribbon presentation* $\mathcal{R}$ is the presentation defined as follows:

- Generating object:

- Generating 1-morphisms:

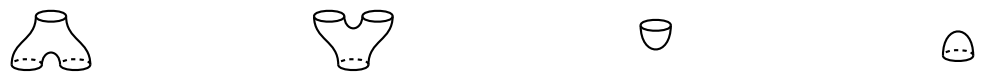

- Invertible generating 2-morphisms:

(1) 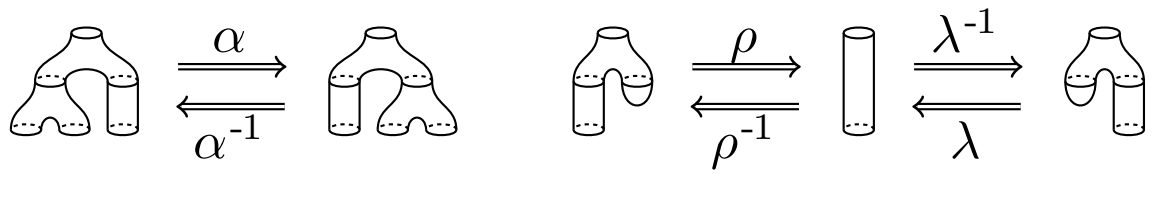

(2) 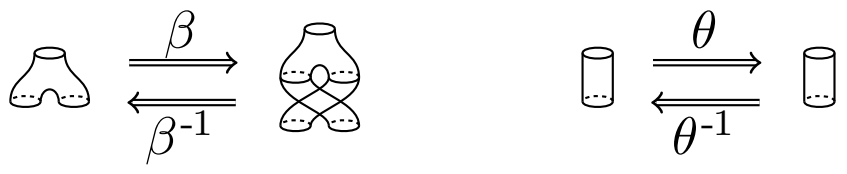

Noninvertible generating 2-morphisms:

(3) 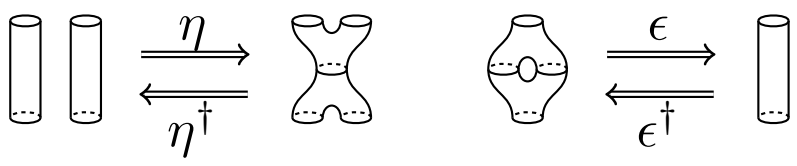

(4) 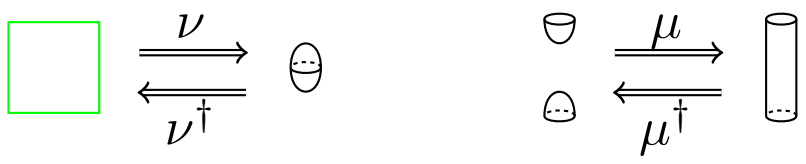

The relations are as follows:

- (Inverses) Each of the invertible generating 2-morphisms $\omega$ satisfies $\omega \circ \omega^{-1} = \mathrm{id}$ and $\omega^{-1} \circ \omega = \mathrm{id}$.
- (Monoidal) The generating 2-morphisms in (1) obey the pentagon and unit equations:

(5) 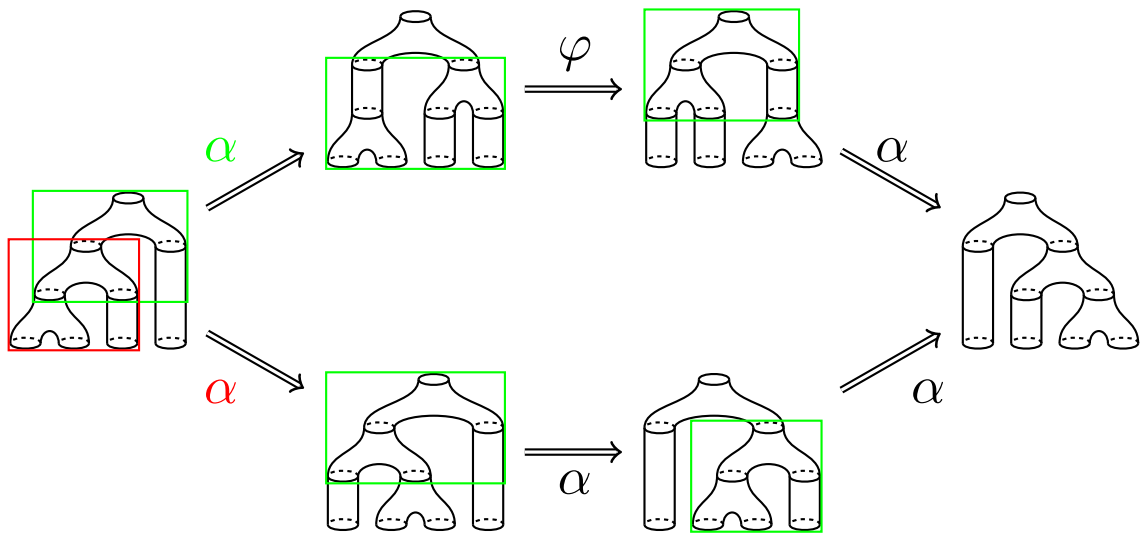

(6) 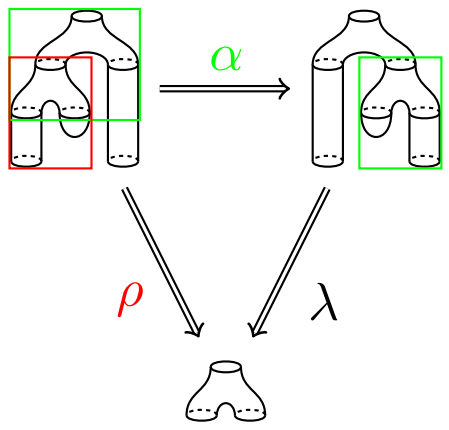

- (Balanced) The data (1) and (2) forms a braided monoidal object equipped with a compatible twist:

(7) 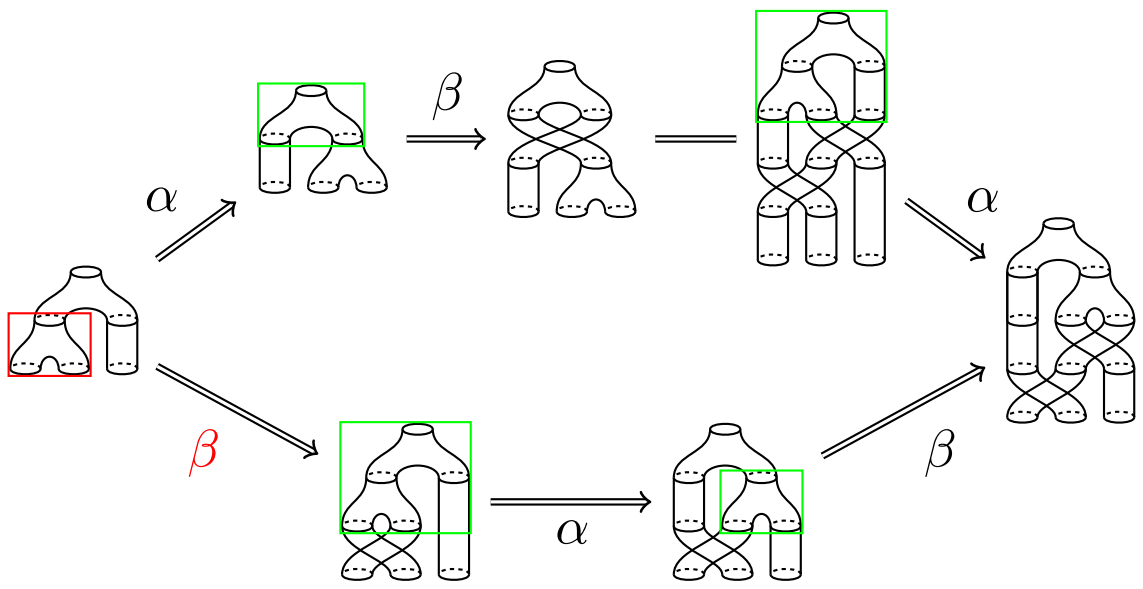

(8)

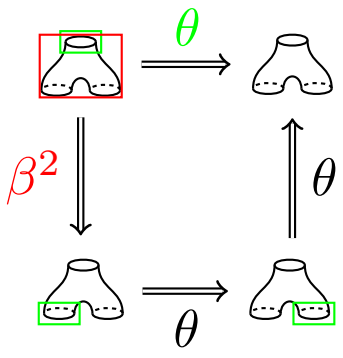

(9)

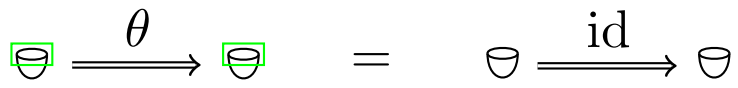

Note that the second hexagon axiom is redundant in the presence of a twist [23].

- (Rigidity) Write $\phi_l$ for the following composite ('left Frobeniusator'):

(10) $\phi_l \quad :=$

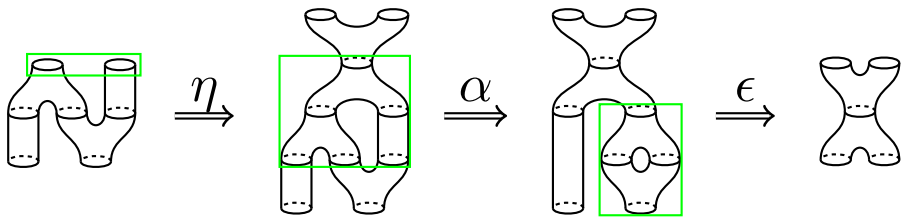

The left rigidity relation says that $\phi_l$ is invertible, with the following explicit inverse:

(11) $\phi_l^{-1} \quad =$

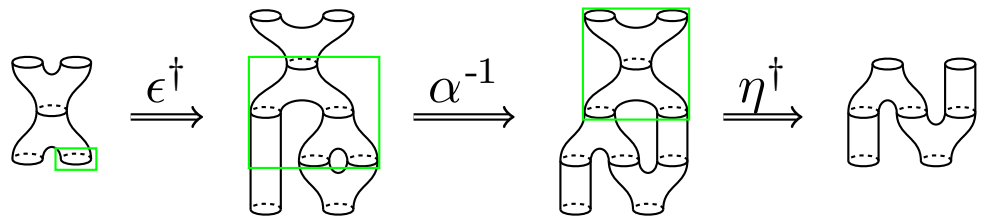

Similarly, write $\phi_r$ for $\phi_l$ rotated about the $z$-axis ('right Frobeniusator'):

(12) $\phi_r \quad :=$

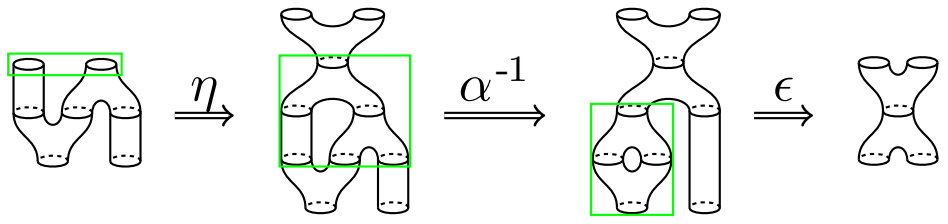

The right rigidity relation says that $\phi_r$ is invertible, with the following explicit inverse :

(13) $\phi_r^{-1} \quad =$

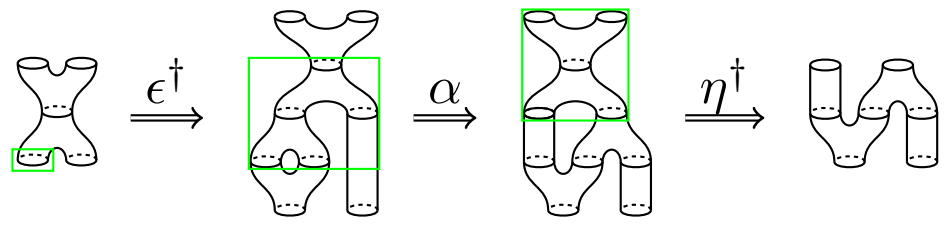

- (Ribbon) The twist satisfies the following equation:

(14)

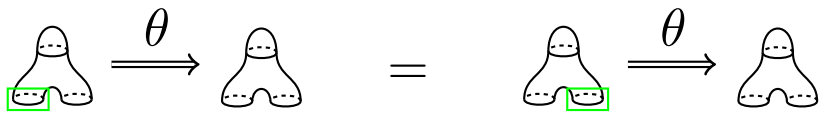

- (Biadjoint) The data (3) expresses 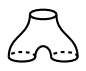 as the biadjoint of 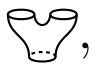, while (4) expresses 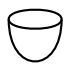 as the biadjoint of 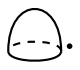. That is, the following equations hold, along with daggers and rotations about the $x$-axis:

(15)

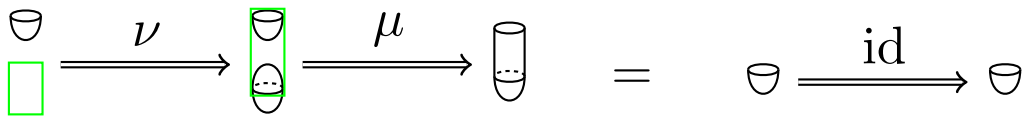

(16) 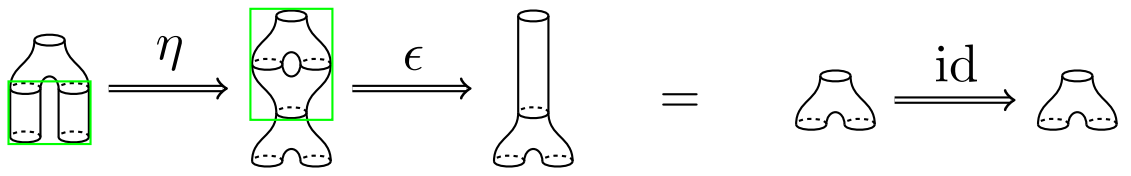

These are 8 equations in total.

- (Pivotality) The following equation holds, together with its rotation about the $z$-axis $(17)^z$:

(17) 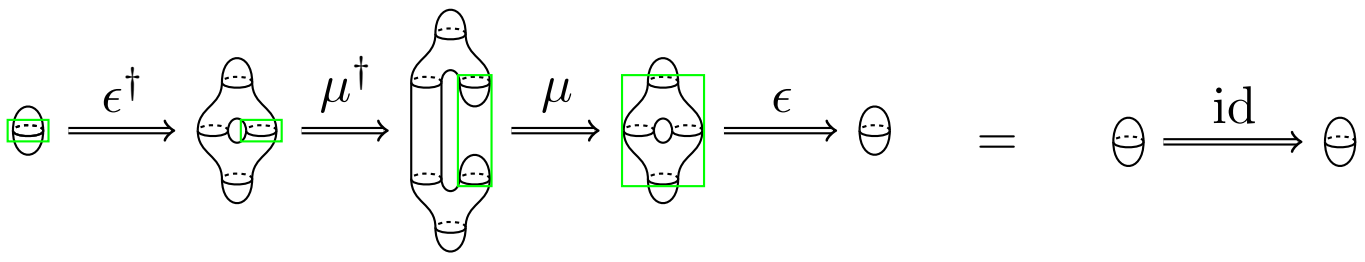

We will often refer to the Monoidal, Balanced, Rigidity, and Ribbon relations collectively as simply the Ribbon relations. The Pivotality axiom is so named because it endows the bicategory $\mathbf{F}(\mathcal{R})$ with a pivotal structure; see Proposition 9 below.

Our results in this paper are purely algebraic; the connection to geometry comes in the comparison with the bordism category in [7]. Nonetheless, it is useful to keep in mind the implicit geometric realizations of our bordisms, as follows:

- the generating object represents a circle;
- the generating 1-morphisms represent the 2-dimensional bordisms naively suggested by the pictured surfaces, namely the pants, copants, cup, and cap;
- the generating 2-morphisms $\alpha$, $\rho$, $\lambda$, and their inverses represent invertible 3-dimensional bordisms induced by the (boundary relative) ambient isotopies suggested by the pictured surfaces.
- the generating 2-morphisms $\eta^\dagger$, $\epsilon$, and $\mu^\dagger$ represent the bordism implementing the addition of a 2-handle about the curves

(18) 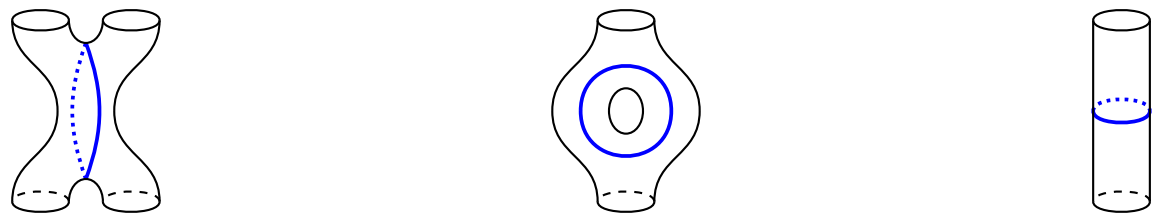

respectively; the 2-morphism $\nu^\dagger$ represents the bordism implementing the addition of a 3-handle;
- the generating 2-morphisms $\beta$ and $\theta$ represent the invertible 3-dimensional bordisms arising as mapping cylinders of the following diffeomorphisms:

(19) 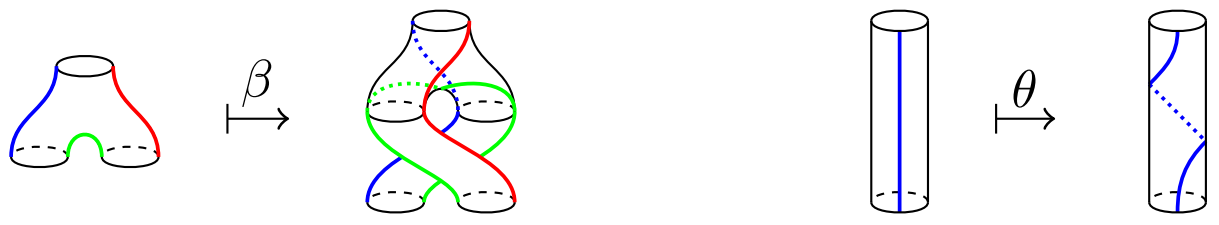

- for a 2-morphism $\sigma$ representing a particular bordism, the 2-morphism $\sigma^\dagger$ represents the time-reversed bordism.

**Definition 4.** The *modular presentation* $\mathcal{M}$ is defined as the ribbon presentation $\mathcal{R}$, plus one extra relation:

- (Modularity) The following equation holds, together with its rotation about the $z$-axis $(20)^z$:

(20)

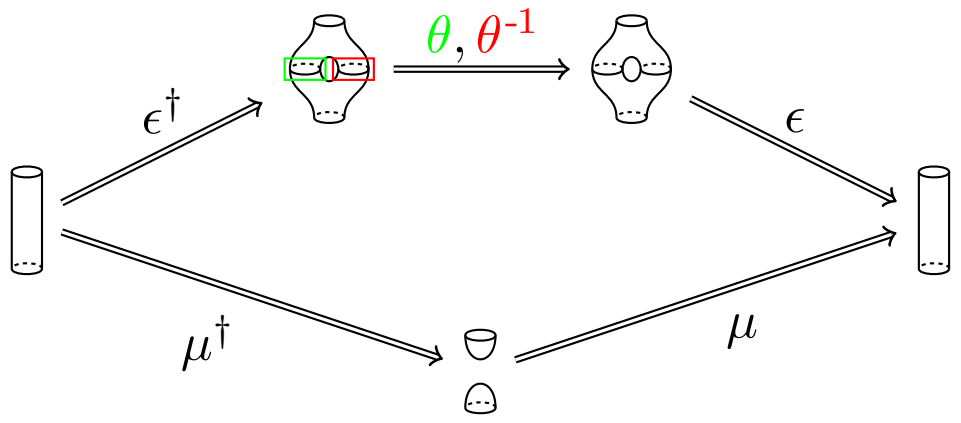

**Remark 5.** This Modularity relation corresponds to the $\mathcal{A}_2$ move from Kerler's bridged link calculus [27] for oriented 3-manifolds with boundary (see also [12, 38]). In [9] we show that a linear representation of $\mathbf{F}(\mathcal{M})$ corresponds roughly to a modular category. The Modularity relation corresponds, in the string diagram calculus for the modular category explained in that paper, to the identity that surrounding a strand by an untwisted unlabelled loop will cut the strand, cf [4, Corollary 3.1.10]. This identity is in turn equivalent to the requirement that the $s$-matrix is invertible, which is the usual notion of modularity. See Remark 27.

**Remark 6.** In [8] we show that there is a symmetric monoidal equivalence from the bicategory $\mathbf{F}(\mathcal{M})$ generated by this modular presentation to $\mathbf{Bord}_{123}^{\mathrm{csig}}$, the bordism bicategory where all $n$-dimensional manifolds come equipped with componentwise bounding $(n+1)$-dimensional manifolds, considered up to $(n+2)$-dimensional bordism, for $n = 1, 2, 3$. The bounding 4-manifolds are determined in this category by their signature, hence the notation $\mathbf{Bord}_{123}^{\mathrm{csig}}$.

**Definition 7.** The *anomaly-free modular presentation* $\mathcal{O}$ is the presentation with the same generators and relations as the modular presentation $\mathcal{M}$, plus one extra relation:

- (Anomaly-freeness) The following equation holds:

(21)

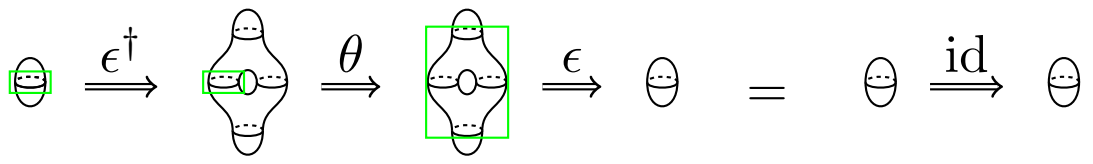

**Remark 8.** The Anomaly-freeness relation corresponds to the blow-up move $\mathcal{O}_1$ in the Kirby calculus for oriented 3-manifolds [12]. When interpreted in a modular category, it corresponds to the requirement that the two Gauss sums $p^+$ and $p^-$ are equal [9], which is the condition for a modular category to be anomaly-free [42].

## 3. Pivotality

In this section we work with the ribbon presentation $\mathcal{R}$, and prove the following result:

**Proposition 9.** *The 2-category $\mathbf{F}(\mathcal{R})$ has a canonical pivotal structure in which $\gamma^{**} = \gamma$ for every 2-morphism $\gamma$. The same holds for the 2-categories $\mathbf{F}(\mathcal{M})$ and $\mathbf{F}(\mathcal{O})$.*

The notion of a pivotal structure on a 2-category is recalled below. The essential relation needed to prove Proposition 9 is the Pivotality relation (17), whose significance can be explained as follows. In the full 2-Morse presentation $\mathcal{F}$, the generators 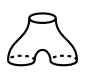 and 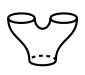 play a symmetric role; thus for instance there is a pentagon

relation (T1) as well as a co-pentagon relation $(\mathrm{T1})^x$. In the ribbon presentation $\mathcal{R}$, the generator 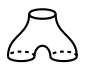 has been singled out as a preferred generator; thus for instance there is only a pentagon relation (5), and no co-pentagon relation. The Pivotality relation (17) is needed in order to restore the symmetry. Formally this is expressed using the language of pivotal structures.

Proposition 9 is used extensively in Section 6 since it allows us to prove various rotated forms of the relations in $\mathcal{F}$ by proving only a single form of the relation and then 'taking duals of both sides'.

3.1. **Duals of 2-morphisms in a bicategory.** Recall the notion of *left and right duals* of 2-morphisms in a bicategory[1]. Let $F, G : A \to B$ be 1-morphisms in a bicategory, and suppose that $F^*, G^* : B \to A$ are right adjoints of $F$ and $G$ respectively, with unit and counit 2-morphisms

$$\eta_F \colon \mathrm{id}_A \Rightarrow F^* \circ F, \quad \epsilon_F \colon F \circ F^* \Rightarrow \mathrm{id}_B$$
$$\eta_G \colon \mathrm{id}_A \Rightarrow G^* \circ G, \quad \epsilon_G \colon G \circ G^* \Rightarrow \mathrm{id}_B\,.$$

Then if $\gamma : F \Rightarrow G$ is a 2-morphism, we can form its *right dual* $\gamma^* : G^* \Rightarrow F^*$ as

$$\gamma^* = (\mathrm{id}_{F^*} * \epsilon_G) \circ (\mathrm{id}_{F^*} * \gamma * \mathrm{id}_{G^*}) \circ (\eta_F * \mathrm{id}_{G^*})$$

where we have suppressed the associator and unitors, and we are using the notation and conventions of Leinster [34] on vertical ($\circ$) and horizontal ($*$) composition in a bicategory. In the string diagram notation for bicategories (see [33] or [6, Chapter 4]), we can express $\gamma^*$ as

(22)

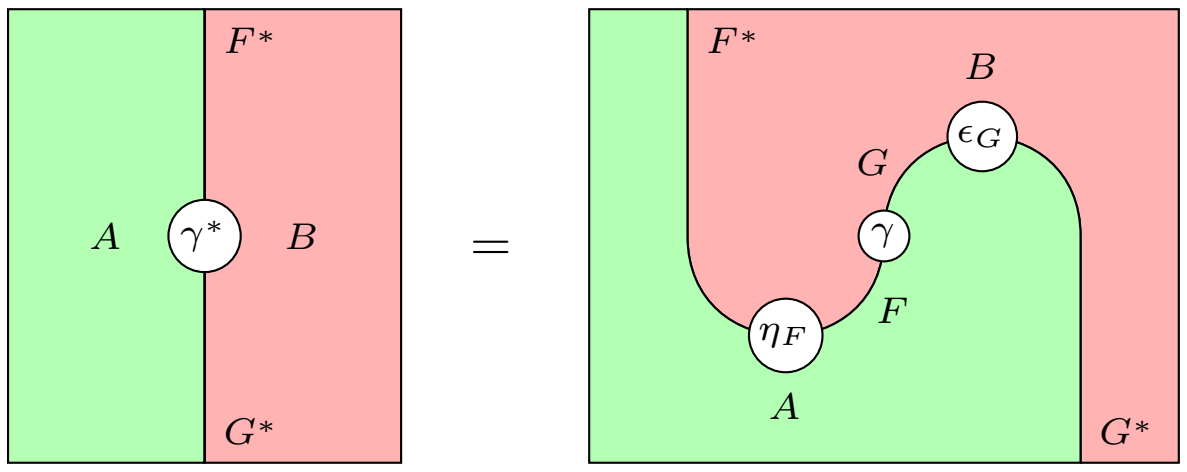

where composition of 1-morphisms runs from right to left and vertical composition of 2-morphisms runs from bottom to top.

Taking right duals of 2-morphisms in a bicategory interacts with vertical composition ($\circ$), horizontal composition ($*$), and tensor product ($\boxtimes$, if present) of 2-morphisms as follows:

- $\mathrm{id}^* = \mathrm{id}$
- $(\gamma \circ \psi)^* = \psi^* \circ \gamma^*$
- $(\gamma * \psi)^* = \psi^* * \gamma^*$
- $(\gamma \boxtimes \psi)^* = \gamma^* \boxtimes \psi^*$

In particular, if $\gamma$ is invertible, then $\gamma^*$ is invertible, with $(\gamma^*)^{-1} = (\gamma^{-1})^*$.

Besides right duals, there are also left duals. In the situation above, suppose that ${}^*F, {}^*G : B \to A$ are *left* adjoints of $F$ and $G$ respectively. Then the *left dual* ${}^*\gamma : {}^*G \Rightarrow {}^*F$ is defined using the unit and counit 2-morphisms

$$n_F \colon \mathrm{id}_A \Rightarrow F \circ {}^*F, \quad e_F \colon {}^*F \circ F \Rightarrow \mathrm{id}_B$$
$$e_G \colon \mathrm{id}_A \Rightarrow G \circ {}^*G, \quad e_G \colon {}^*G \circ G \Rightarrow \mathrm{id}_B\,.$$

[1]This notion is sometimes referred to as the *calculus of mates* [25].

which witness $^*F$ and $^*G$ as *left* adjoints of $F$ and $G$ respectively:

(23) 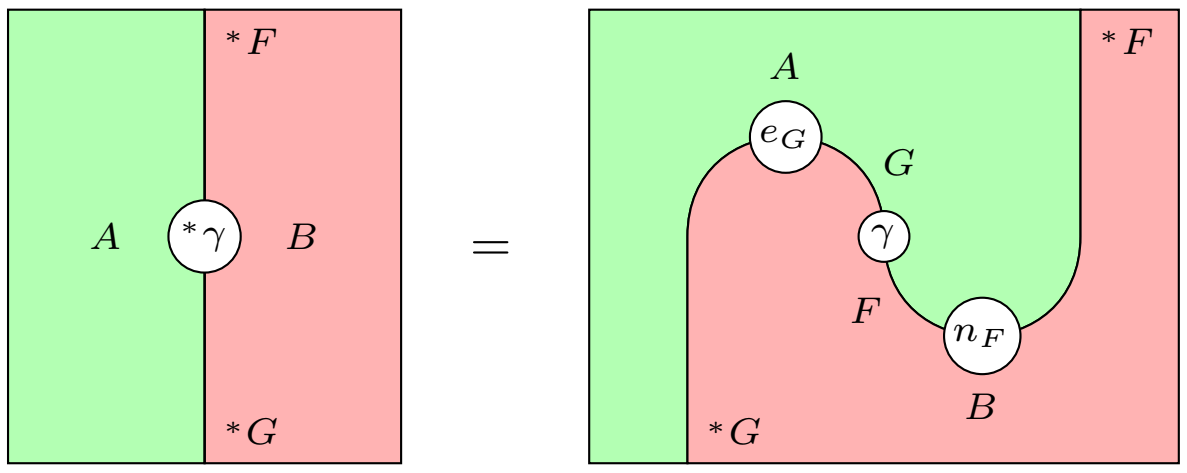

The same properties as above hold for left duals, and we always have

- $^*(\gamma^*) = (^*\gamma)^* = \gamma$.

However, the equation $\gamma^* = {}^*\gamma$ (which is equivalent to $\gamma^{**} = \gamma$) need not hold; indeed in general they have different source and target 1-morphisms.

In general, a bicategory $\mathbf{C}$ is said to have *left and right adjoints* if every 1-morphism $F\colon A \to B$ has a left adjoint $F^*$ and a right adjoint $^*F$. If we make a choice of these adjoints (together with their unit and counit maps) for each $F$, then we obtain functors

$$(-)^*,\ {}^*(-) : \mathbf{C} \to \mathbf{C}^{\mathrm{op}}$$

where $\mathbf{C}^{\mathrm{op}}$ is the bicategory with the same objects as $\mathbf{C}$ but whose 1-morphisms and 2-morphisms have reversed their sources and targets.

**Definition 10.** Let $\mathbf{C}$ be a bicategory with chosen left and right adjoints for each 1-morphism. A *pivotal structure* on $\mathbf{C}$ is a choice of natural isomorphism of functors $(-)^* \cong {}^*(-)$.

We remark that a pivotal structure can be equivalently reformulated as an *even-handed structure* in the sense of [6].

3.2. **Duals of 2-morphisms in $\mathbf{F}(\mathcal{R})$.** Let us apply these definitions in the context of the bicategory $\mathbf{F}(\mathcal{R})$. First we need to choose left and right adjoints for each 1-morphism.

The Biadjoint relations (16) witness 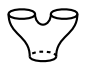 as a right adjoint of 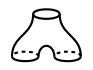 via the unit and counit maps $\eta$ and $\epsilon$, and as a left adjoint via $\epsilon^\dagger$ and $\eta^\dagger$. Similarly the unit and counit maps $\nu$ and $\mu$ witness 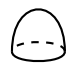 as a right adjoint of 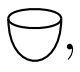, and as a left adjoint via $\mu^\dagger$ and $\nu^\dagger$.

The right adjoint 1-morphism of a composite of 1-morphisms in $\mathbf{F}(\mathcal{R})$ is constructed inductively according to the rules $(S \circ T)^* = T^* \circ S^*$ and $(S \boxtimes T)^* = S^* \boxtimes T^*$, and similarly for the left adjoint 1-morphisms. The accompanying unit and counit maps are also defined inductively in the same way, by composing those of the generators. Finally, the right and left duals of 2-morphisms are defined by the equations (22) and (23) respectively. Thus, we have defined functors

$$(-)^*,\ {}^*(-) : \mathbf{F}(\mathcal{R}) \to \mathbf{F}(\mathcal{R})^{\mathrm{op}}\,.$$

For example, according to this definition the right dual of

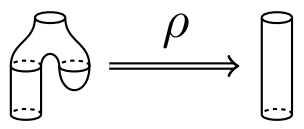

is

$$\rho^* \quad = $$

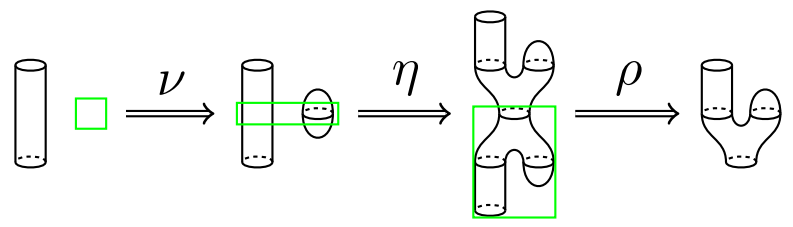

while the left dual is

$$^*\rho \quad = $$

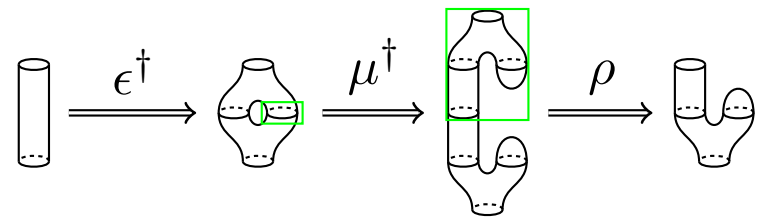

By definition, a pivotal structure on $\mathbf{F}(\mathcal{R})$ is a natural isomorphism $(-)^* \cong {}^*(-)$. In fact, the pivotal structure we will define will be the *identity* natural isomorphism. Thus, to show that this is a pivotal structure amounts to showing that $\gamma^* = {}^*\gamma$ for every 2-morphism in $\mathbf{F}(\mathcal{R})$ (i.e. we must check naturality with respect to $\gamma$). It is sufficient to show this on the generating 2-morphisms, and, moreoever, since every equation between 2-morphisms which holds in $\mathbf{F}(\mathcal{R})$ also holds in $\mathbf{F}(\mathcal{M})$ and $\mathbf{F}(\mathcal{O})$, this will establish Proposition 9 for those bicategories too.

We record the following for later use.

**Lemma 11.** *In $\mathbf{F}(\mathcal{R})$, the following equations hold:*

(i) $\theta^* = \theta$

(ii) 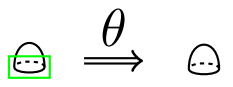 $=$ id

*Proof.* Item (i) follows since the identity 1-morphism 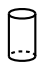 trivially has itself as a right dual. Item (ii) follows from taking right duals of both sides of the Ribbon relation (9). □

3.3. **Rotations of the ribbon generators.** It is useful to have a shorthand for the following composites of generators of $\mathcal{R}$.

**Definition 12.** We define the following composite 2-morphisms in $\mathbf{F}(\mathcal{R})$:

(24) $\quad \check{\alpha} := (\alpha^{-1})^* =$ 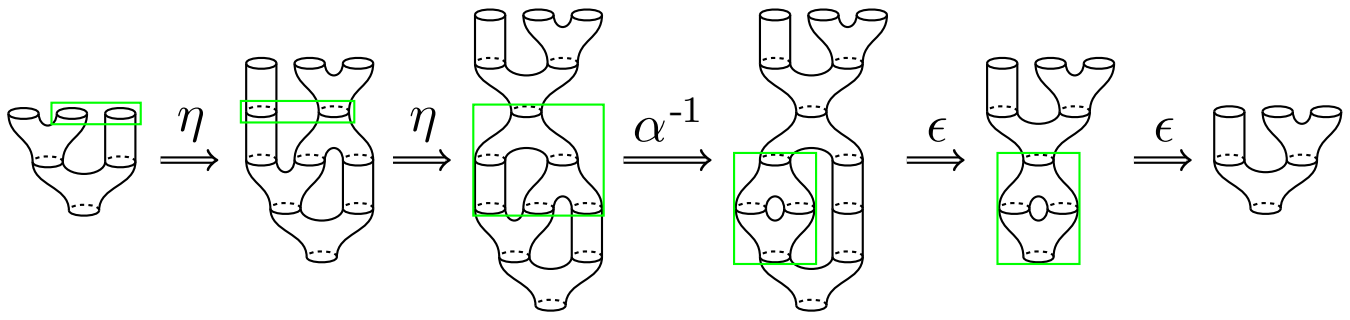

(25) $\quad \check{\rho} := (\rho^{-1})^* =$ 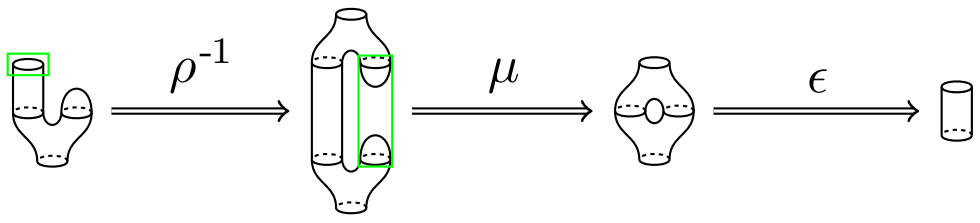

(26) $\quad \check{\lambda} := (\lambda^{-1})^* =$ 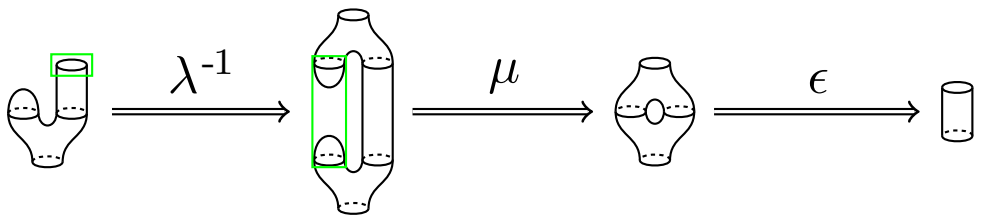

(27) $\check{\beta} := \beta^* =$ 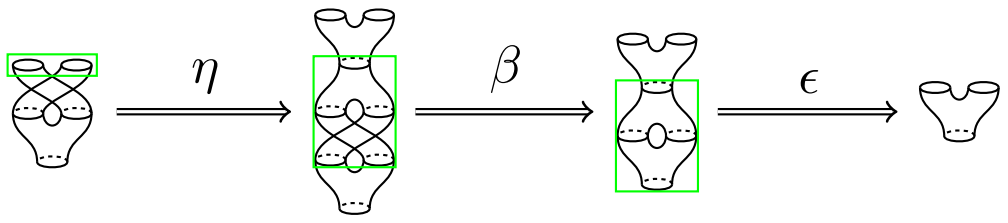

**Remark 13.** It follows from Section 3.1 that the inverses of these 2-morphisms are as follows:

$$\check{\alpha}^{-1} = \alpha^* \qquad \check{\rho}^{-1} = \rho^* \qquad \check{\lambda}^{-1} = \lambda^* \qquad \check{\beta}^{-1} = \beta^* \tag{28}$$

3.4. **Pivotality on the cylinder.** We will use the following lemma many times in our calculations.

**Lemma 14** (Local nature of creating a hole)**.** *In* $\mathbf{F}(\mathcal{R})$*, the following equation holds for* $\epsilon^\dagger$*, as well as its rotated version* $(29)^z$*:*

(29) $\epsilon^\dagger \quad =$ 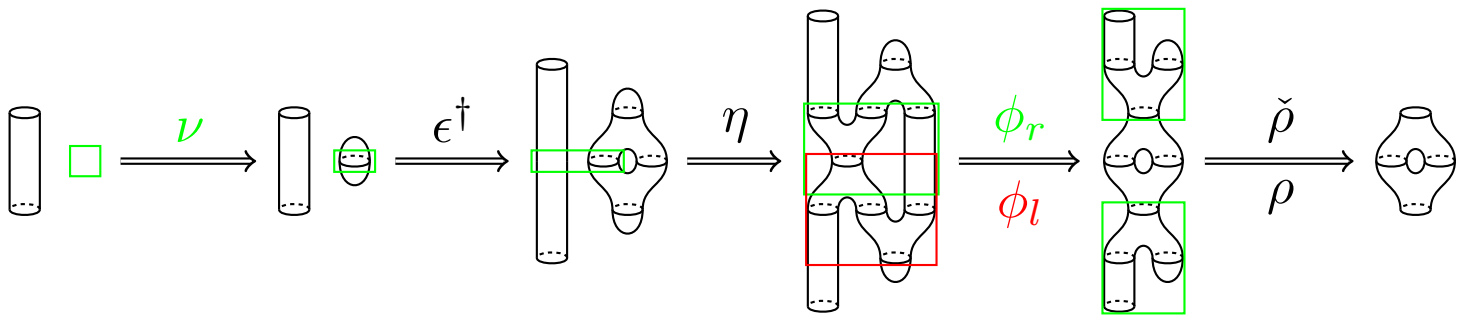

*Proof.* Substitute in the definition of $\check{\rho}$ as $(\rho^{-1})^*$, move the $\rho^{-1}$ to the beginning of the composite and move the $\mu$ forward to cancel with the $\nu$, yielding the upper path below:

(30) 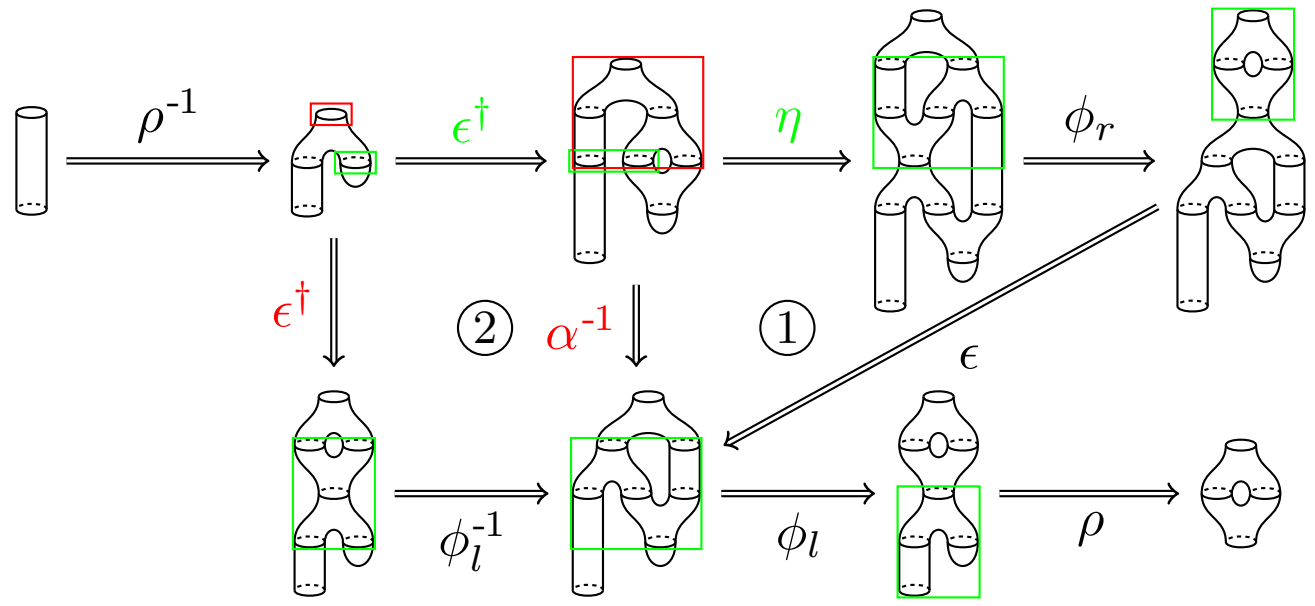

Diagram ① follows from inserting the definition of $\phi_r$, and using the $\eta$–$\epsilon$ adjunction equation twice. Diagram ② commutes by inserting the definition of $\phi_l^{-1}$, and using the $\epsilon^\dagger$–$\eta^\dagger$ adjunction equation once. The counterclockwise composite equals $\epsilon^\dagger$, by cancelling the $\phi_l^{-1}$ with the $\phi_l$ and cancelling the $\rho^{-1}$ with the $\rho$. The relation $(30)^z$ can be proved in a similar way. □

The following lemma shows that the Pivotality relation (17), which was initially defined on the sphere, also holds on the cylinder.

**Lemma 15.** *In* $\mathbf{F}(\mathcal{R})$*, the following composite is the identity, together with its rotated version* $(31)^z$*:*

(31) 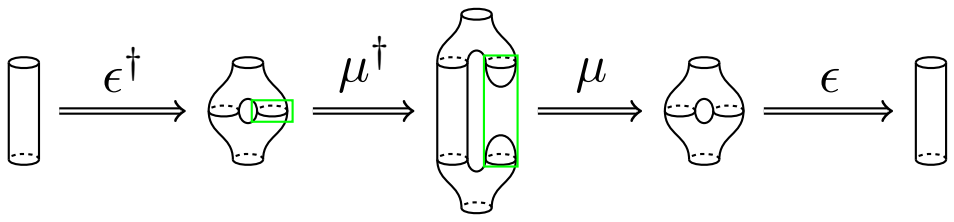

*Proof.* Replace $\epsilon^\dagger$ using Lemma 14, and move the $\mu \circ \mu^\dagger$ terms and the $\epsilon$ term forward. We obtain the clockwise composite below:

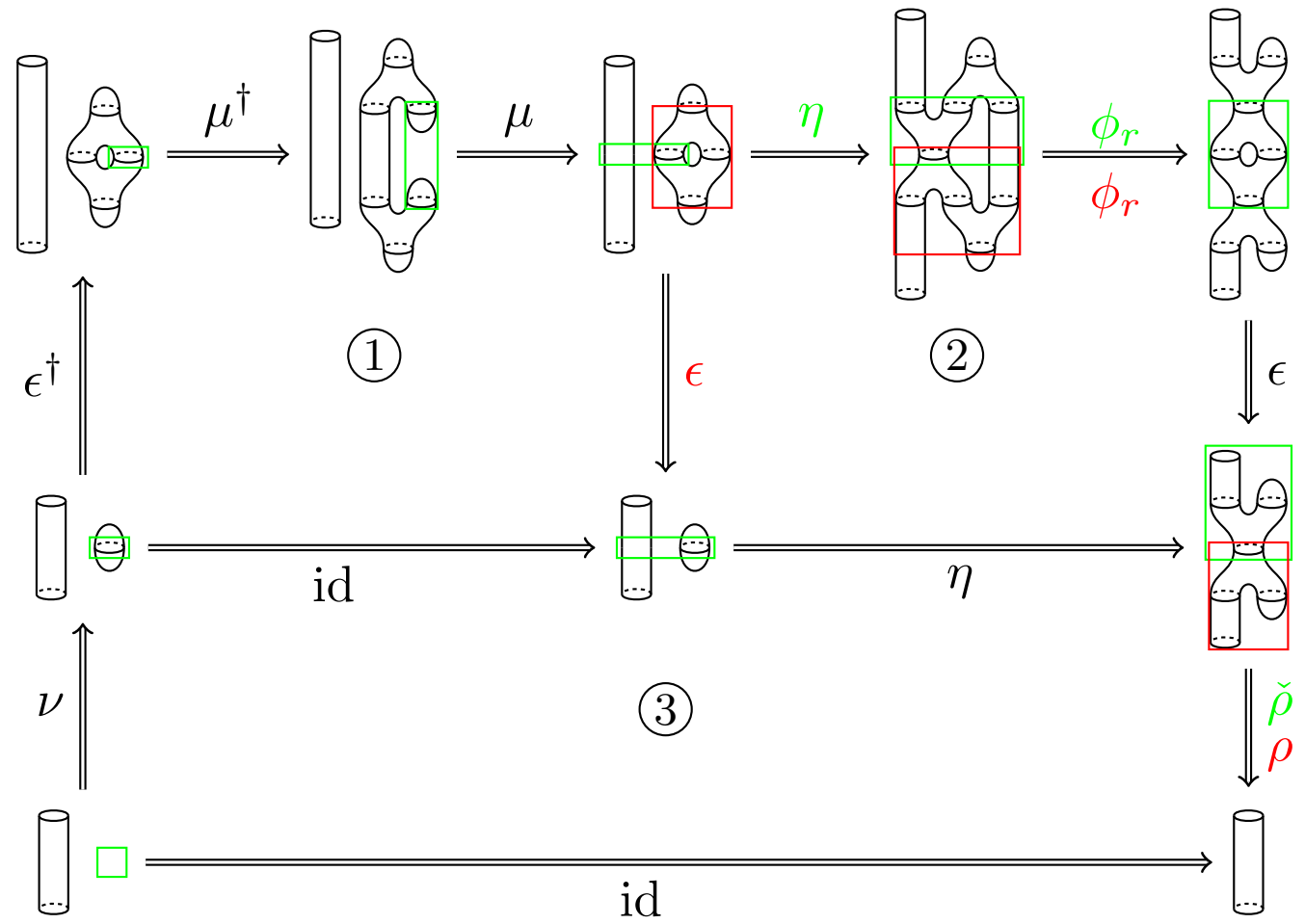

Diagram ① is the Pivotality relation. Diagram ② follows from expanding the Frobeniusators and using the $\eta$–$\epsilon$ adjunction equations. Diagram ③ follows from expanding the definition of $\check{\rho}$ as $(\rho^{-1})^*$ and then using the $\nu$–$\mu$ and $\eta$–$\epsilon$ adjunction equations. The relation $(31)^z$ is proved in a similar way. □

3.5. **Proof of pivotality.** As discussed in Section 3.2, to prove Proposition 9 it is sufficient to show that $\gamma^{**} = \gamma$ for each generating 2-morphism $\gamma$ of $\mathcal{R}$, which we accomplish in a series of lemmas. To start with, note that if $\gamma$ is one of the noninvertible generators, then $\gamma^{**} = \gamma$ follows from the adjunction relations, since the adjunction relations ensure that both $\gamma^*$ and ${}^*\gamma$ are $\gamma^\dagger$.

**Lemma 16.** *In* $\mathbf{F}(\mathcal{R})$, $\rho^{**} = \rho$ *and* $\lambda^{**} = \lambda$.

*Proof.* We have $[\rho^{**} = \rho] \Leftrightarrow [\rho^* = {}^*\rho] \Leftrightarrow [(\rho^{-1})^* \circ {}^*\rho = \mathrm{id}]$. Let us expand the last equation here:

$(\rho^{-1})^* \circ {}^*\rho \quad =$

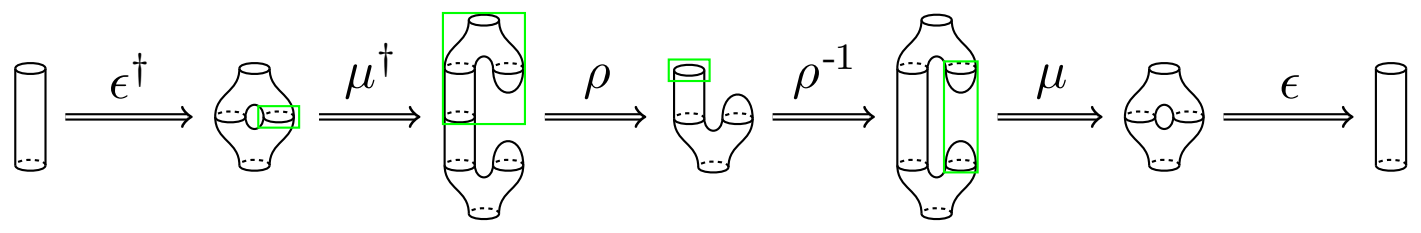

$=$

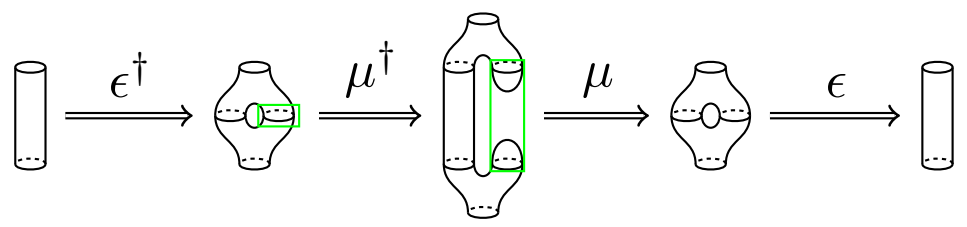

$=$ id (by Lemma 15)

The proof that $\lambda^{**} = \lambda$ is similar. □

**Lemma 17.** *In* $\mathbf{F}(\mathcal{R})$, *the following relation holds for* $\eta^\dagger$, *together with its rotated version* $(32)^z$:

(32) $\eta^\dagger$ =

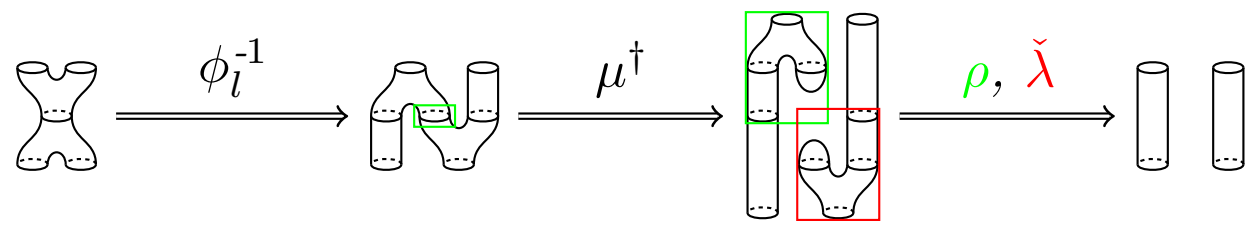

*Proof.* In the counterclockwise composite, expand out $\phi_l^{-1}$ and rearrange using naturality to arrive at the upper clockwise composite below:

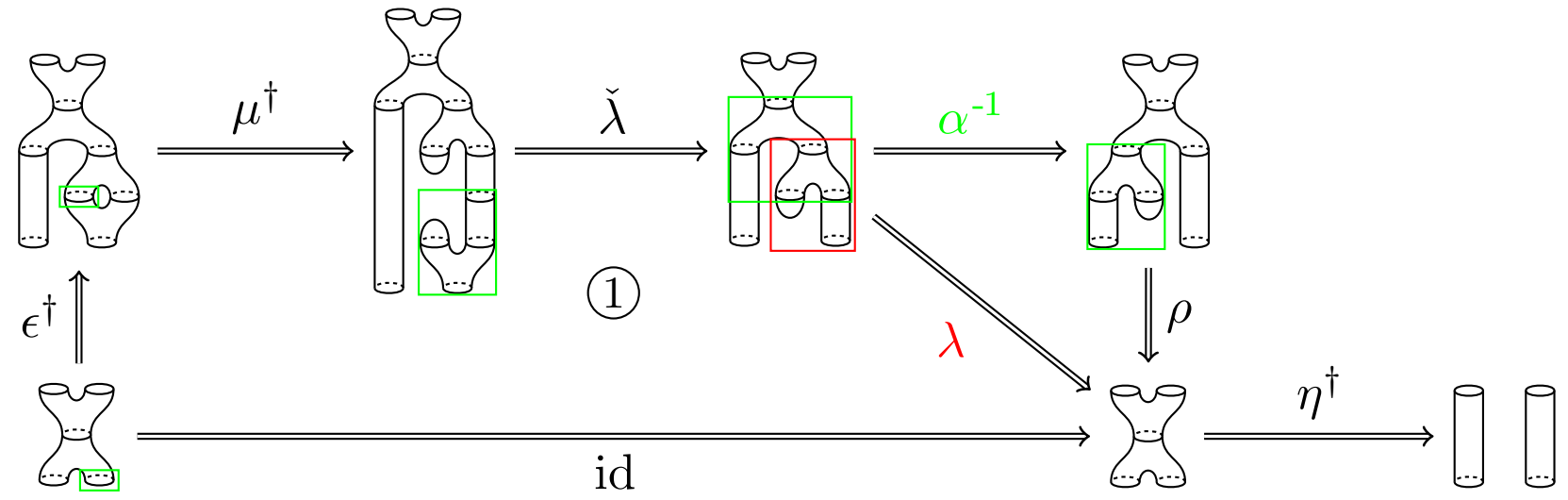

Diagram ① follows from expanding $\check{\lambda}$ as ${}^*(\lambda^{-1})$ using Lemma 16, and then applying the adjunction equations. The triangle diagram is the monoidal unit diagram (6). □

**Lemma 18.** *In* $\mathbf{F}(\mathcal{R})$, *the equations* $\phi_r = (\phi_l^{-1})^*$ *and* $\phi_l = (\phi_r^{-1})^*$ *hold, and thus* $\phi_l^{**} = \phi_l$ *and* $\phi_r^{**} = \phi_r$.

*Proof.* To show that $\phi_r = (\phi_l^{-1})^*$, unravel the definition of $(\phi_l^{-1})^*$, to obtain the clockwise composite below:

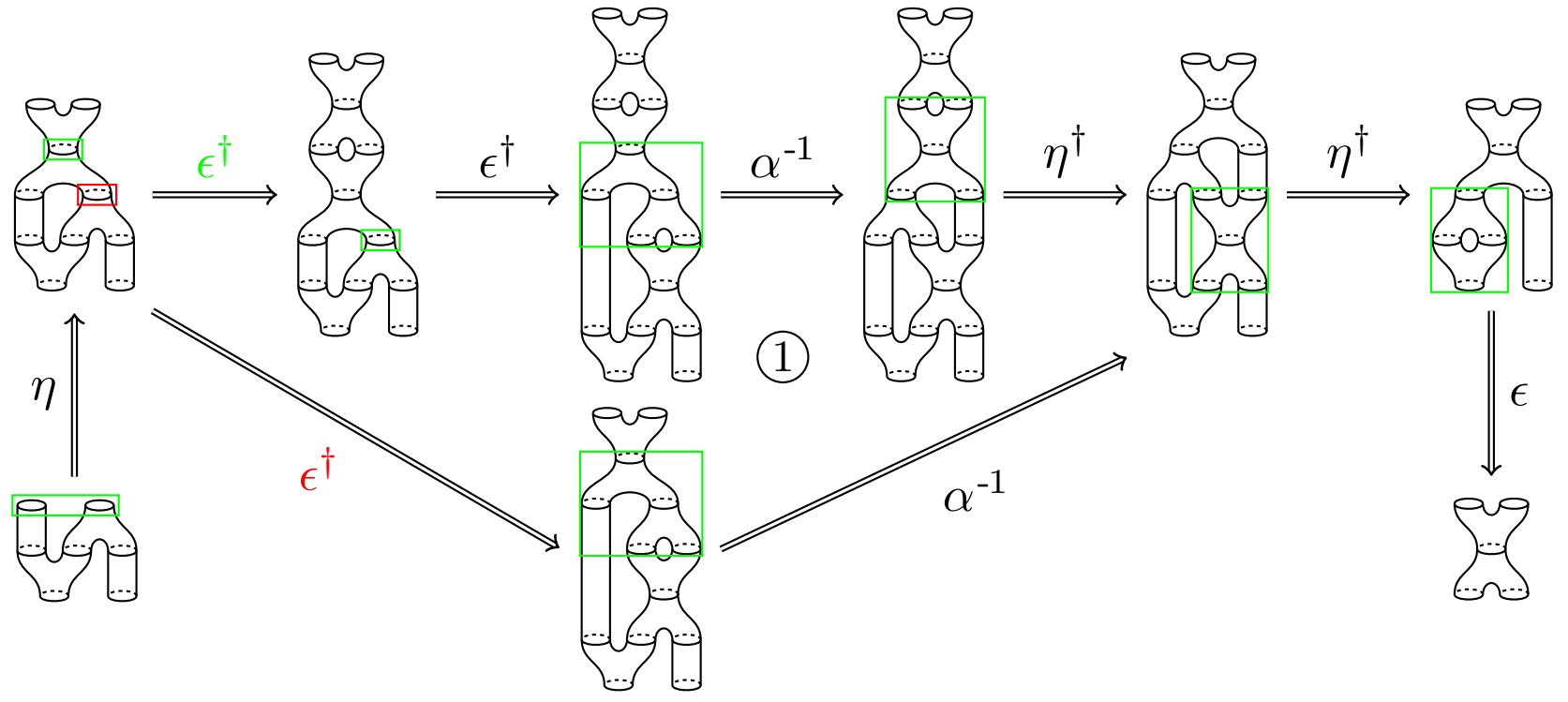

Diagram ① follows from commuting the first $\epsilon^\dagger$ to just before the $\eta^\dagger$ and then using the $\epsilon^\dagger$–$\eta^\dagger$ adjunction equation. The result now follows from commuting the $\eta^\dagger$ to the left of the $\alpha^{-1}$ to cancel with the $\epsilon^\dagger$. What remains is the definition of $\phi_r$. The equation $\phi_l = (\phi_r^{-1})^*$ follows similarly. Taking inverses of both sides of each of these equations, we obtain that $\phi_r^{-1} = \phi_l^*$ and $\phi_l^{-1} = \phi_r^*$. Hence $\phi_l^{**} = (\phi_r^{-1})^* = (\phi_r^*)^{-1} = \phi_l$ and similarly $\phi_r^{**} = \phi_r$. □

**Lemma 19.** *In* $\mathbf{F}(\mathcal{R})$, $\alpha^{**} = \alpha$.

*Proof.* Unravelling the definition of $\alpha^{**}$, we obtain the following:

$\alpha^{**}$ =

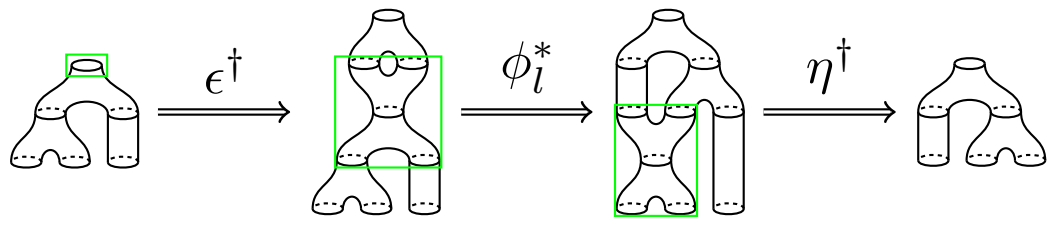

Using Lemma 18, we replace $\phi_l^*$ with $\phi_r^{-1}$. The result then follows from expanding out $\phi_r^{-1}$ and using the $\epsilon^\dagger$–$\eta^\dagger$ adjunction equation twice. □

**Lemma 20.** *In* $\mathbf{F}(\mathcal{R})$*,* $\beta^{**} = \beta$*.*

*Proof.* This is equivalent to $(\beta^{-1})^{**} \circ \beta = \text{id}$. Start with $(\beta^{-1})^{**} \circ \beta$, and move the $\epsilon^\dagger$ inside the $(\beta^{-1})^{**}$ to the front. We arrive at the composite on the top row below:

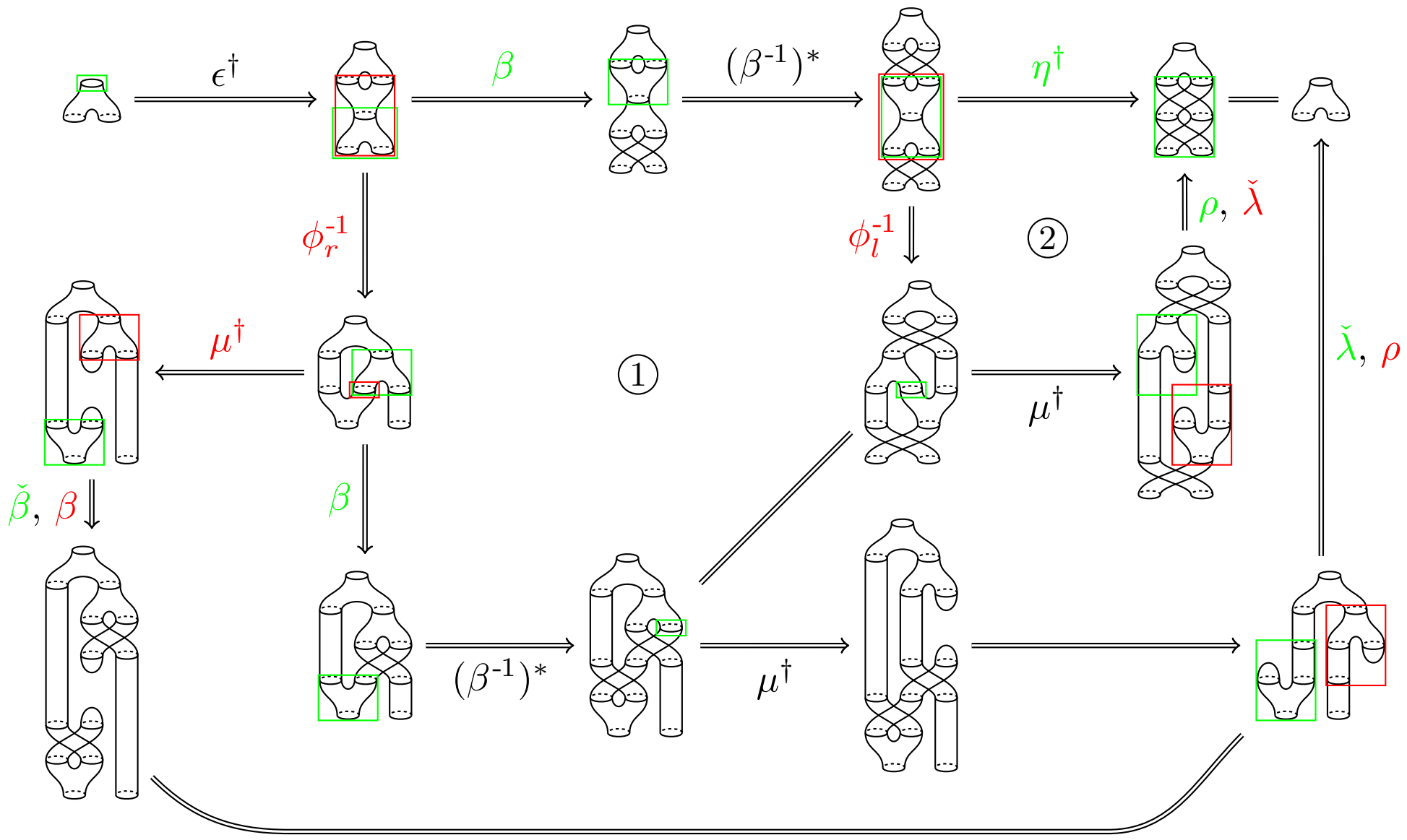

Diagram ① is (FFG1) in Appendix A, which is proved below in Section 6.9 using only the expansion of the Frobeniusators, the adjunction equations, and coherence of the hexagon relation, and hence we can apply it here. Diagram ② is Lemma 17. In the outer counterclockwise composite, the $\beta, \check{\beta}$ terms fall away by coherence since braiding with the unit has no effect. The following composite remains:

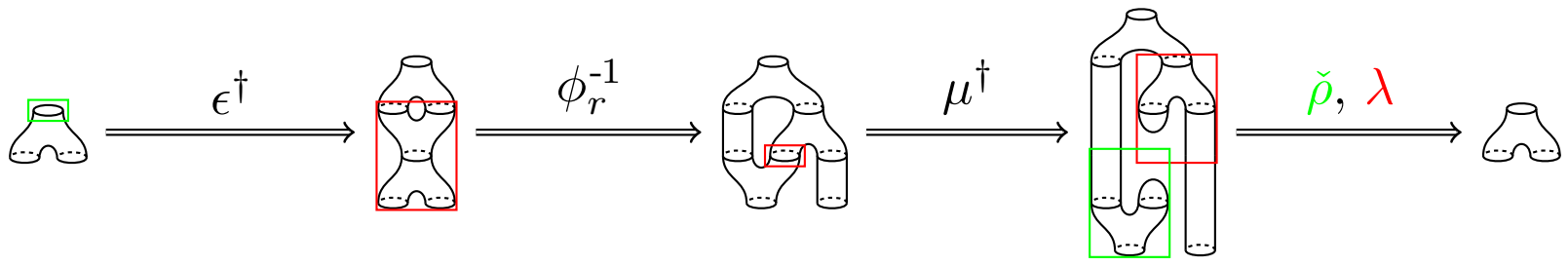

We apply Lemma 17 again to transform the last three arrows into an $\eta^\dagger$, which cancels with the $\epsilon^\dagger$ and we are done. □

## 4. Modularity

In this section we add the Modularity relation (20) to the relations considered thus far; that is we derive certain equations between 2-morphisms which hold in $\mathbf{F}(\mathcal{M})$. Note that since every relation in the modular presentation $\mathcal{M}$ is also a relation in the anomaly-free presentation $\mathcal{O}$, each of these equations will also hold in $\mathbf{F}(\mathcal{O})$.

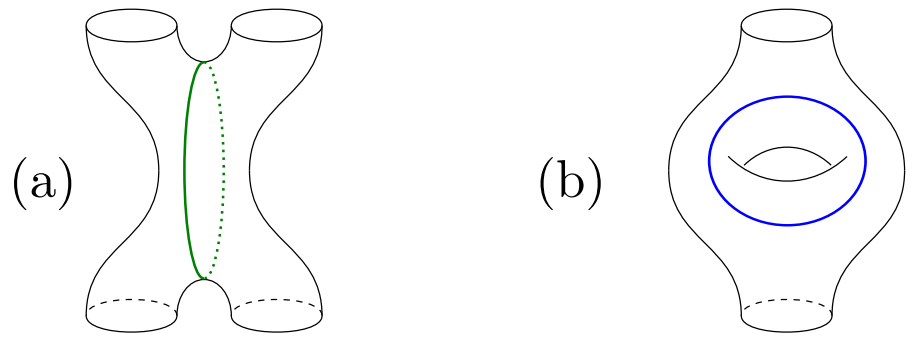

FIGURE 1. The 2-morphisms II and $A$ represent right-handed Dehn twists around the simple closed curves (a) and (b) respectively.

The first use of the Modularity relation will be to establish representatives in $\mathbf{F}(\mathcal{M})$ for the Dehn twists about the closed curves (a) and (b) in Figure 1. The mapping cylinders representing these Dehn twists are by definition present in the bordism bicategory $\mathbf{Bord}^{\mathrm{or}}_{123}$, but we have not yet seen explicit representatives for them in $\mathbf{F}(\mathcal{M})$.

In Section 4.1 we define 2-morphisms II and $A$ as certain composites of generating 2-morphisms in $\mathbf{F}(\mathcal{M})$; these 2-morphisms represent the 'missing Dehn twists'. We similarly define a composite 2-morphism III representing a diffeomorphism of the twice-punctured torus that switches the meridian and longitude. As the generating 2-morphisms used in these composites are not all invertible, we also need to prove that the 2-morphisms II, $A$, and III are nevertheless still invertible in $\mathbf{F}(\mathcal{M})$. In Section 4.2 we show that these newly defined 2-morphisms are self-dual in the sense that $\mathrm{II}^* = \mathrm{II}$, $A^* = A$ and $\mathrm{III}^* = \mathrm{III}$. In Section 4.3 we derive an alternative decomposition of III and use it to derive some known mapping class group identities as consequences of the relations in $\mathbf{F}(\mathcal{M})$.

### 4.1. **The missing diffeomorphisms.**

**Definition 21.** We define the following composites in $\mathbf{F}(\mathcal{M})$[2]:

(33) II :=
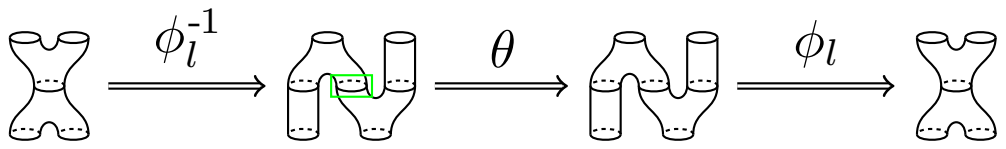

(34) $A$ :=
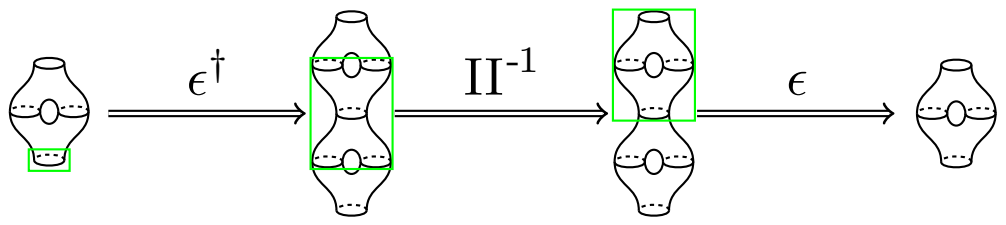

(35) III :=
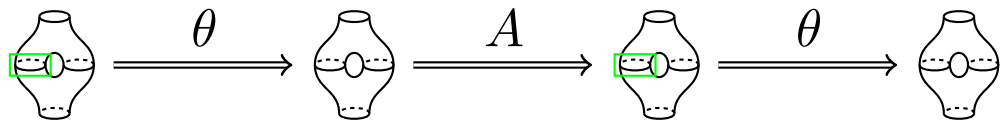

(36) III′ :=
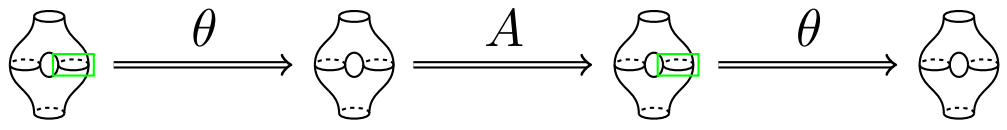

[2] The notation here is chosen to evoke the Hatcher-Thurston presentation of the mapping class group [20].

These composites are intended to represent mapping cylinders of diffeomorphisms, so we need to show that they are invertible. The Rigidity relation implies that II is invertible. It is not *a priori* clear that $A$, III and III$'$ are invertible, because they are constructed from noninvertible generators. We now show that the invertibility of $A$ follows from the Modularity relation, and is in fact equivalent to it. We will need the following two lemmas.

**Lemma 22.** *In* $\mathbf{F}(\mathcal{M})$*, the following identity holds:*

$A \quad =$ 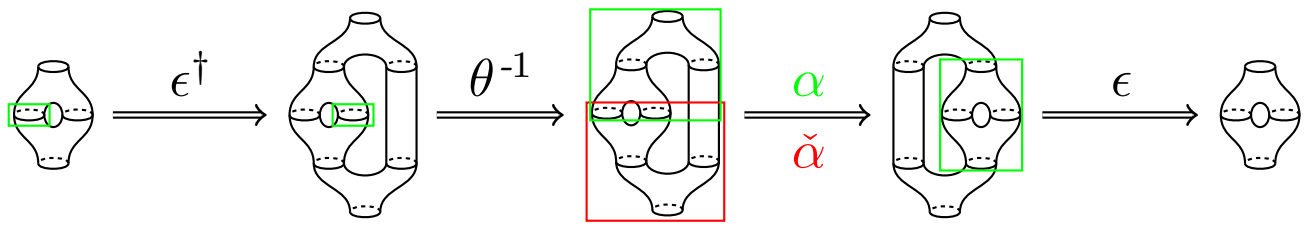

*Proof.* Start with the definition of $A$ from (34) in the clockwise composite below:

(37) 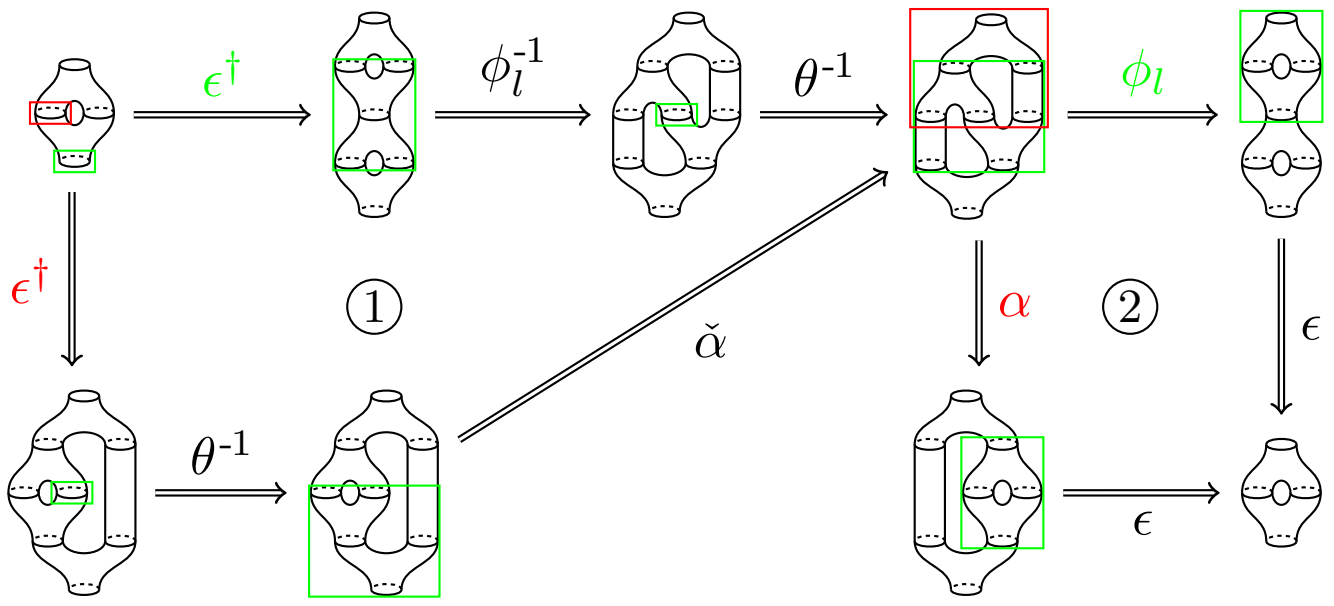

Diagram ① commutes by expanding $\breve{\alpha}$ as ${}^*(\alpha^{-1})$ using Lemma 19, and then using the $\eta^\dagger$–$\epsilon^\dagger$ adjunction equations. Diagram ② commutes by expanding $\phi_l$ and then using the $\eta$–$\epsilon$ adjunction equations. ☐

**Lemma 23.** *In* $\mathbf{F}(\mathcal{M})$*, the following diagram commutes, together with* $(38)^\dagger$, $(38)^z$ *and* $(38)^{z\dagger}$*:*

(38) 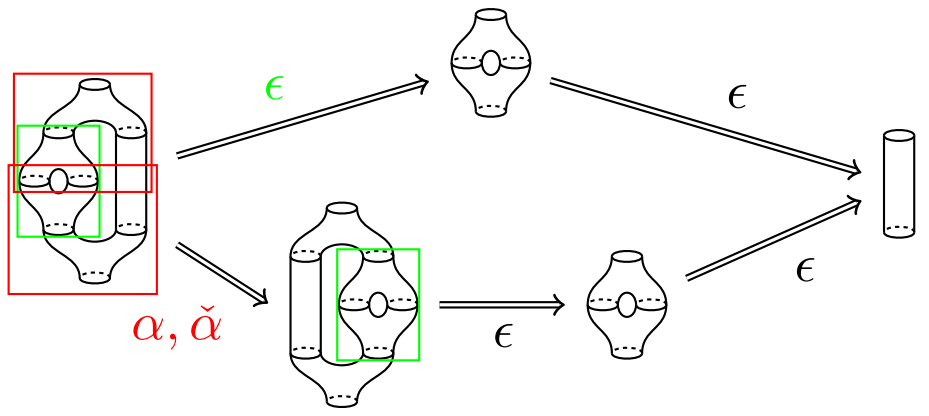

*Proof.* This result follows from the expansion of $\breve{\alpha} := (\alpha^{-1})^*$ in terms of $\eta$, $\alpha^{-1}$, and $\epsilon$, and the $\eta$–$\epsilon$ adjunction equations. The other versions follow similarly. ☐

The *once-punctured torus* and *twice-punctured torus* refer to the following composite 1-morphisms in $\mathbf{F}(\mathcal{M})$:

once-punctured torus := 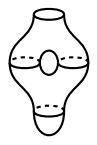 twice-punctured torus := 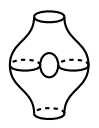

**Remark 24.** On the once-punctured torus III = III$'$. This follows from the Ribbon relation (14).

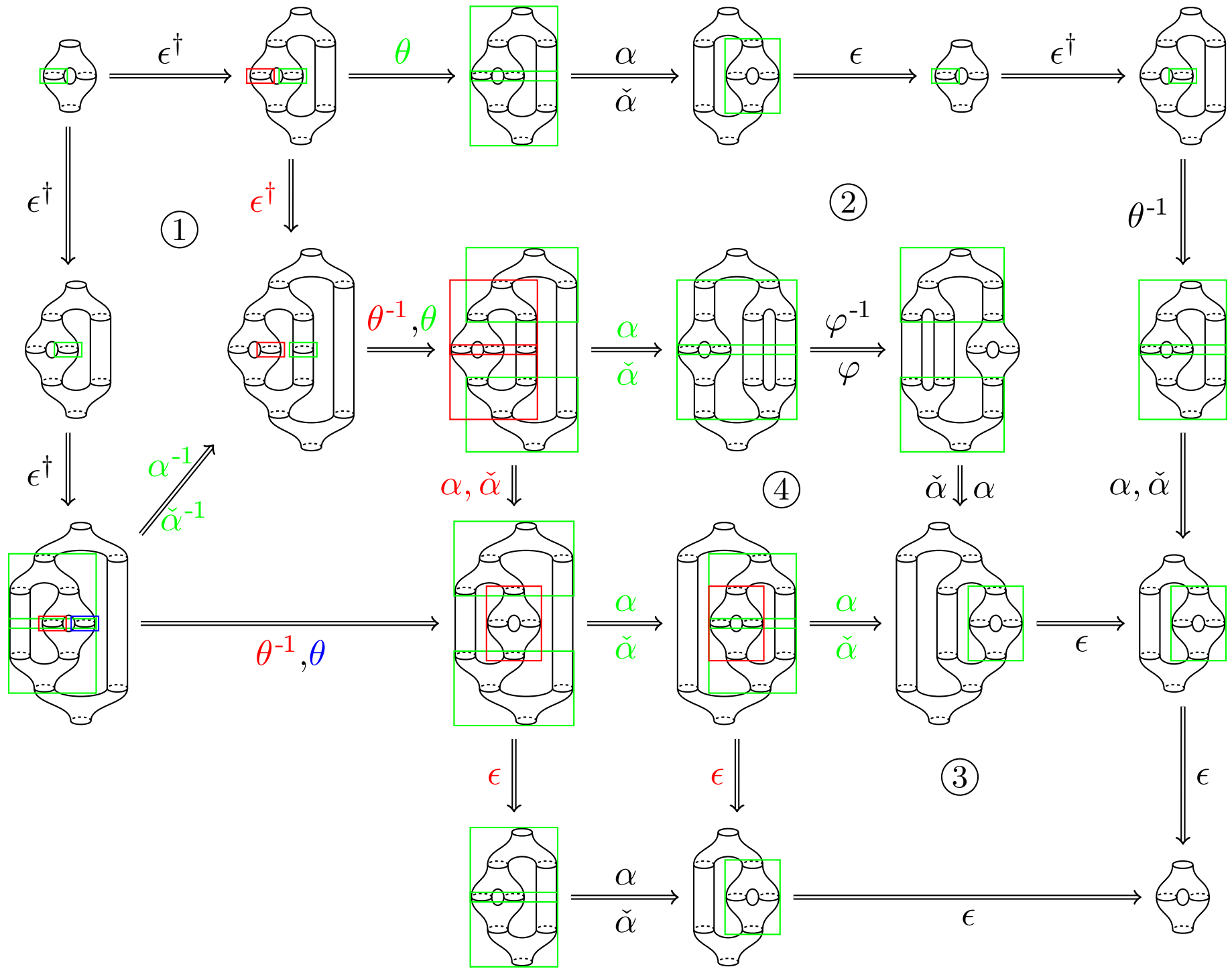

FIGURE 2. Simplifying $A \circ A^\dagger$.

We now show that the invertibility of $A$ is equivalent to the Modularity relation.

**Proposition 25.** *In the presence of the Ribbon and Biadjoint relations, the following are equivalent:*

*(i) The Modularity relation* (20)
*(ii) $A$ is invertible on the twice-punctured torus, with inverse $A^\dagger$*
*(iii) $A$ is invertible on the once-punctured torus, with inverse $A^\dagger$*

*Proof.* (i) ⇒ (ii). Start with $A \circ A^\dagger$ as in the clockwise composite in Figure 2, where we have used Lemma 22, and simplify as shown. Diagrams ① and ③ are Lemma 23, ② is naturality and ④ is the pentagon relation. In the outer counterclockwise composite, apply the Modularity relation (20) to arrive at the uppermost clockwise

composite below:

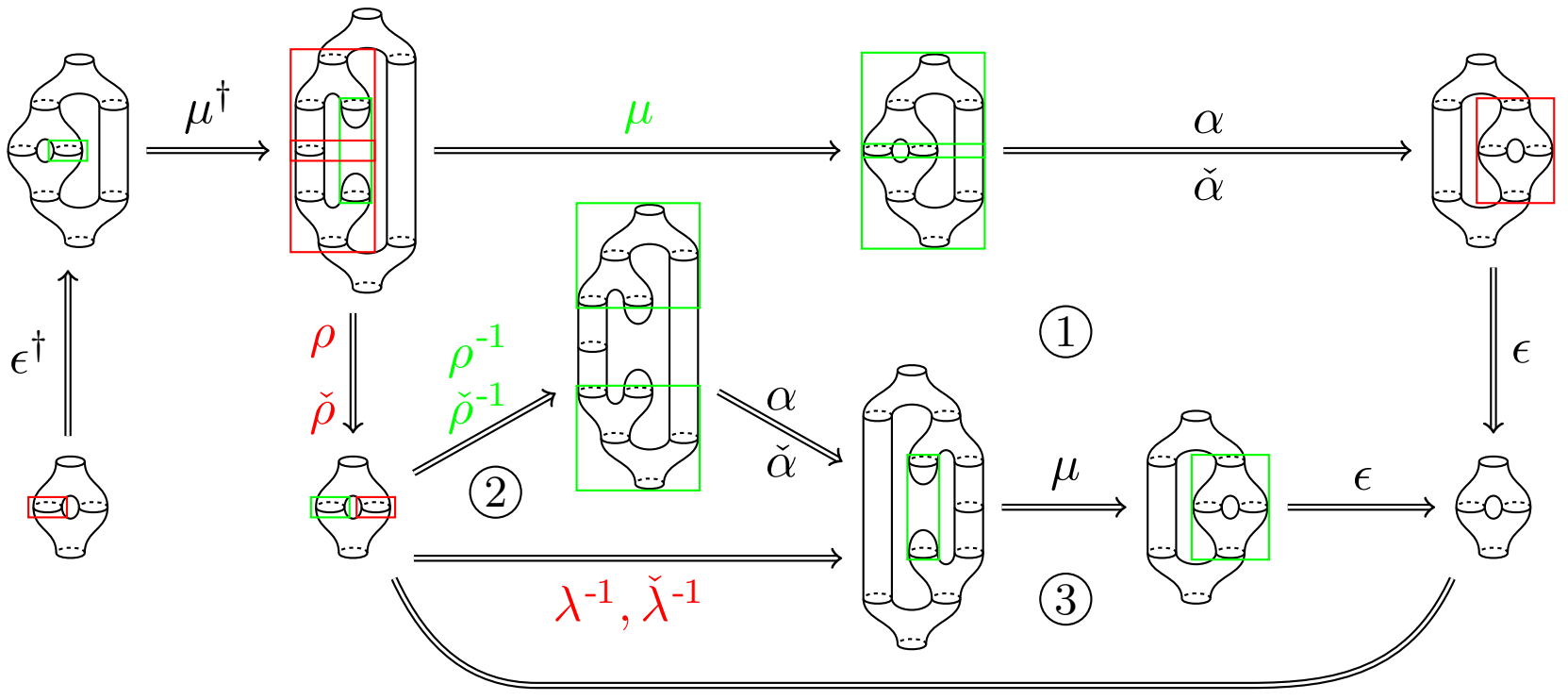

In diagram ① we insert the identity and use naturality, while diagram ② follows from the Monoidal unit relation (6). Diagram ③ is a consequence of Pivotality in the form of $(31)^z$, since the clockwise composite of diagram (C) is simply ${}^*\lambda \circ (\lambda^{-1})^*$. For a similar reason, Pivotality in the form of (31) establishes that the resulting counterclockwise composite equals the identity. Hence $A \circ A^\dagger = \mathrm{id}$. A similar argument establishes $A^\dagger \circ A = \mathrm{id}$.

(ii) ⇒ (i). If $A$ is invertible with inverse $A^\dagger$, then III will be invertible with inverse $\mathrm{III}^\dagger$, and so the Charge Conjugation diagram (Proposition 37 below) holds, and we can pre- and post-compose it as follows:

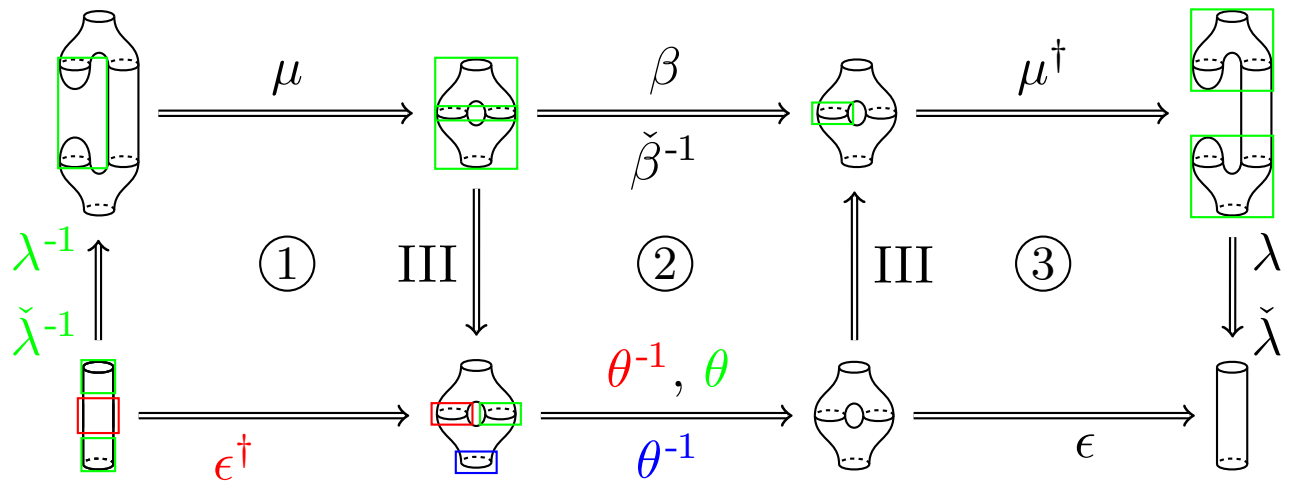

Diagram ② is Proposition (37). Diagram ① is proved as follows. Expand III as $\theta_l \circ A \circ \theta_l$ as in (35). The initial $\theta_l$ falls away due to the Ribbon relation (9). Rearrange the rest by naturality to arrive at the following composite:

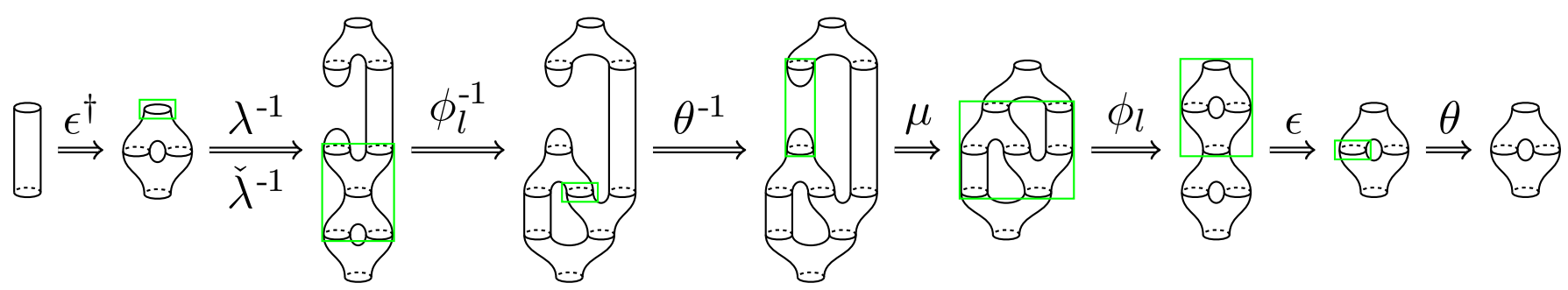

The $\theta^{-1}$ commutes to just after the $\epsilon^\dagger$ using naturality and the Ribbon relation (14). This allows the $\phi_l^{-1}$ to cancel with the $\phi_l$. An application of Pivotality in the form of Lemma 16 then reduces the composite to $\theta_l \circ \theta_l^{-1} \circ \epsilon^\dagger$, so that the $\theta$ terms cancel, proving (a). Diagram ③ is proved in a similar way.

Expanding out the $\check{\lambda}^{-1}$, $\check{\beta}^{-1}$ and $\check{\lambda}$ in the top row and using the Ribbon relations and Pivotality, the top row becomes which is the bottom composite of the Modularity relation (20). In the lower composite, the lower $\theta^{-1}$ can commute past the $\epsilon$ to become a $\theta^{-1}$ between the $\mu^\dagger$ and the $\mu$, where it falls away due to the Ribbon relation (9). Hence the lower composite is the top composite of (20).

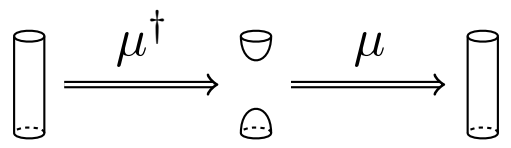

(ii) $\Rightarrow$ (iii) is clear. To prove (iii) $\Rightarrow$ (ii), write $\Gamma$ for the invertible 2-morphism $\alpha^{-1} \circ \phi_r^{-1} \circ \lambda_{\text{bot}}^{-1}$. Then we claim that

(39) $A =$ 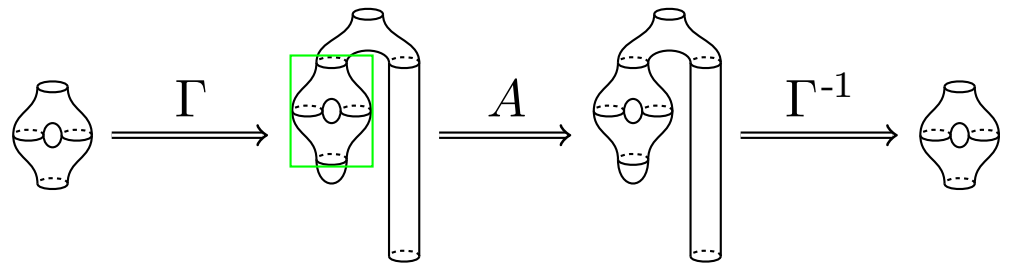 ,

and hence $A$ is invertible on the twice-punctured torus. Indeed, equation (39) is an easy consequence of the relation (T2) in $\mathcal{F}$ from Appendix A.7 which we prove in Section 6.7, and which follows simply from expanding out $A$. $\square$

**Corollary 26.** *The 2-morphisms* III *and* III$'$ *are invertible in* $\mathbf{F}(\mathcal{M})$.

**Remark 27.** The usual condition for a semisimple ribbon category to be modular is that the $\tilde{s}$-matrix with entries

(40) $\tilde{s}_{ij} :=$ 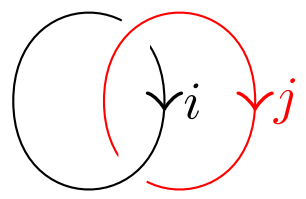

is invertible. Here we are following the notation of [4]. The Modularity relation (20) corresponds (as explained in paper [9]) to the following equation:

(41) 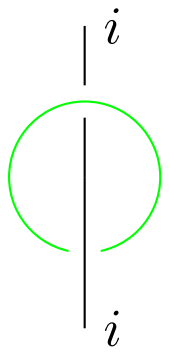

$= \quad p^+ p^- \delta_{i,0}$

Here the unlabelled strand refers to a sum over the simple objects, each weighted by their dimension. It is known that in a semisimple ribbon category, (41) follows from (40) [4, Theorem 3.1.19]. Kerler has also shown the converse [28, Lemma 17]. Hence, we can interpret Proposition 25 as giving a step-by-step *cobordism proof* of

Kerler's result, so that an analog of the identity will be true in *any* target bicategory, not just the usual target of 2–vector spaces. Indeed Kerler was led to his formula from similar topological considerations [28, page 62].

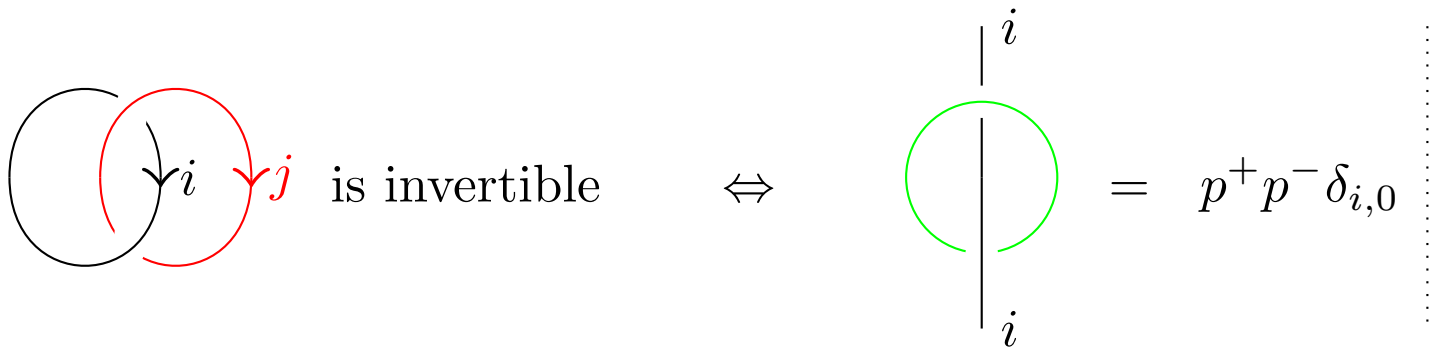

4.2. II **and** III **are self-dual.** We will need the following duality properties of II and III for later use in Section 6.

**Proposition 28.** *In* $\mathbf{F}(\mathcal{M})$*,* $A^* = A$*.*

*Proof.* Recall the following equation for $A$:

(42) $A$ =

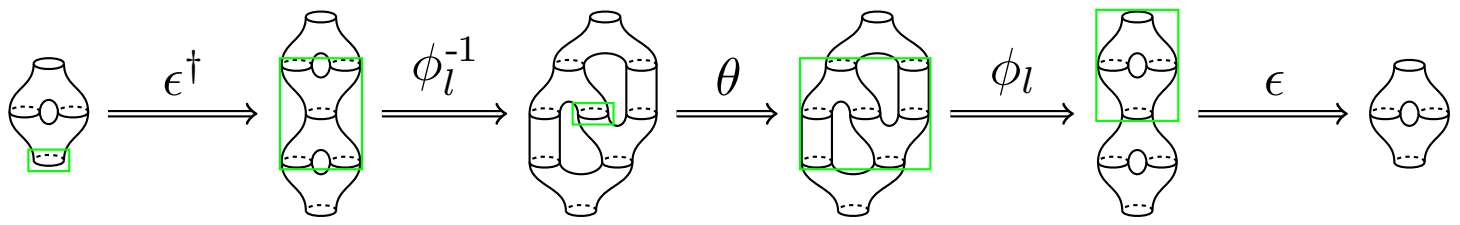

Similarly, if we use the properties of right duals from Section 3.1, combined with Lemma 18, we have the following equation for $A^*$:

(43) $A^*$ =

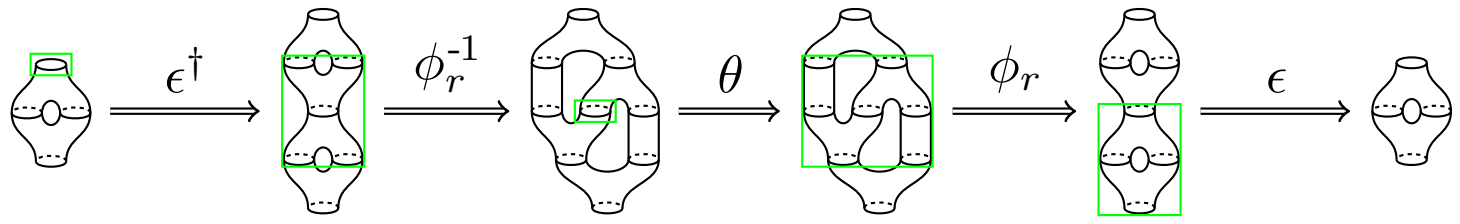

We must show that these are equal. Indeed, analogous manipulations as in the proof of Lemma 22 will show that the composite on the right hand side of (43) is also equal to the right hand side of Lemma 22. Hence $A^* = A$. □

**Corollary 29.** *In* $\mathbf{F}(\mathcal{M})$*, the following equation holds:*

(44) $A$ =

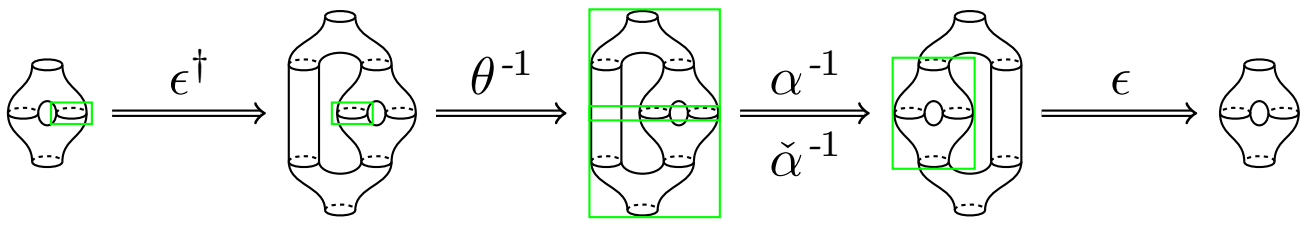

*Proof.* Follows from taking the right dual of both sides of Lemma 22, and using $\alpha^{**} = \alpha$ from Lemma 19. □

**Corollary 30.** *In* $\mathbf{F}(\mathcal{M})$*,* $\mathrm{III}^* = \mathrm{III}$ *and similarly* $(\mathrm{III}')^* = \mathrm{III}'$*.*

**Proposition 31.** *In* $\mathbf{F}(\mathcal{M})$*,* $\mathrm{II}^* = \mathrm{II}$*; equivalently,* $\phi_r \circ \theta \circ \phi_r^{-1} = \phi_l \circ \theta \circ \phi_l^{-1}$*.*

*Proof.* The fact that the two statements are equivalent follows from Lemma 18, since II is by definition $\phi_l \circ \theta \circ \phi_l^{-1}$ so that $\mathrm{II}^*$ is $\phi_r \circ \theta \circ \phi_r^{-1}$. To prove that these are equal, expand out $\phi_r^{-1}$ using $\eta^\dagger$ and $\epsilon^\dagger$, and expand $\epsilon^\dagger$ using Lemma 14. This allows us to move the $\theta$ to near the beginning, and we arrive at the following expression for $\phi_r \circ \theta \circ \phi_r^{-1}$:

(45)

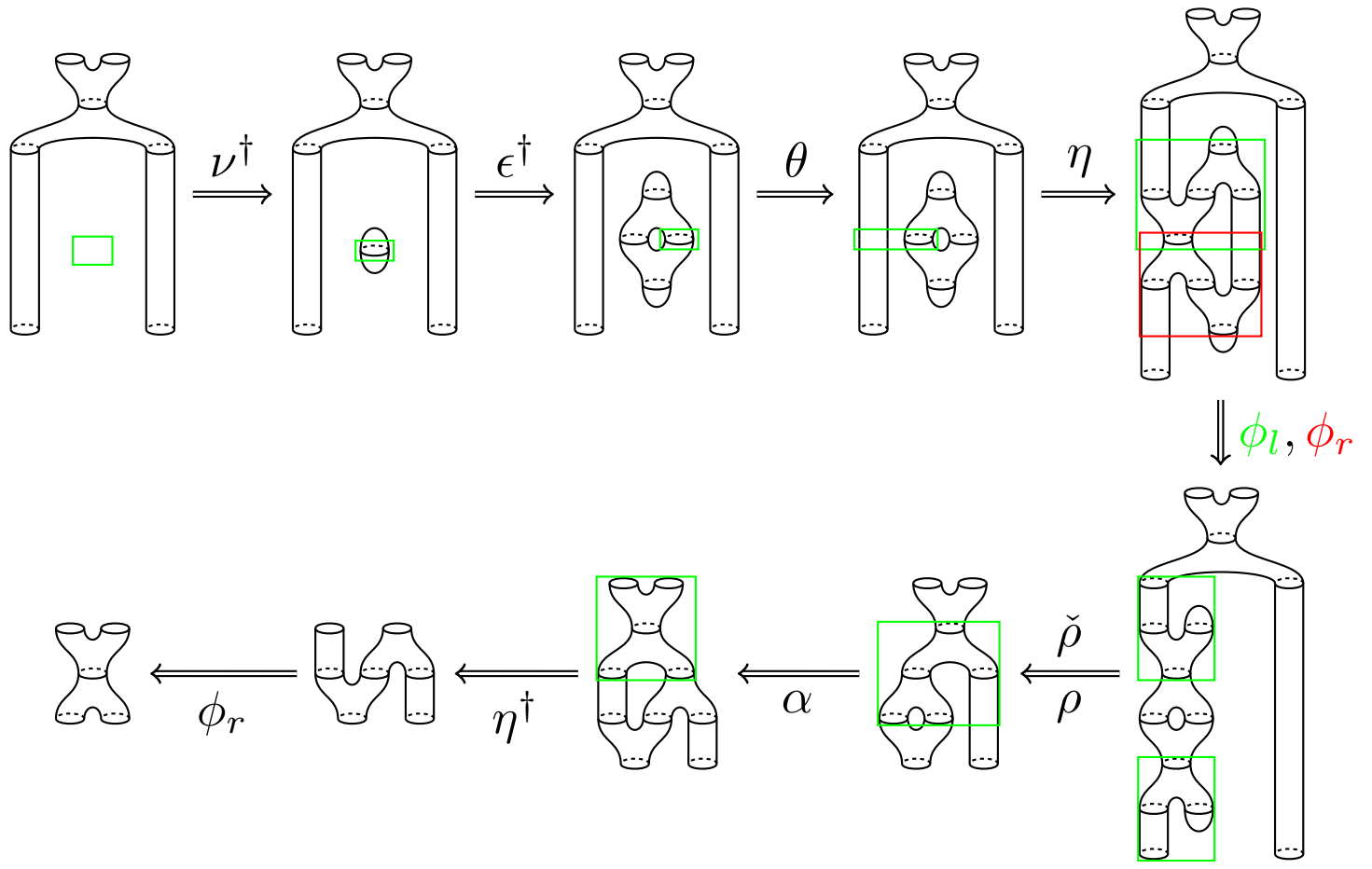

Doing the same manipulations for $\phi_l \circ \theta \circ \phi_l^{-1}$, we arrive at $(45)^z$. To prove that these two expressions are equal, use the Ribbon relation (14) to move the $\theta$ around the bend, expand out the Frobeniusators and use the adjunction equations. □

**Corollary 32.** *In* $\mathbf{F}(\mathcal{M})$, $\mathrm{II} = \mathrm{II}^z$, $A = A^z$, *and* $\mathrm{III}' = \mathrm{III}^z$.

*Proof.* Unravelling the definitions in Appendix B, $\phi_l^z = \phi_r$ and vice-versa. Hence the rotation of ($\mathrm{II} = \phi_l \circ \theta \circ \phi_l^{-1}$) about the $z$-axis is $\phi_r \circ \theta \circ \phi_r^{-1}$. Proposition 31 shows that these are equal. This in turn implies that $A = A^z$ and that $\mathrm{III}' = \mathrm{III}^z$. □

4.3. **Mapping class group identities.** In this section we establish that the analogues of certain identities in the mapping class group of a torus also hold in $\mathbf{F}(\mathcal{M})$. We first need the following.

**Lemma 33.** *In* $\mathbf{F}(\mathcal{M})$, III *decomposes in 'braided form' as follows:*

(46)

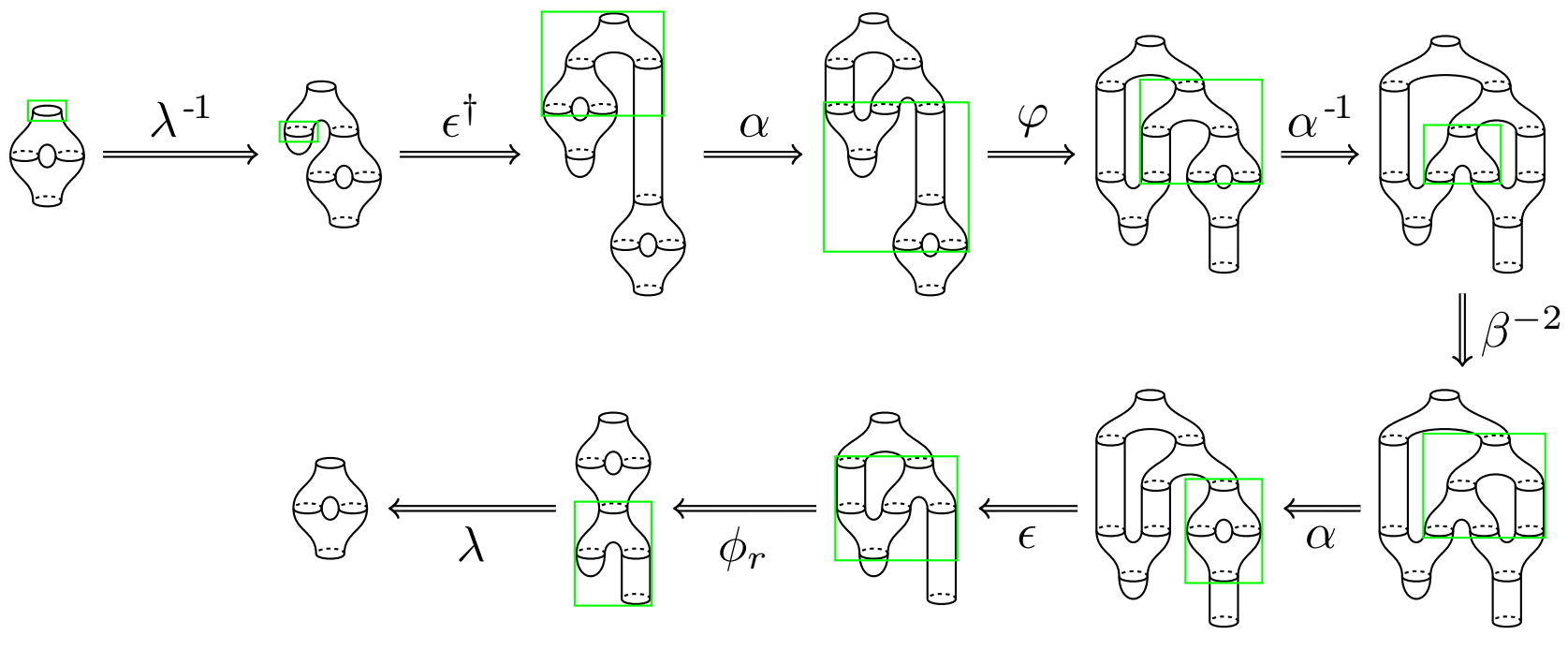

*Moreover,* $\mathrm{III}^{-1}$ *decomposes in precisely the same way, except the* $\beta^{-2}$ *term is replaced by a* $\beta^2$*.*

*Proof.* We must show that the following diagram commutes:

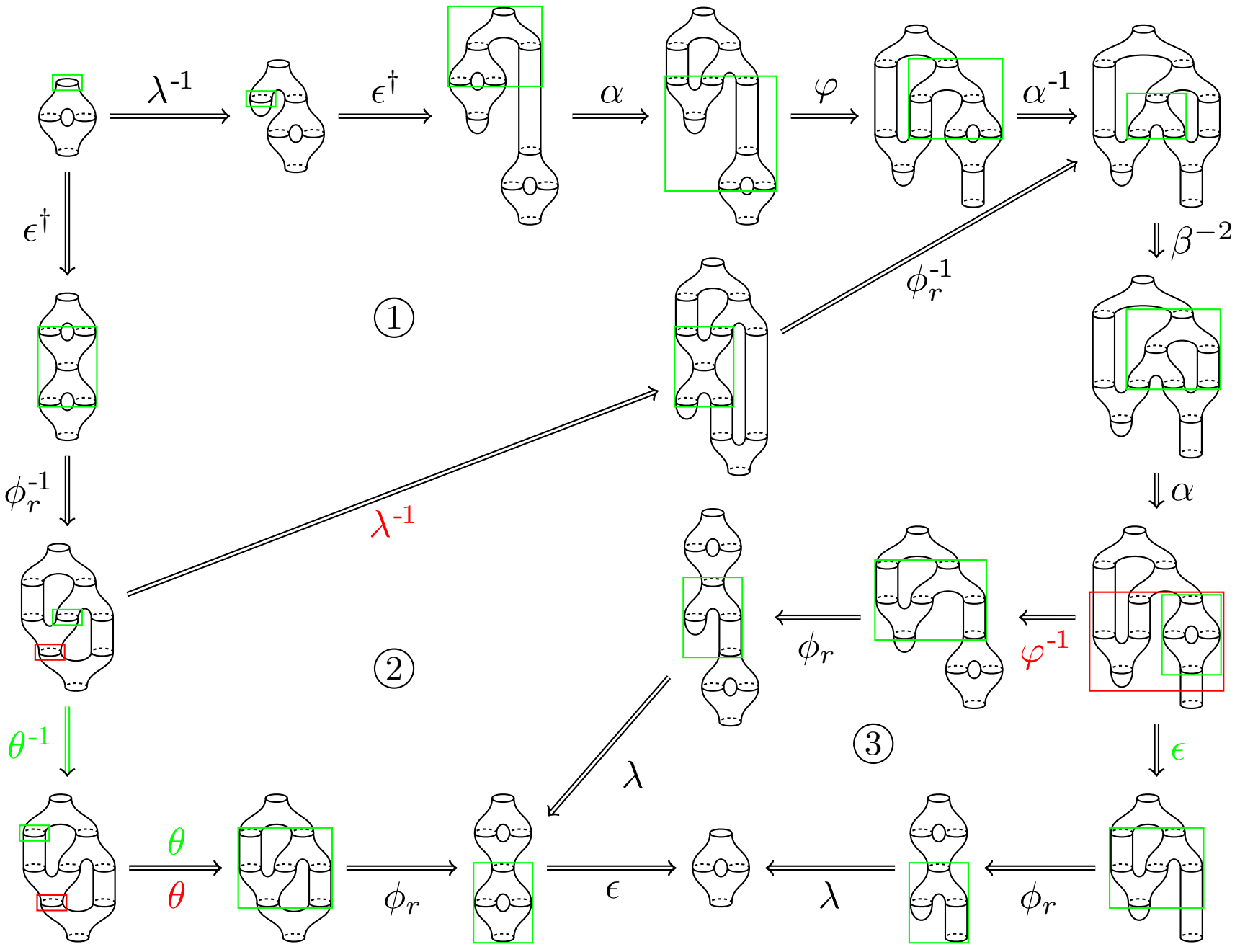

The outer counter clockwise composite is the definition of III, but we have used the form of $A$ given in the right hand side of (43) and have rearranged the $\theta$ terms somewhat. Diagram ① follows from naturality, expansion of the Frobeniusator and the pentagon equation (5). Diagram ② follows from reversing the direction of the $\theta$ terms and then moving them clockwise by naturality (and the Ribbon relation for the red $\theta$) to cancel with the $\beta^{-2}$ using the Balancing relation. Diagram ③ follows from expansion of the Frobeniusator and the unit axioms. This proves the first claim.

To prove the second claim, recall that Proposition 25(ii) showed that $\mathrm{III}^{-1} = \mathrm{III}^{\dagger}$, a consequence of the Modularity relation. If we had started with $\mathrm{III}^{\dagger}$ in the counterclockwise composite above, the only thing that would have changed is that the $\theta$ terms would have been replaced with their inverses. Hence the same argument which proved Diagram ② would still be valid, as long as we replace the $\beta^{-2}$ with a $\beta^2$. ☐

**Definition 34.** We define the *charge conjugation 2-morphism* $c$ as the following composite in $\mathbf{F}(\mathcal{M})$:

$c \quad :=$ 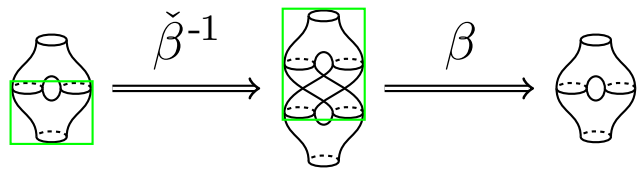

On the torus, the 2-morphism $c$ represents the hyperelliptic involution element in the mapping class group. We show this in two stages. Firstly, it is an involution.

**Lemma 35.** *In* $\mathbf{F}(\mathcal{M})$*,* $c^2 = \mathrm{id}$ *on the torus. On the once-punctured torus,* $c^2 = \theta_{\mathrm{top}}$.

*Proof.* We give the following argument:

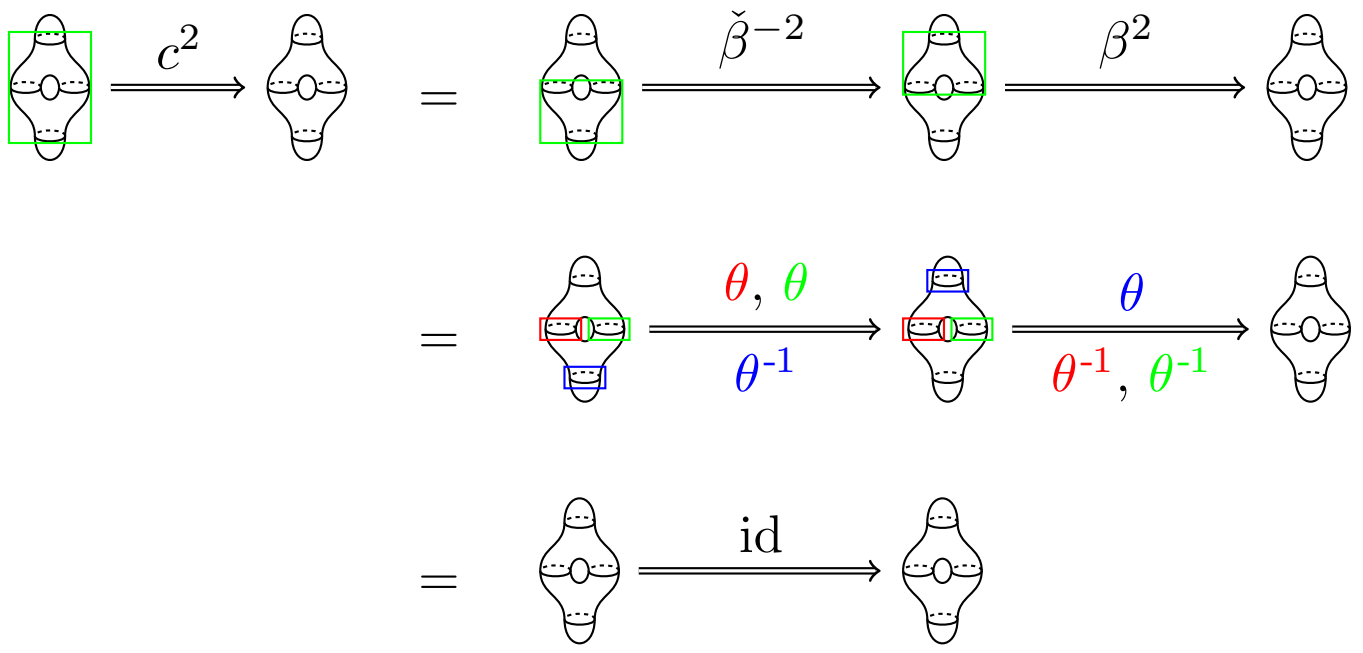

Here we have used the Ribbon relations (8) and (9) as well as Lemma 11. On the once-punctured torus, the $\theta$ on the belt of the pants will remain. □

**Remark 36.** The reason for calling $c$ 'charge conjugation' is as follows. In the TQFT $Z$ corresponding to a modular category $C$ [9] over an algebraically closed field, $Z(S^1 \times S^1)$ is spanned by the isomorphism classes of simple objects $X_i$ of $C$, and $Z(c)$ is the linear map which sends $X_i \mapsto X_i^*$. See also [4, Theorem 3.1.7].

**Proposition 37** (Charge conjugation)**.** *In* $\mathbf{F}(\mathcal{M})$ *the following equation holds, together with its rotated version* $(47)^x$*:*

(47)

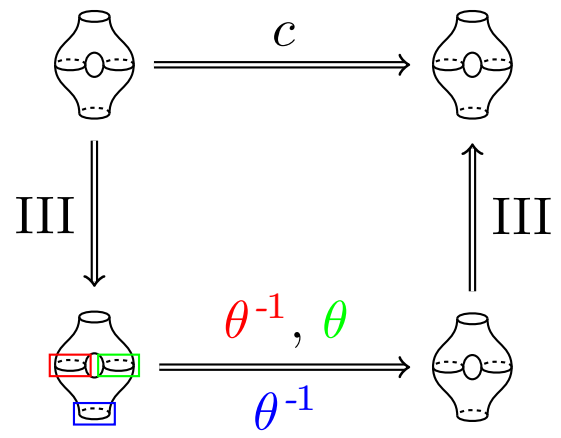

*Proof.* Using the invertibility of III, it suffices to show that the composite $X$, obtained by starting at the top right and proceeding counterclockwise to the bottom right, equals III$^{-1}$. We will show this using the braided decomposition of III and III$^{-1}$ from Lemma 33. Using (i) naturality to move the $\beta^{-1}$ and $\check{\beta}$ past some terms, (ii) the Ribbon relation (14) to move the red $\theta^{-1}$ around the 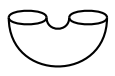 segment, and (iii)

the definition of $\check{\beta}$ as $\beta^*$, we find that $X$ is equal to the following composite:

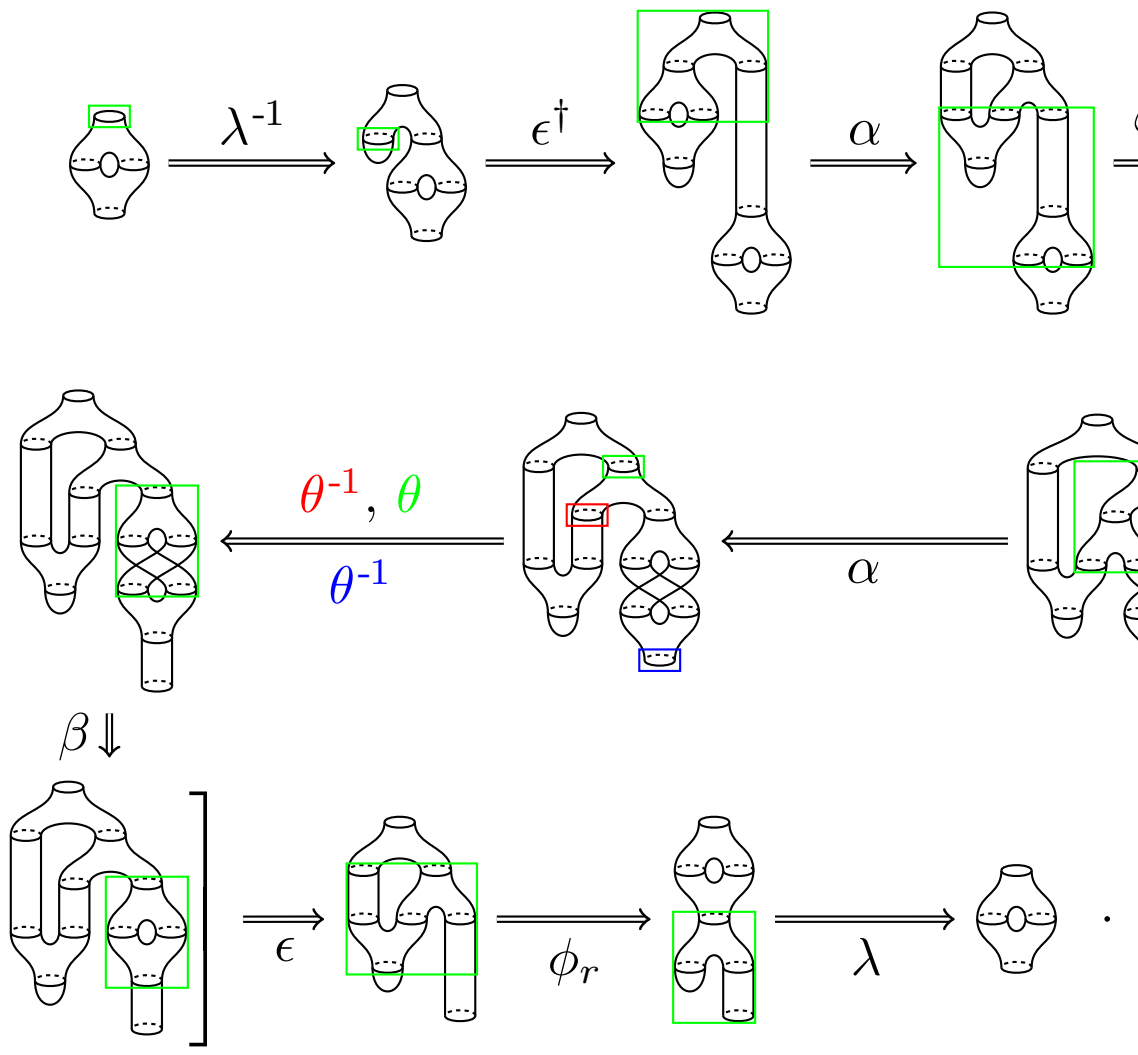

Using naturality and the Ribbon relations, we can eliminate all the Dehn twists from the diagram. We are left with showing that bracketed composite above equals the corresponding segment from the braided decomposition of $\mathrm{III}^{-1}$ from Proposition 33. This amounts to checking that the following diagram commutes, where the clockwise composite is the bracketed composite above:

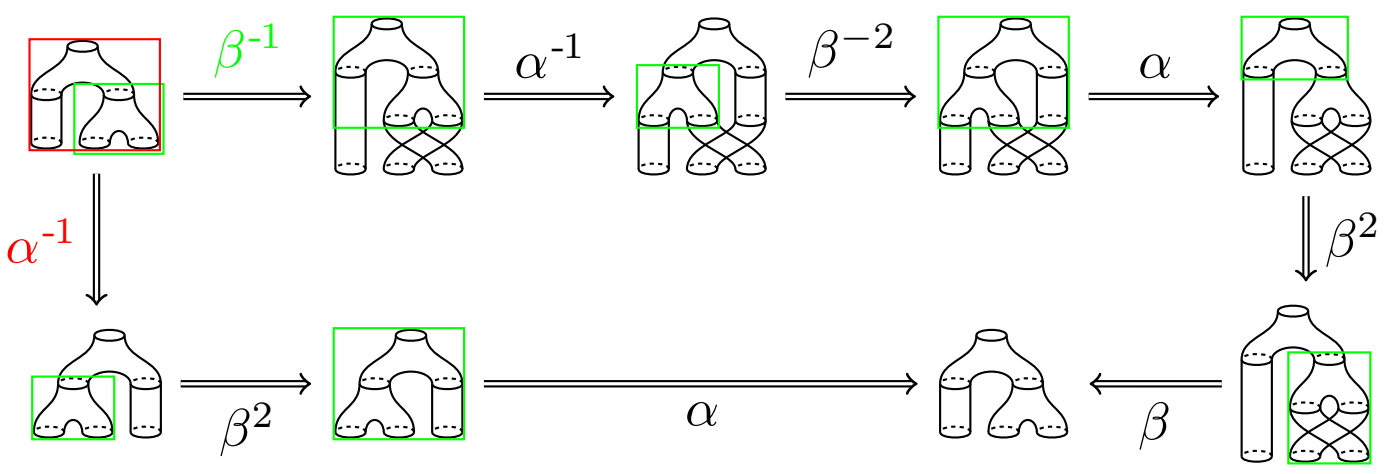

This diagram commutes by braided monoidal coherence, since both sides produce the same underlying braid. □

**Corollary 38.** (Compare [4, Theorem 3.1.18].) *In* $\mathbf{F}(\mathcal{M})$, $\mathrm{III}^2 = c$ *on the once-punctured torus and hence* $\mathrm{III}^4 = \theta_{\mathrm{top}}$ *on the once-punctured torus.*

*Proof.* Cup off the copants of (47). The Ribbon relations (9) and (14) cause the $\theta$ terms to fall away and hence $\mathrm{III}^2 = c$. □

**Corollary 39.** (Compare [4, Theorem 3.1.18].) *In* $\mathbf{F}(\mathcal{M})$, $\mathrm{III}^4 = \mathrm{id}$ *on the torus.*

## 5. The anomaly

In this section we continue to work with the modular presentation $\mathcal{M}$. Our goal is to study the nature of the following 'anomaly' 2-morphism $x$ in the 2-category $\mathbf{F}(\mathcal{M})$.

**Definition 40.** In $\mathbf{F}(\mathcal{M})$, we define the *anomaly on the sphere* $x'$ as the following composite 2-morphism on the 2-sphere:

$x' \quad :=$ 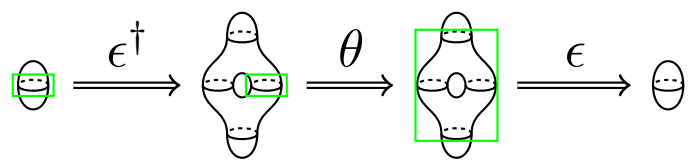

We define the *anomaly on the cylinder* $x$ as the following composite:

$x \quad :=$ 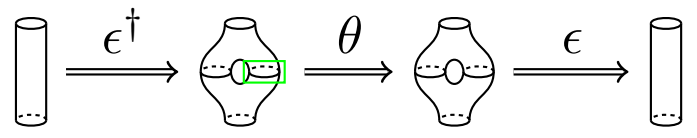

Thus, the anomaly-free modular presentation $\mathcal{O}$ differs from the modular presentation $\mathcal{M}$ by adding the relation $x' = \mathrm{id}$.

The key result in this section is establishing that the relation (T4) in the 2-Morse presentation $\mathcal{F}$ (see Appendix A) follows from the relations in the anomaly-free modular presentation $\mathcal{O}$. (T4) is the braid relation for Dehn twists, a relation in the mapping class group of a thrice-punctured torus. What we are able to prove directly in the modular presentation $\mathcal{M}$, in Theorem 49, is an anomalous form of this equation, where one side has been postcomposed by attaching the composite 2-morphism $x$—attaching $x$ is the analogue of the blow-up move in the Kirby calculus. We will see that $x = \mathrm{id}$ follows from the axiom $x' = \mathrm{id}$ of the presentation $\mathcal{O}$, and so the relation (T4) follows.

In [8] we show that there is a symmetric monoidal equivalence

$$A : \mathbf{F}(\mathcal{M}) \xrightarrow{\sim} \mathbf{Bord}_{123}^{\mathrm{csig}} \tag{48}$$

where $\mathbf{Bord}_{123}^{\mathrm{csig}}$ is the *componentwise signature extension* of the oriented bordism bicategory $\mathbf{Bord}_{123}^{\mathrm{or}} \simeq \mathbf{F}(\mathcal{O})$ (see Remark 6). Under this equivalence, attaching a copy of $x'$ to a 3-manifold $M$ corresponds to attaching a copy of $\mathbb{CP}^2$ to its bounding 4-manifold and hence changing the signature of that bounding manifold by 1.

**Remark 41.** Let $S$ be a 1-morphism in $\mathbf{F}(\mathcal{M})$ representing a connected closed surface. The equivalence (48) shows that the automorphism group of $S$ in $\mathbf{F}(\mathcal{M})$ is isomorphic to the automorphism group of $A(S)$ in $\mathbf{Bord}_{123}^{\mathrm{csig}}$. This latter group, after unravelling the definitions, is the *signature central extension* $\Gamma^{\mathrm{sig}}(\Sigma)$ of the oriented mapping class group of the surface represented by $S$, in the sense of [30, 35].

5.1. **Centrality of the anomaly.** The following lemma shows that the anomaly on the cylinder $x$ is the same as creating a 2-sphere, applying the anomaly $x'$ on the 2-sphere, and then attaching the 2-sphere to the cylinder.

**Lemma 42.** *The following equation holds in* $\mathbf{F}(\mathcal{M})$*, together with its rotated version* $(49)^z$*:*

(49) 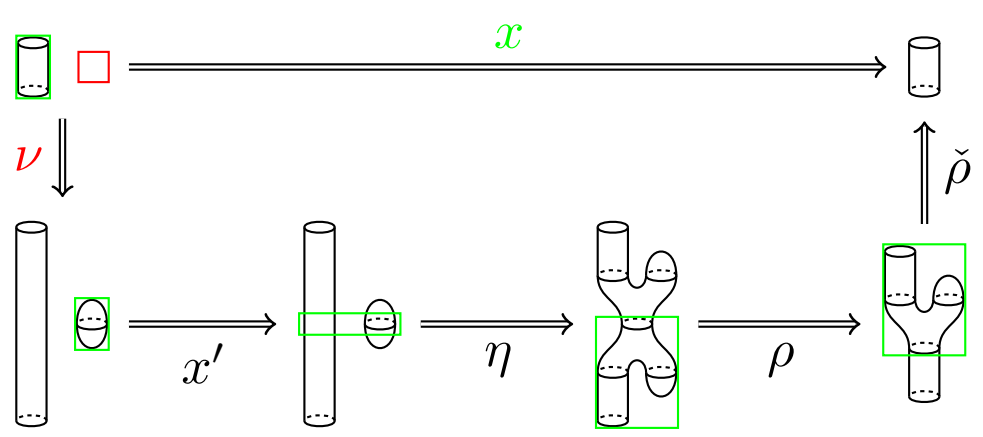

*Proof.* Using Lemma 14 and naturality, we can rewrite the definition of $x$ as the outer clockwise composite below:

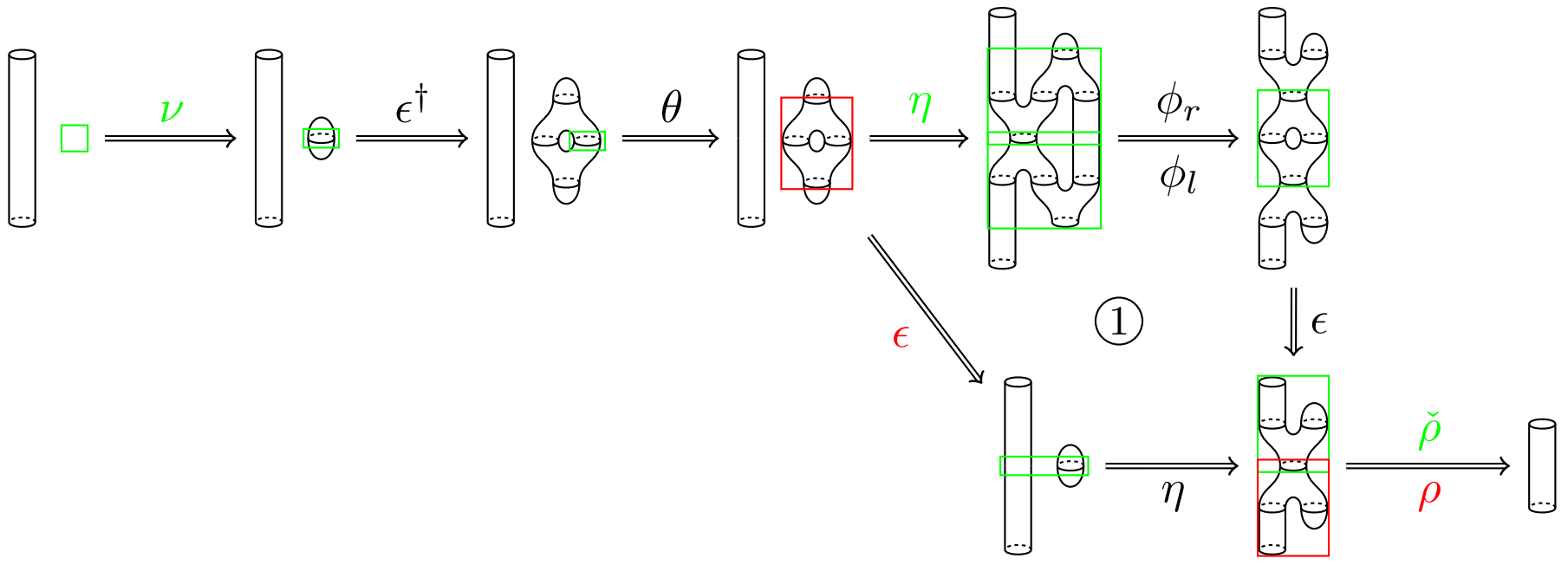

Diagram ① commutes from expansion of the Frobeniusators and the $\eta$–$\epsilon$ adjunction equations. The rotated version $(49)^z$ is proved similarly. ☐

**Corollary 43.** *For the modular presentation $\mathcal{M}$, the anomaly-freeness relation* (21) $x' = \mathrm{id}$ *on the sphere is equivalent to the equation* $x = \mathrm{id}$ *on the cylinder.*

*Proof.* Substituting the anomaly-freeness relation (21) into the expression in Lemma 42 gives

$x \quad =$ 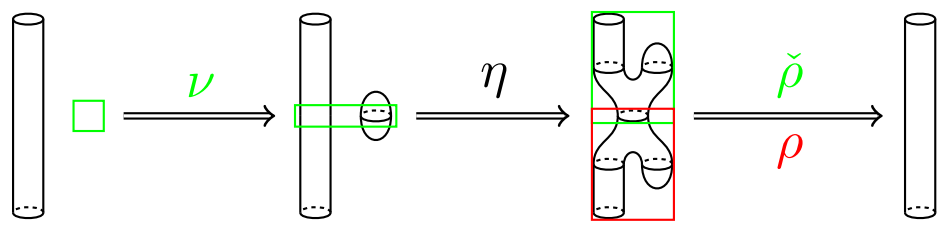

which equals the identity, by expanding $\check{\rho}$ as $(\rho^{-1})^*$, and using naturality and the adjunction equations. Conversely, if $x = \mathrm{id}$, then we can cap off the top and bottom to get $x' = \mathrm{id}$. ☐

**Lemma 44.** *The anomaly on the cylinder $x$ satisfies the following equation in* $\mathbf{F}(\mathcal{M})$*:*

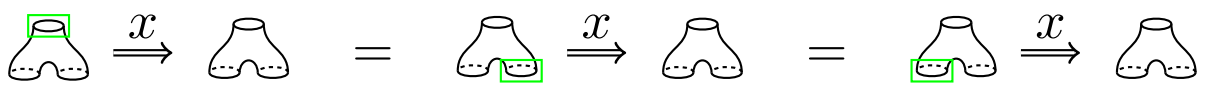

*Similar equations hold on* 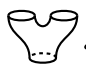*.*

*Proof.* We prove the first equation, the others are similar. Using Lemma 42, this comes down to proving that attaching a sphere to the belt of a pair-of-pants is the same as attaching it to the right leg:

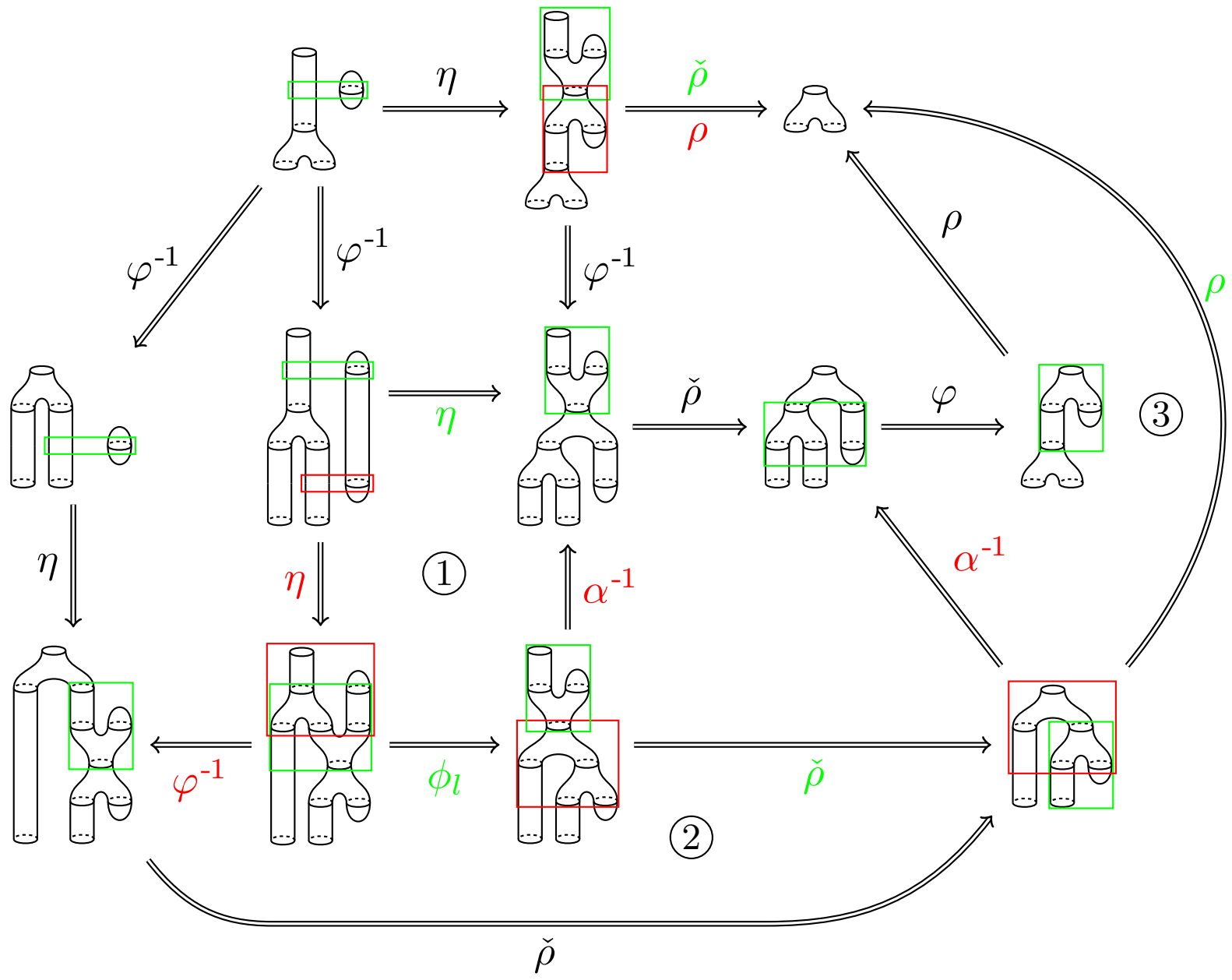

Diagram (1) follows from expansion of $\phi_l$ and the $\eta$–$\epsilon$ adjunction equation. Diagram (2) is diagram (B1) from Appendix A which is proved in Section 6.5 using only the Monoidal, Pivotality, and Biadjoint relations. Diagram (3) follows from monoidal coherence. □

**Corollary 45.** *In* $\mathbf{F}(\mathcal{M})$*, the anomaly on the cylinder* $x$ *is central with respect to horizontal and vertical*[3] *composition of 2-morphisms.*

*Proof.* Using Lemma 44, we can commute $x$ past the generators, irrespective of where it is situated. For instance,

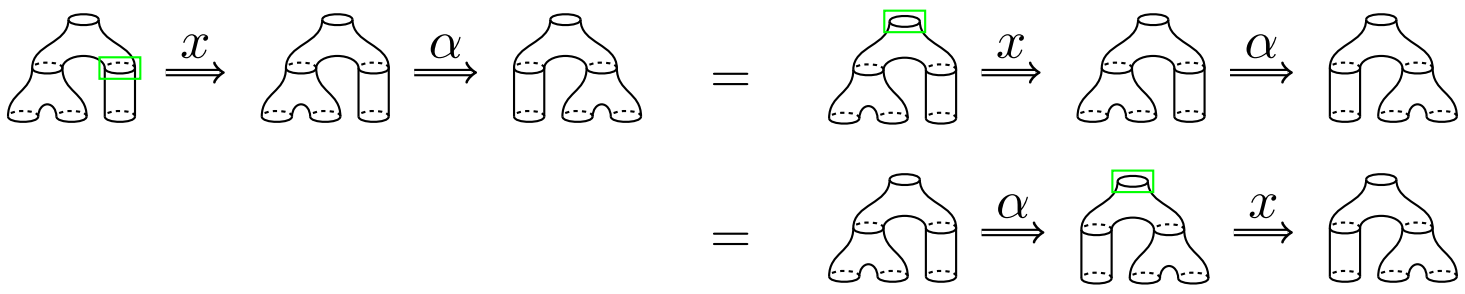

□

5.2. **Invertibility of the anomaly.**

**Lemma 46.** *In* $\mathbf{F}(\mathcal{M})$*,* $x'^{\dagger}x' = \mathrm{id} = x'x'^{\dagger}$*. Equivalently,* $x^{\dagger}x = \mathrm{id} = xx^{\dagger}$*.*

*Proof.* The equivalence of the two assertions follows from Lemma 42. We prove it for $x$, using colours to keep track of things:

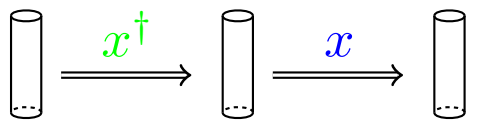

[3] Recall from Section C.2 that due to the graphical calculus conventions in this paper, what we call *horizontal composition* of 2-morphisms is more often called *vertical composition*. That is, if $\alpha\colon f \Rightarrow g$ and $\beta\colon g \Rightarrow h$, then their horizontal composite is a 2-morphism $\beta \circ \alpha\colon f \Rightarrow h$.

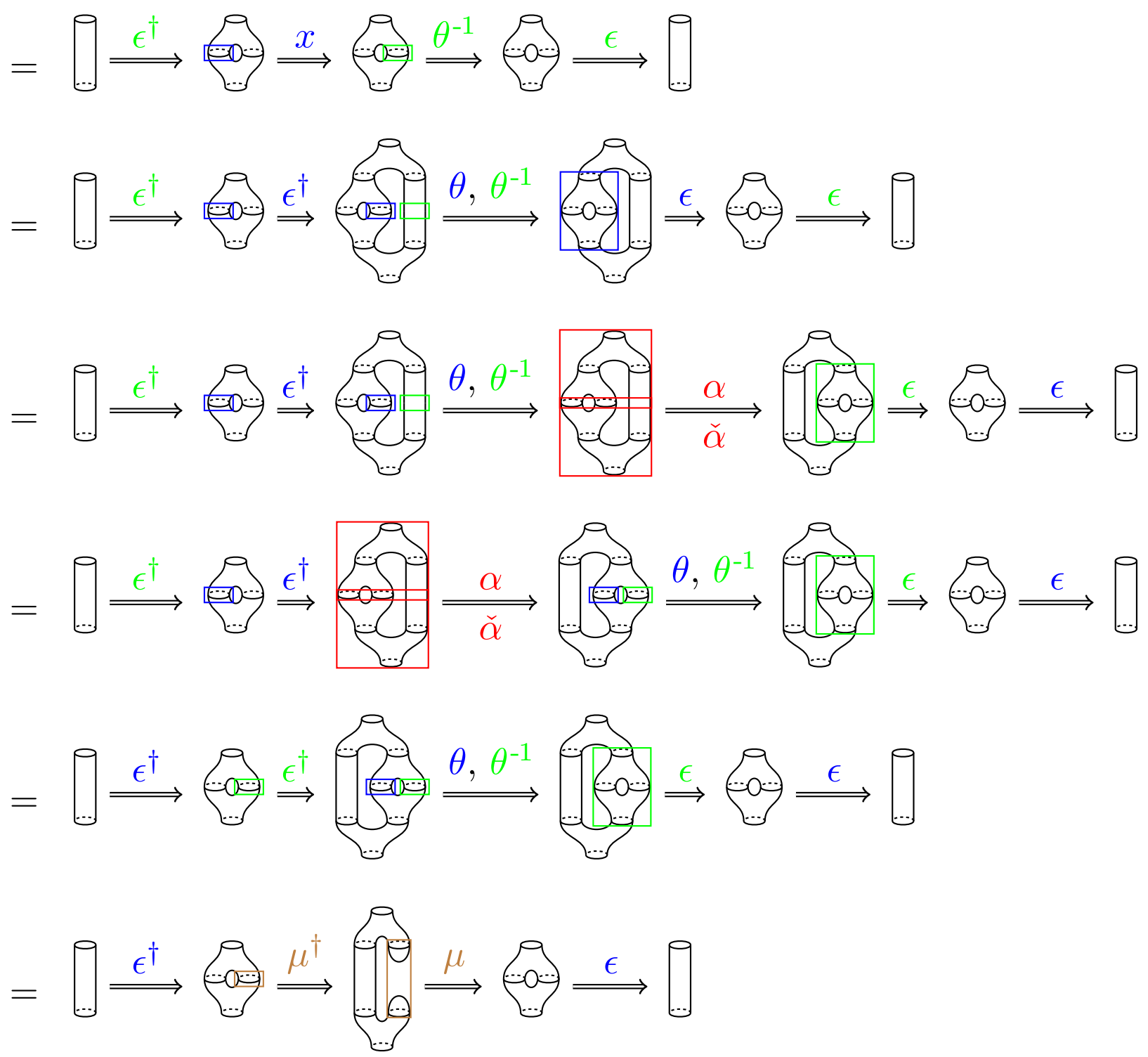

$$= \text{id}$$

The first equality uses the centrality of $x$ from Lemma 44, the second and fourth use naturality, the third and fifth use Lemma 23, the sixth uses the Modularity relation (20), and the final equality uses Pivotality in the form of Lemma 15. Since $x$ is central, this also shows $x^\dagger x = \text{id}$. □

5.3. **A central extension.** Having established that $x$ is invertible (Lemma 46) and central (Corollary 45), we can now use a result from [9] to establish a certain central extension of automorphism groups. Write $T$ as a shorthand for the torus 1-morphism in $\mathbf{F}(\mathcal{M})$:

$$T :=$$ 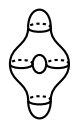

We then have the following result.

**Proposition 47.** *The projection functor $\pi : \mathbf{F}(\mathcal{M}) \to \mathbf{F}(\mathcal{O})$ induces a central extension of groups:*

$$0 \to \mathbb{Z} \xrightarrow{i} \mathrm{Aut}_{\mathbf{F}(\mathcal{M})}(T) \xrightarrow{\pi} \mathrm{Aut}_{\mathbf{F}(\mathcal{O})}(T) \to 1$$

*Here, $i\colon \mathbb{Z} \to \mathrm{Aut}_{\mathbf{F}(\mathcal{M})}(T)$ sends $n \mapsto x^n$.*

*Proof.* Indeed, $\pi$ is surjective by definition and has kernel generated by $i(1) = x$, using the centrality of $x$ from Corollary 45. Hence, we must show that $x$ has infinite order.

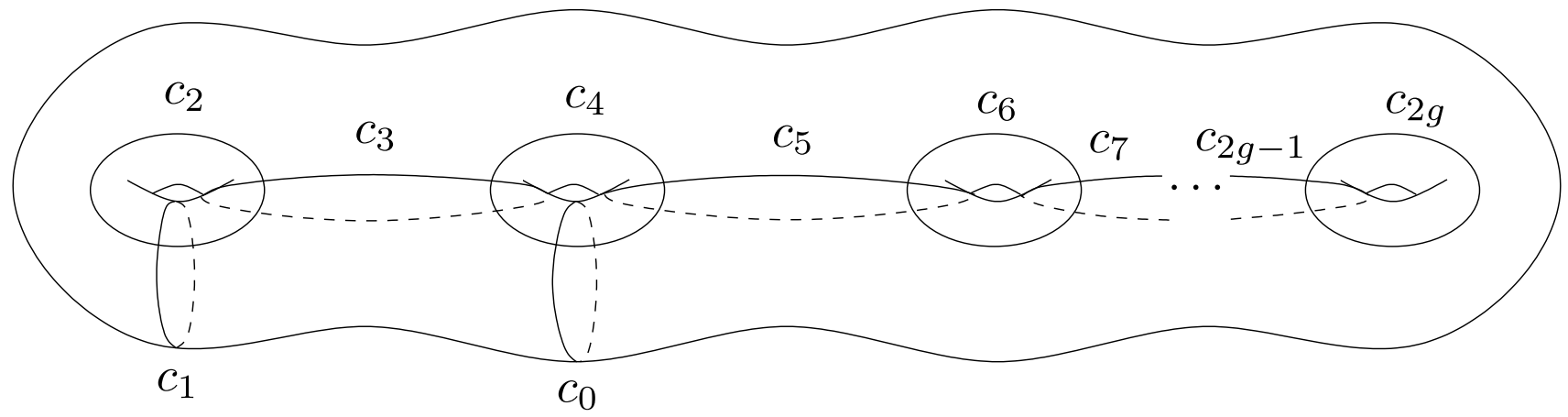

FIGURE 3. The Humphries generators for $\Gamma(\Sigma_g)$ are the Dehn twists $a_i$ about these simple closed curves $c_i$. Taken from [11].

We show this as follows. One of the results in [9] is that every modular category $C$ gives rise to a certain representation $Z_C^+$ of $\mathbf{F}(\mathcal{M})$.[4] Moreover, in this representation $Z_C^+(x)$ corresponds to multiplication by the constant $\sqrt{p^+/p^-}$ calculated from $C$, where $p^+$ and $p^-$ are the Gauss sums of $C$. Let $C_k$ be the modular category arising from SU(2) Chern-Simons at level $k$. Then [21]

$$\sqrt{\frac{p^+}{p^-}} = \exp\left(2\pi i \frac{3k}{8(k+2)}\right)$$

which has arbitrarily high order for large enough $k$. This proves the claim. □

**Remark 48.** In fact, $\mathrm{Aut}_{\mathbf{F}(\mathcal{M})}(\Sigma)$ is the *signature central extension* [30, 35] of $SL(2,\mathbb{Z})$. Showing that requires, however, the result that $\mathbf{F}(\mathcal{F}) \simeq \mathbf{Bord}_{123}^{\mathrm{or}}$ from [7], whose proof is geometric and requires Cerf theory. The statement and proof of Proposition 47 is purely algebraic and more elementary.

5.4. **The anomalous braid relation for Dehn twists.** In Figure 3 we display the Humphries generators $a_i$ for the mapping class group $\Gamma(\Sigma_g)$ of a closed genus $g$ oriented surface. (Note that one should rotate the picture in Figure 3 ninety degrees to match with the standard genus $g$ surface in our bordism category.) As part of Wajnryb's finite presentation [11] of $\Gamma(\Sigma_g)$, there are the *braid relations*

$$a_i a_j a_i = a_j a_i a_j \quad \text{if } c_i \text{ intersects } c_j. \tag{50}$$

In our bicategorical presentation, the three different types of Dehn twists $a_i$ in Figure 3 are represented as follows:

$$c_0, c_1 \quad \leftrightarrow \quad \theta \tag{51}$$

$$c_3, c_5, \ldots, c_{2g-1} \quad \leftrightarrow \quad \mathrm{II} \tag{52}$$

$$c_2, c_4, \ldots, c_{2g} \quad \leftrightarrow \quad A \tag{53}$$

It turns out that the braid relations (50) are indeed valid relations in $\mathbf{F}(\mathcal{O})$, but are only correct up to a factor of $x$ in $\mathbf{F}(\mathcal{M})$ (they become *anomalous*). The following theorem, whose proof will occupy us for the remainder of the section, proves this for Dehn twists of type (51) with Dehn twists of type (53).

This is a crucial result, for the following reason. We will show in Section 6.7 that in our bicategorical setup, the remaining braid relation, for Dehn twists of type (52) with Dehn twists of type (53), follows from the (51)-intersect-(53) braid relation. The (52)-intersect-(53) braid relation is in fact a relation in the full

[4] Note that the proof of this result in [9] is entirely algebraic and does not depend on the equivalence between $\mathbf{F}(\mathcal{M})$ and $\mathbf{Bord}_{123}^{\mathrm{csig}}$, proved in [8], or on this Proposition 47 itself.

2-Morse presentation $\mathcal{F}$, namely (T4) in Appendix A. The most difficult part of showing that the relations in the modular presentation $\mathcal{M}$ generate the relations in the 2-Morse presentation $\mathcal{F}$ is proving relation (T4).

**Theorem 49.** *The following equation holds in* $\mathbf{F}(\mathcal{M})$*, together with its rotated version* $(54)^z$*:*

(54)

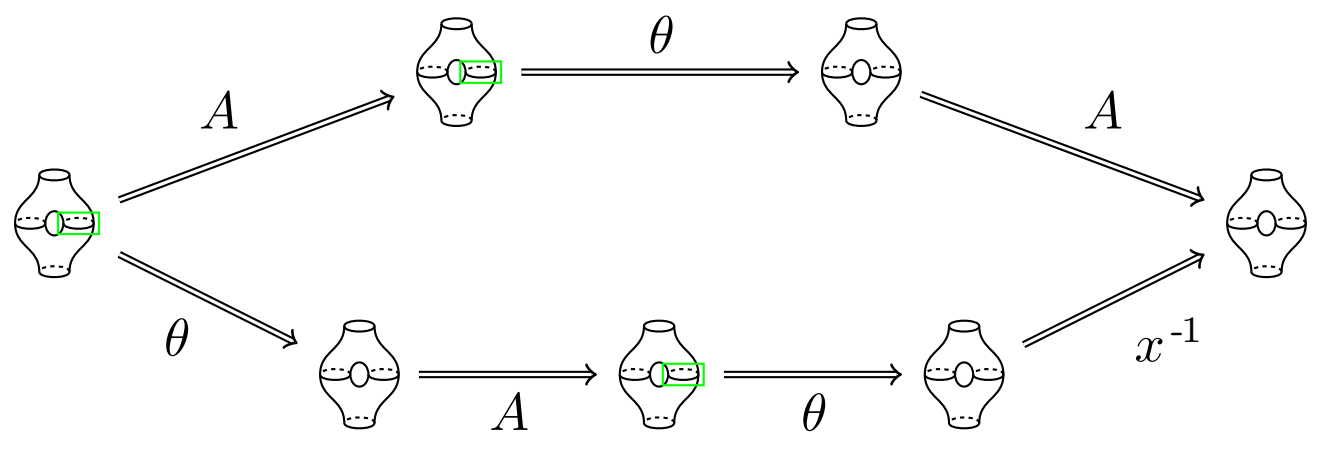

*Proof.* Expand out the $A$ terms using the form of $A$ given in Corollary 29. Then we compute as follows, where the clockwise composite below is the clockwise composite in (54), and the counterclockwise composite below is the counterclockwise composite in (54) (we have moved the $x^{-1}$ term to the front of the composite using Corollary 45):

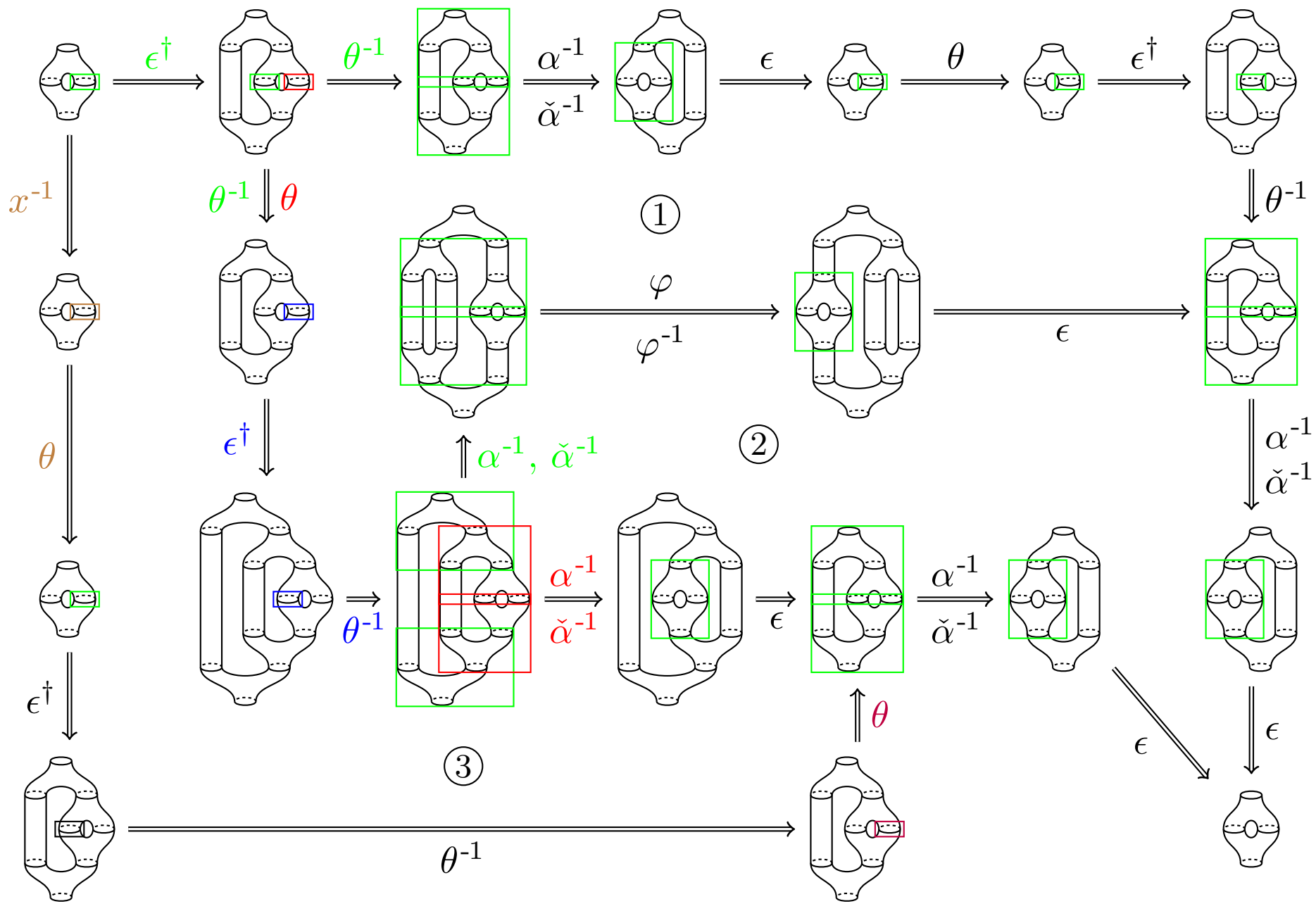

Diagram ① commutes by naturality, diagram ② commutes by a short calculation involving naturality and the pentagon equation, and diagram ③ follows from Lemma 50 and Lemma 56 below. By commuting the final $\theta$ to the end, the counterclockwise composite is thus $\theta_{\text{right}} \circ A \circ \theta_{\text{right}} \circ x^{-1}$. The version $(54)^z$ is proved in a similar way. □

**Lemma 50.** *Diagram ③ from the proof of Theorem 49 will follow if the following diagram commutes:*

(55)

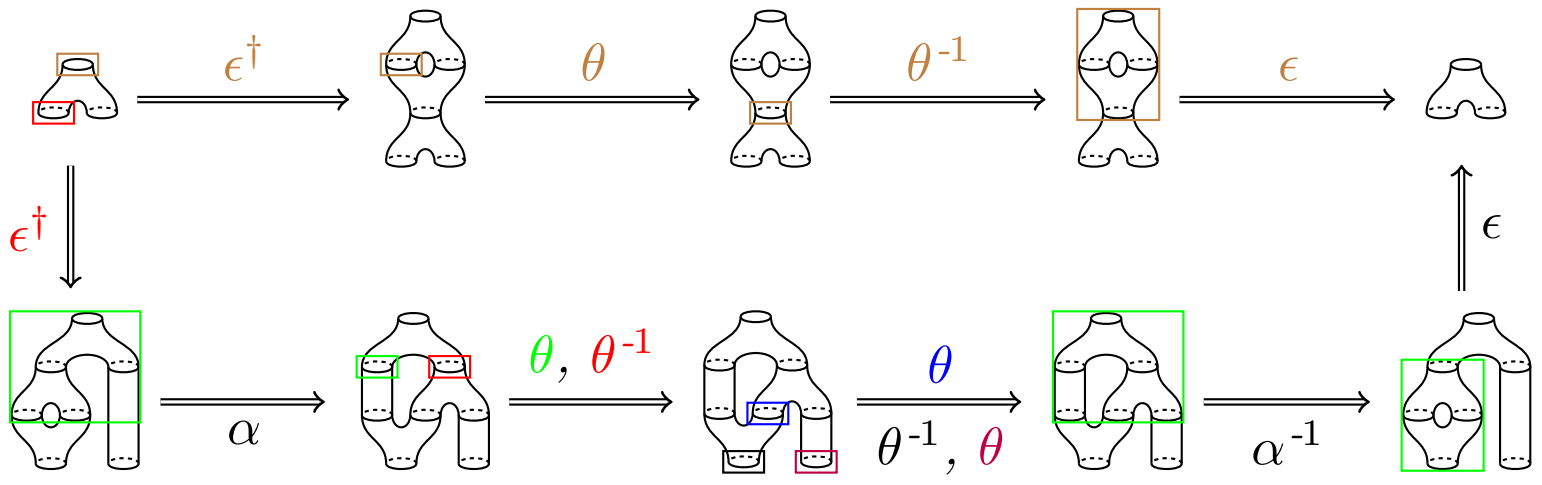

*Proof.* Diagram ③ will follow if the following two composites from 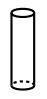 to 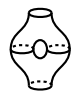 are equal:

(56)

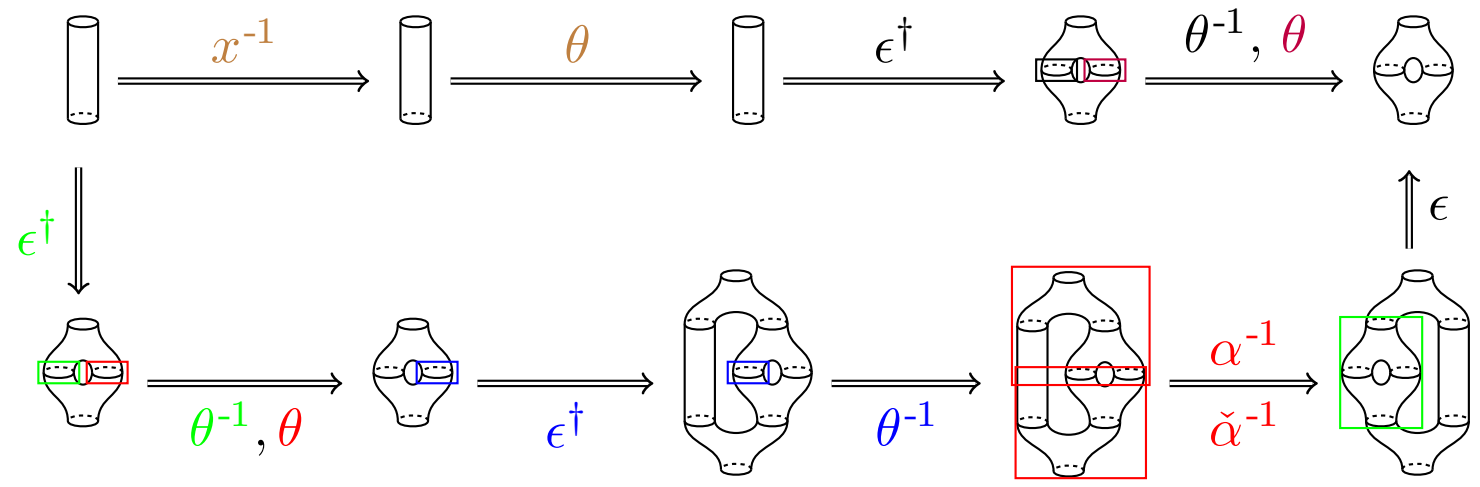

Since the dagger of every relation in $\mathbf{F}(\mathcal{M})$ is also a relation, proving (56) is equivalent to proving its dagger, which is an equality between two composites from 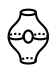 to 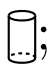; let us call it $X$.

The $\eta$–$\epsilon$ adjunction equations give a bijection of hom-sets

$$\text{2-Hom}_{\mathbf{F}(\mathcal{M})}\left( \right. $$ 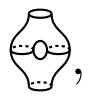, 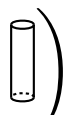 $$\left. \right) \cong \text{2-Hom}_{\mathbf{F}(\mathcal{M})}\left( \right.$$ 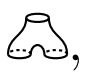, 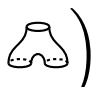 $$\left. \right).$$

Hence the relation $X$ translates into a relation between two composites from 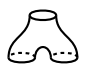 to 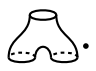. After simplification, this gives exactly diagram (55). □

5.5. **Finishing the proof of the anomalous braid relation.** In this subsection we complete the proof of Theorem 49 by showing that equation (55) indeed holds in $\mathbf{F}(\mathcal{M})$. Our proof is lengthy, even though (55) appears simple. We have been unable to find a simpler proof.

Since this is the key technical calculation in this paper, we pause for a moment to explain the method we used to perform it. In [9] we show that a linear representation of $\mathbf{F}(\mathcal{M})$ corresponds to a modular category $C$ equipped with a choice of square root $p$ of the global dimension. Hence one may ask: what does (55) look like in a linear representation of $\mathbf{F}(\mathcal{M})$? In [9] we describe an *internal string diagram calculus* for working with linear representations of $\mathbf{F}(\mathcal{M})$, where the 'ordinary' string diagrams for the modular category $C$ are thought of as living inside the generators of $\mathbf{F}(\mathcal{M})$.

Here is what (55) looks like in this notation:

(57)

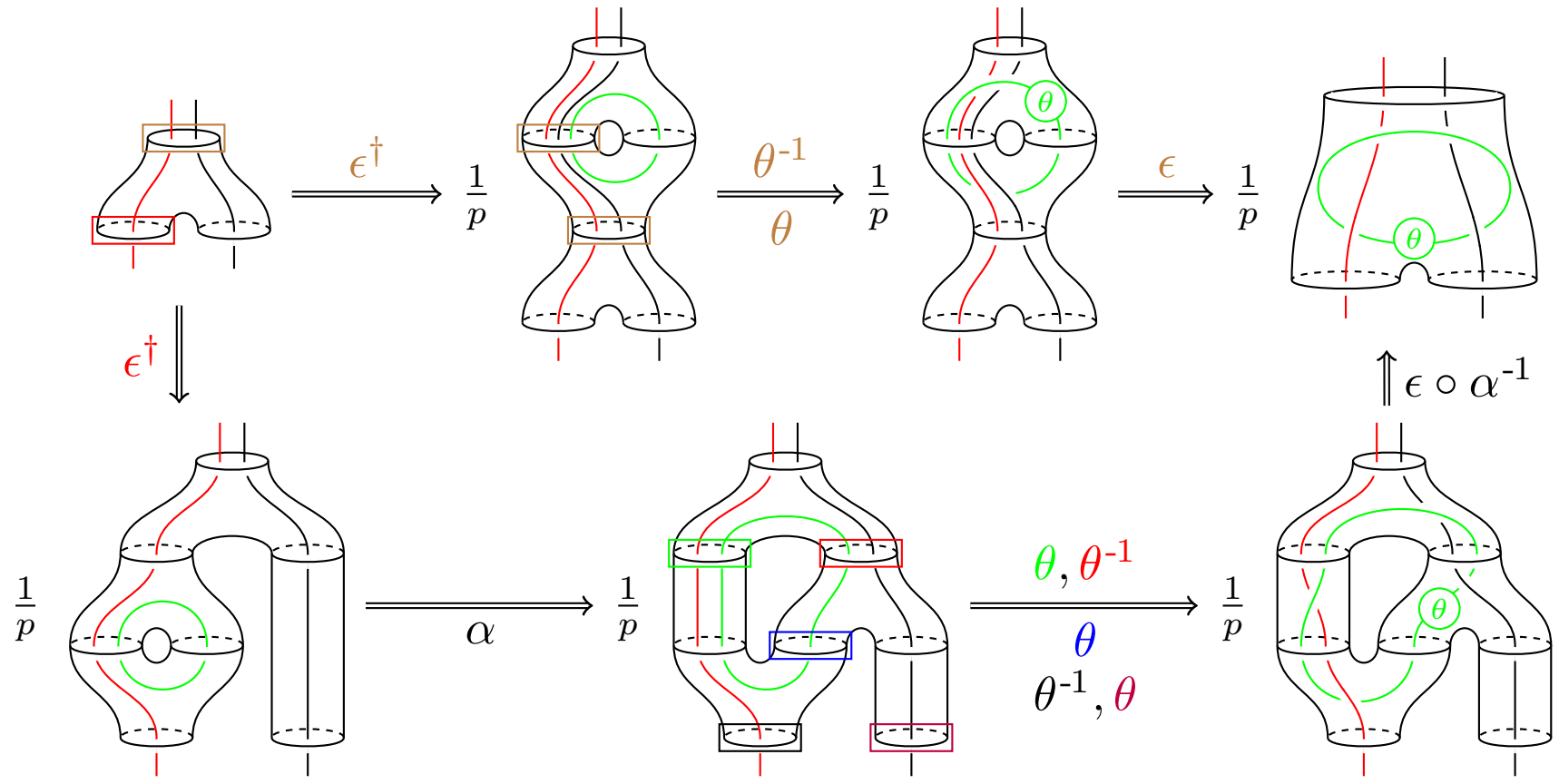

In this form, it states that 'twisting an unlabelled loop around two strands lying on its left all at once is the same as inserting the loop in the middle, and then consecutively twisting around the left and then the right strand'. Note that (57) does not prove the relation (55) inside $\mathbf{F}(\mathcal{M})$; rather it proves that (55) is true in any *linear representation* of $\mathbf{F}(\mathcal{M})$. It is known that there are equations which hold in every linear representation of $\mathbf{F}(\mathcal{M})$ but which do not hold in $\mathbf{F}(\mathcal{M})$; this follows from a result of Funar [13].

Nevertheless, we would like to mimic the internal string diagram calculation (57) inside $\mathbf{F}(\mathcal{M})$. To do this, we will need to use Lemma 14 to write $\epsilon^\dagger$ in the form 'create a torus, then attach it', and rewrite the twists in terms of elementary braid operations.

We begin with the clockwise direction of (55).

**Lemma 51.** *The clockwise direction of* (55) *can be written as follows:*

(58)

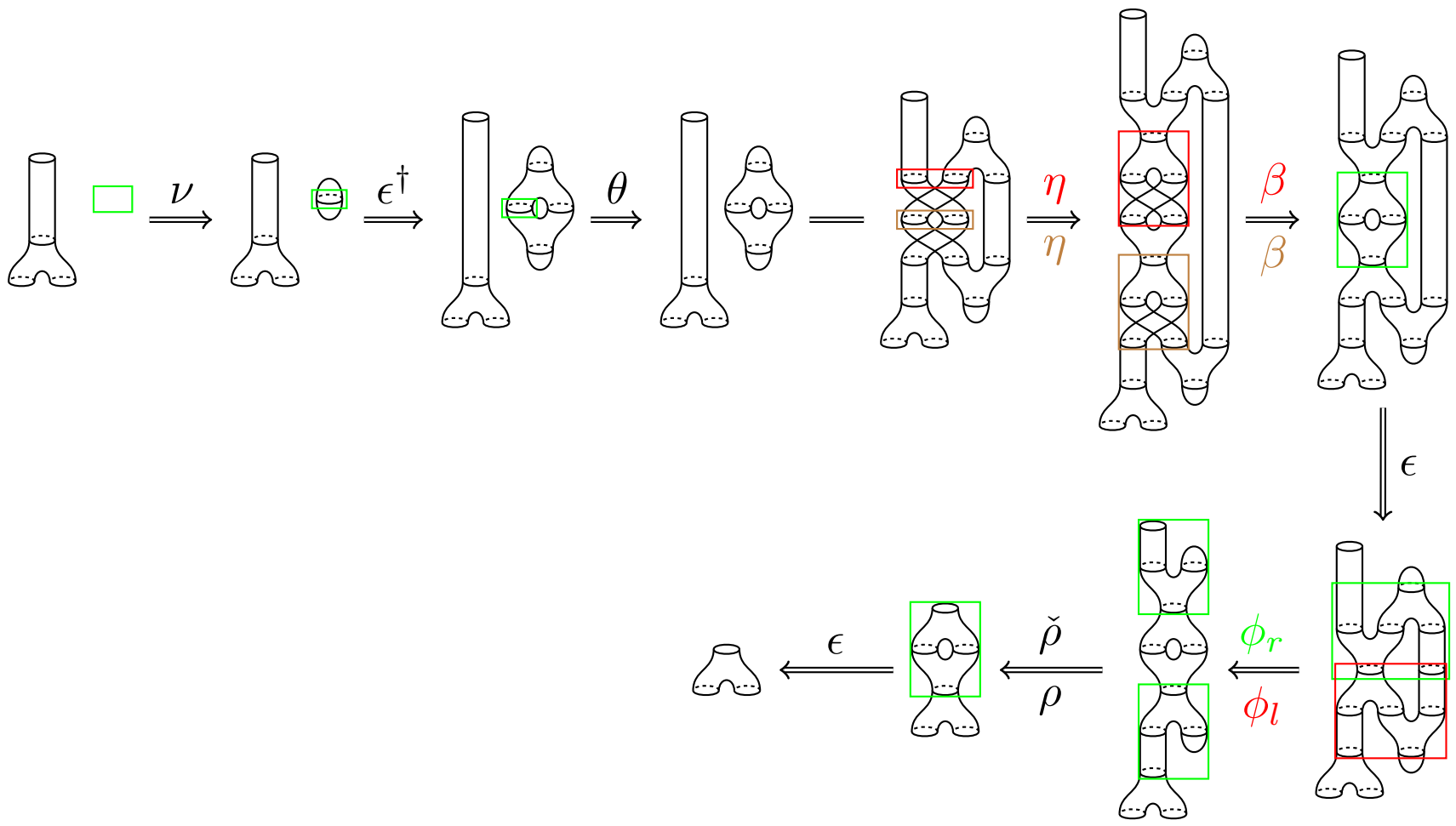

*Proof.* Starting from the expression above, replace the middle segment using the commutative diagram

(59)
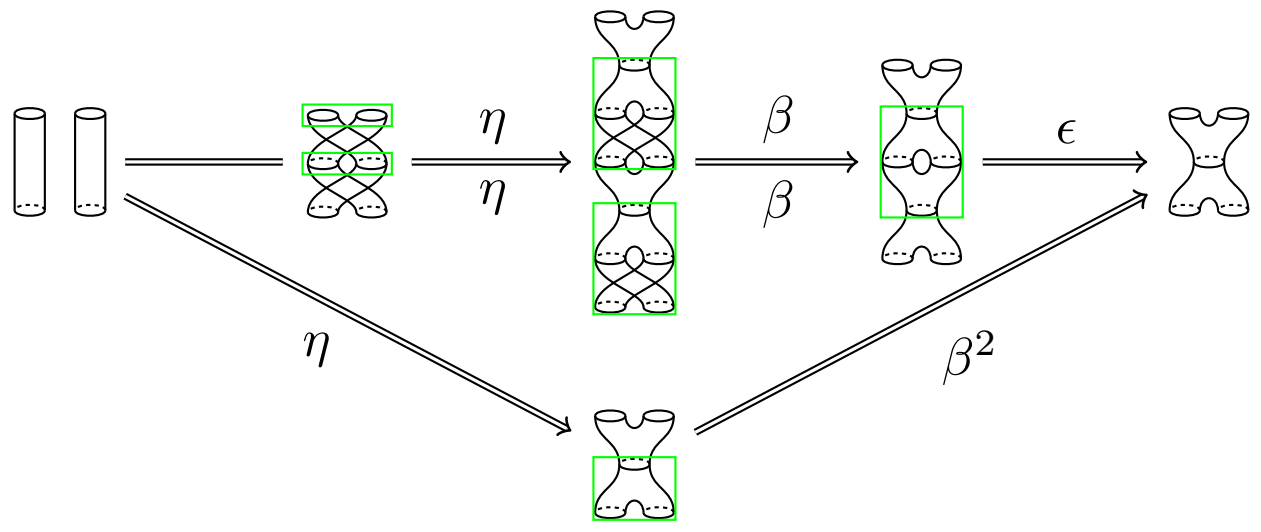

which follows from the $\eta$–$\epsilon$ adjunction equations. Then replace the $\beta^2$ with a $\theta$ and two $\theta^{-1}$s, one of which loops round and cancels with the initial $\theta$. The result is essentially the clockwise composite of (55), after expanding $\epsilon^\dagger$ using Lemma 14. □

**Lemma 52.** *The counterclockwise direction of* (55) *can be written as follows:*

(60)
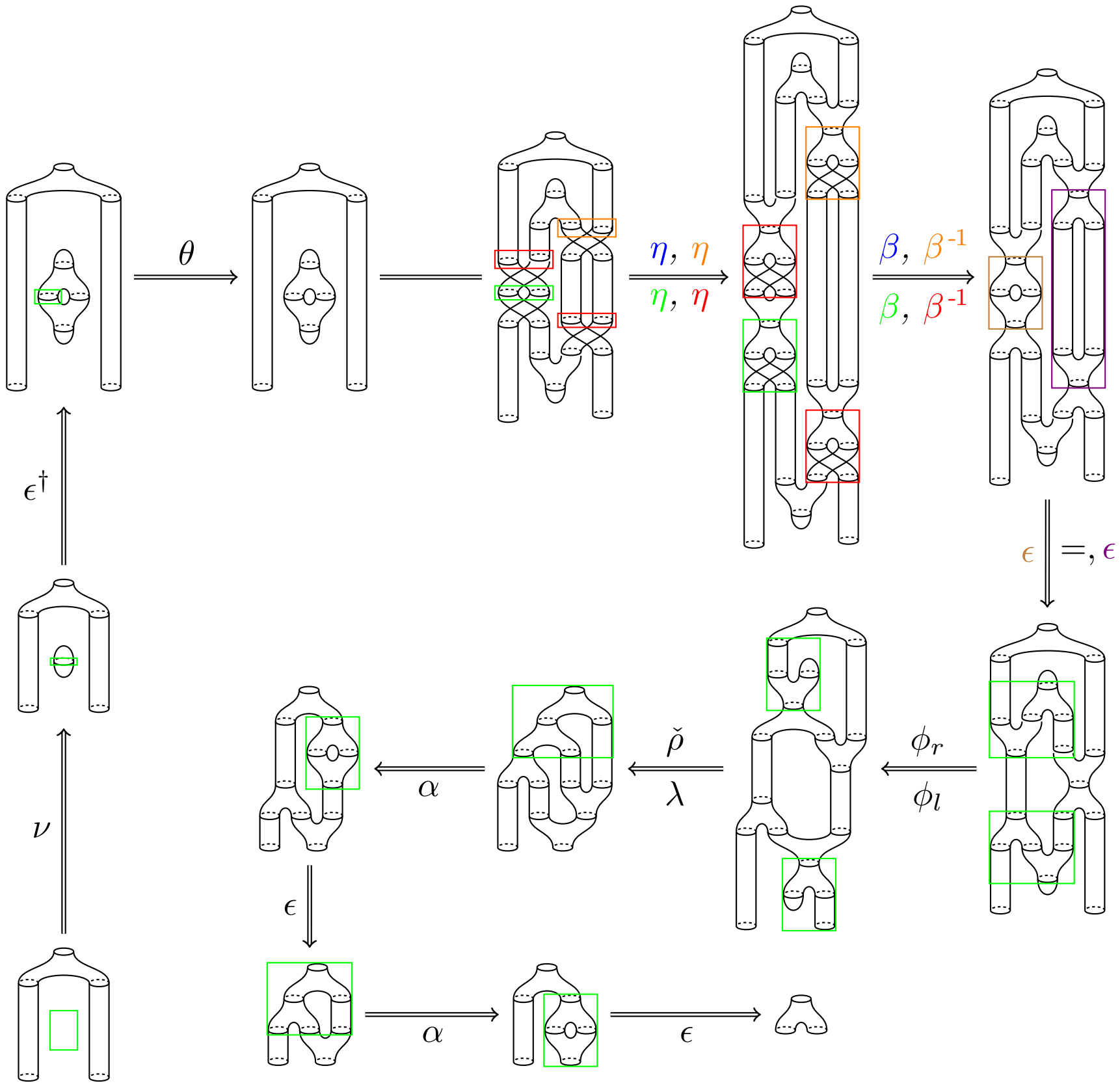

*Proof.* Begin with the above composite. Rearrange the top right segment beginning with the $\eta$ terms and ending with the pair of (brown and black) $\epsilon$ moves as follows. Let the middle pair of $\eta, \eta$ and $\beta, \beta$ moves (in red and green) act first, followed by the left (brown) $\epsilon$. Call this composite $X$. Then let the outer pair of $\eta$, $\eta$ and $\beta^{-1}$, $\beta^{-1}$ moves (in orange and red) act next, followed by the right (violet) $\epsilon$. Call this composite $Y$. Replace $X$ with a single $\eta$, $\beta^2$ move using equation (59), and similarly replace $Y$ with a single $\eta$, $\beta^{-2}$ move. Move the $\theta$ at the top left of the diagram forward so that it acts at the same time as the $\beta^2$ move.

Now, the initial segment of the above diagram takes the form 'create a torus, then attach it to the left leg of a pair of pants'. Use Lemma 14 to replace this with an $\epsilon^\dagger$ on the left leg of the pair of pants, followed by $\phi_r^{-1} \circ \check{\rho}^{-1} \circ \phi_l^{-1} \circ \rho^{-1}$. Cancel the $\phi_r^{-1}$ and $\check{\rho}^{-1}$ with the $\phi_r$ and $\check{\rho}$. Let the $\eta$ from the $Y$ segment cancel with the $\epsilon$ forming the third-last move in (60). We arrive at the following:

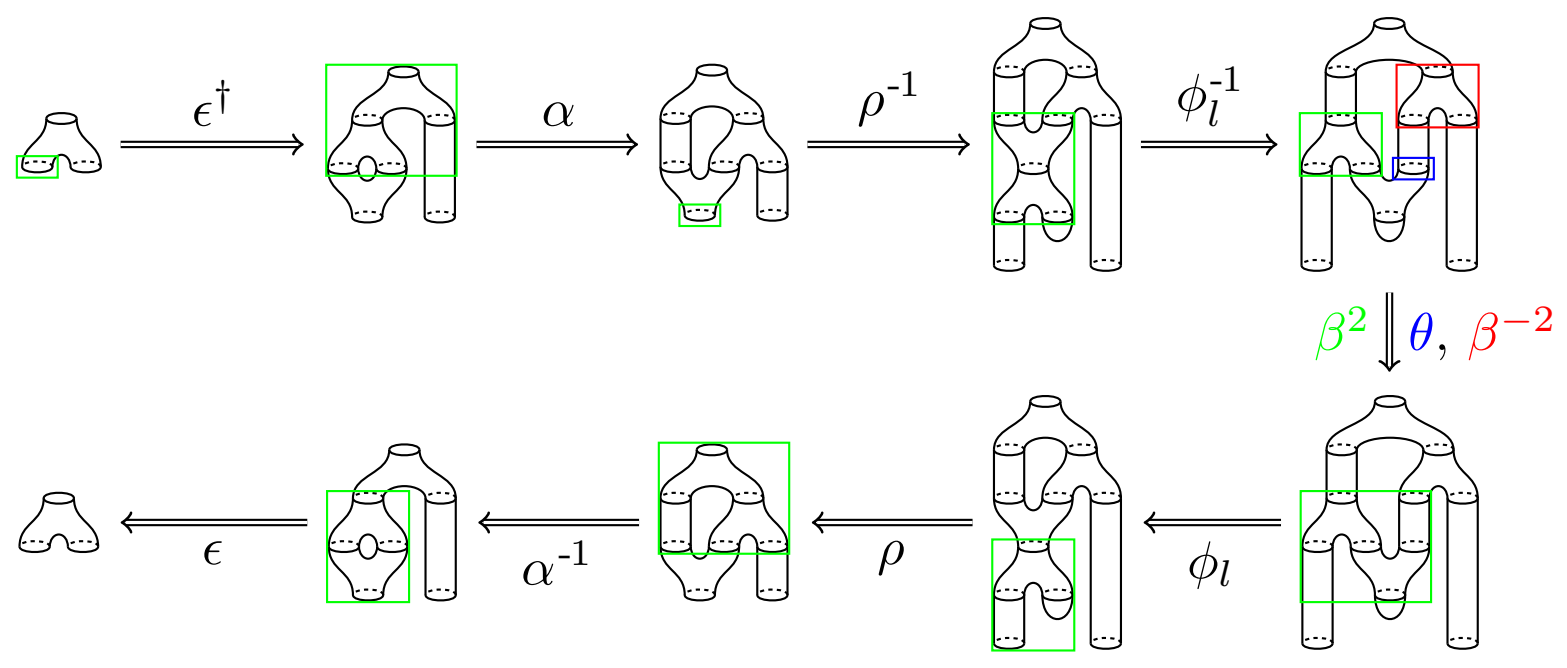

After replacing the $\beta^2$ and $\beta^{-2}$ with twists, we can cancel the $\phi_l$ and $\rho$ terms and we arrive at the counterclockwise composite in (55). □

To complete the proof that (58) equals (60), we will need to 'pull the torus together with its braiding moves through the pair of pants'. We will need the following three lemmas.

**Lemma 53.** *In* $\mathbf{F}(\mathcal{M})$*, the following diagram commutes:*

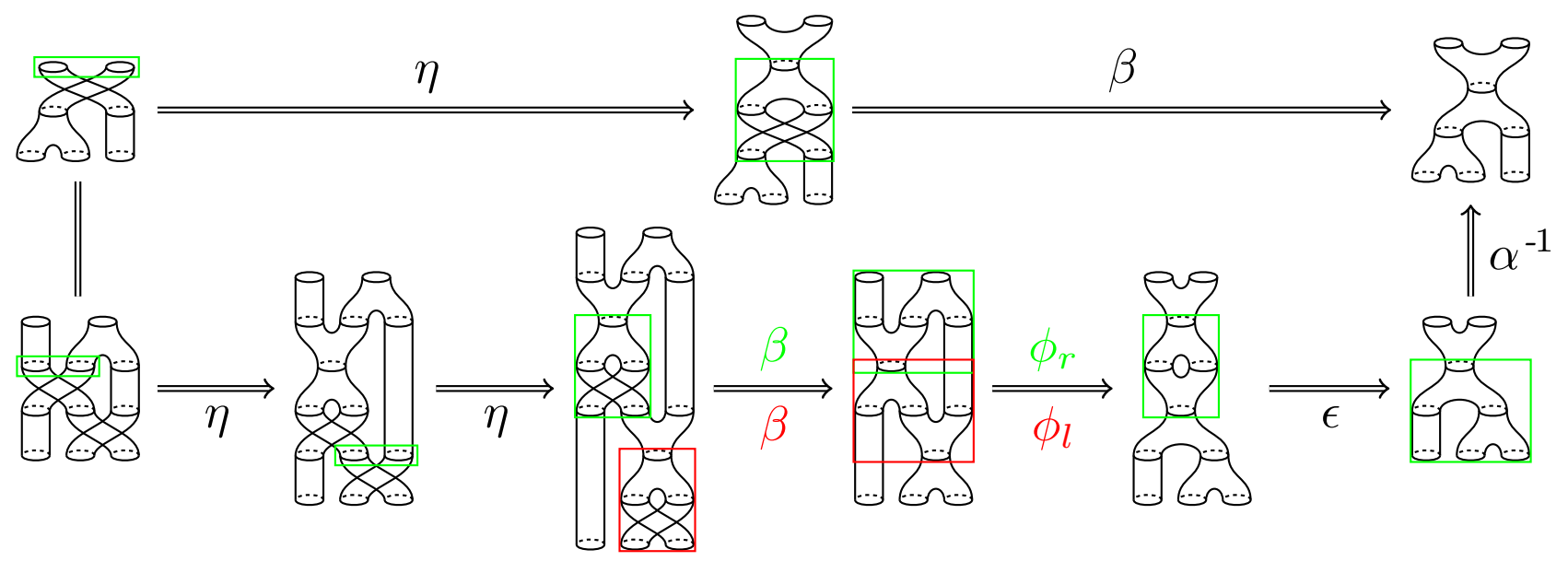

*Proof.* Expand out the Frobeniusators, use the hexagon equation (7) once, and apply the $\eta$–$\epsilon$ adjunction equation twice. □

**Lemma 54.** *In* $\mathbf{F}(\mathcal{M})$*, the following diagram commutes:*

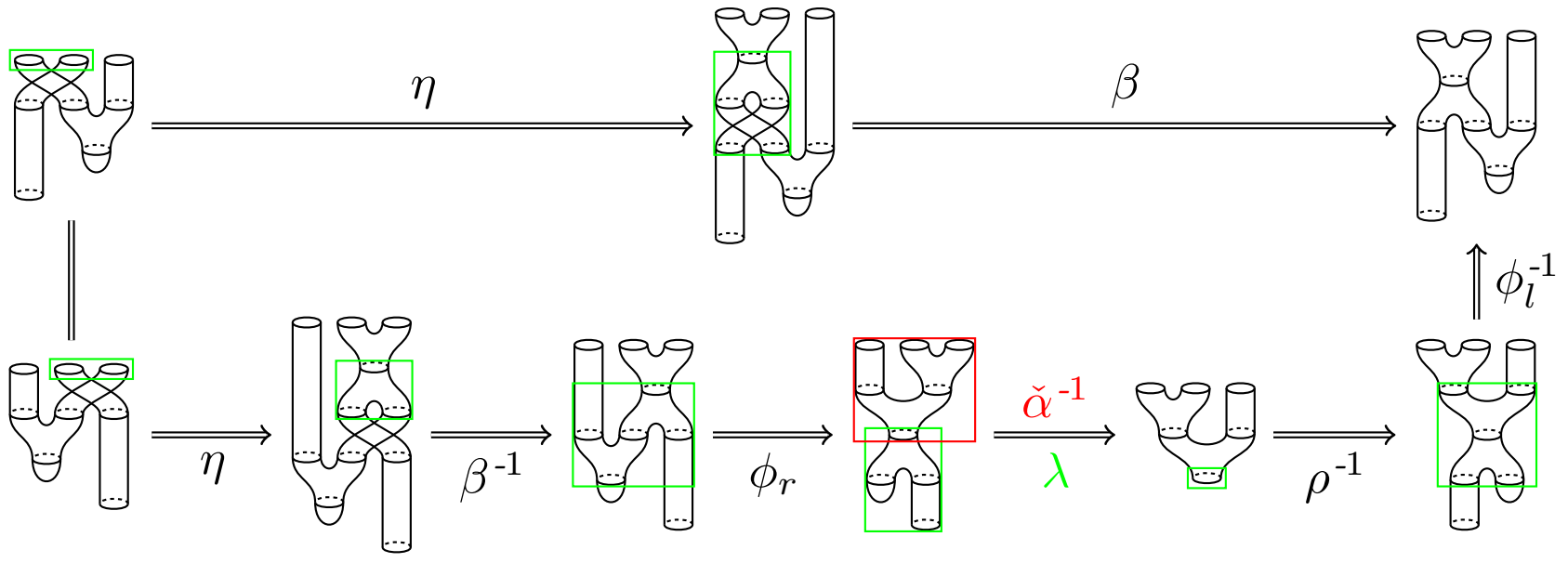

*Proof.* Consider the following diagram:

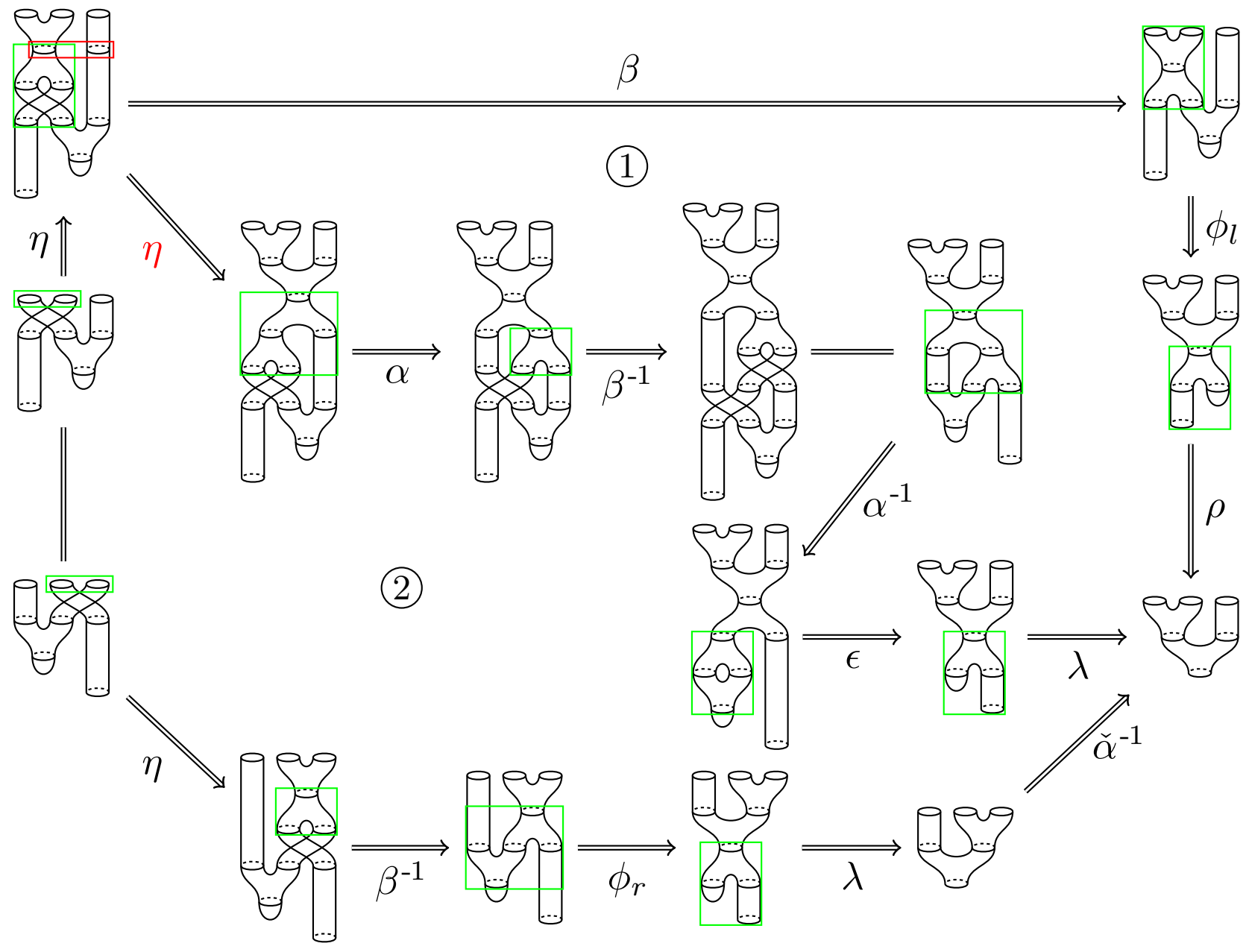

Diagram ① follows from expanding $\phi_l$ and the hexagon equation. Diagram ② follows from reversing the direction of $\check{\alpha}^{-1}$ to become a $\check{\alpha}$, and then moving it back to just after the (red) $\eta$. This allows the initial $\eta$ to be moved to the right, and the diagram follows. □

**Lemma 55.** *In* $\mathbf{F}(\mathcal{M})$*, the following diagram commutes:*

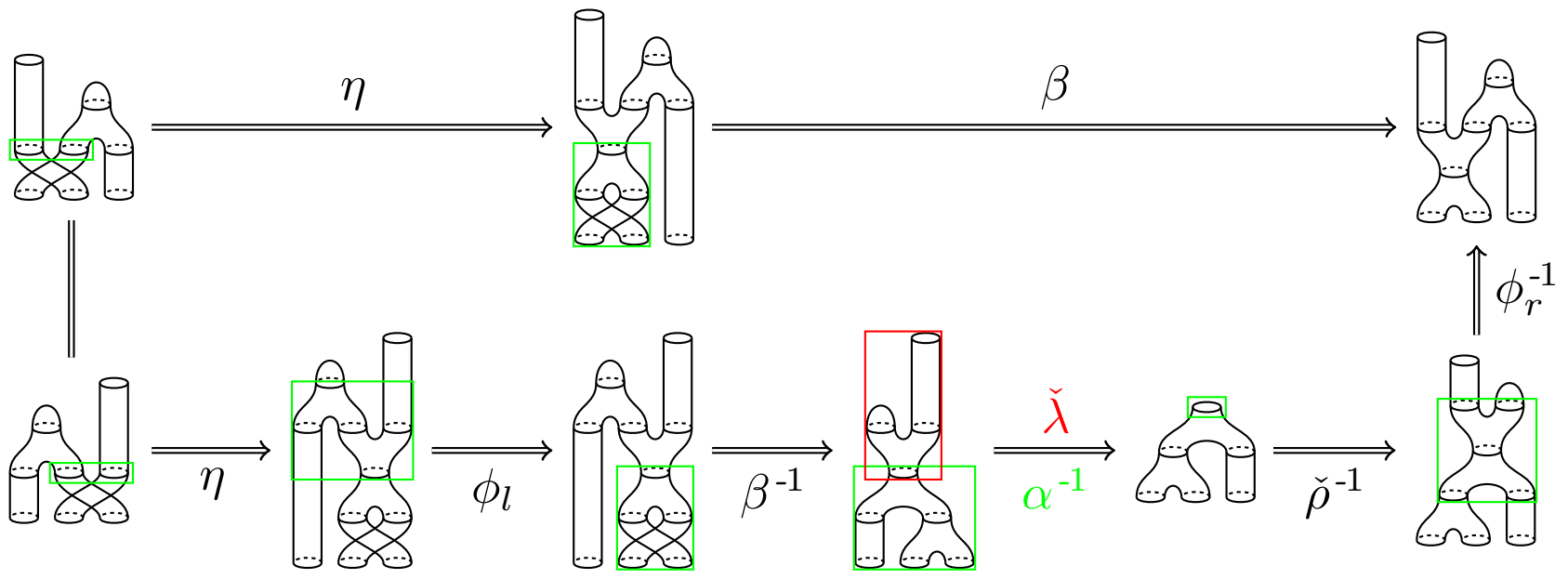

*Proof.* Reverse the direction of $\beta$ and move the resulting $\beta^{-1}$ round; then use the hexagon equation to replace $\beta^{-1} \circ \alpha^{-1} \circ \beta^{-1}$ with $\alpha^{-1} \circ \beta^{-1} \circ \alpha^{-1}$. Write the resultant $\alpha^{-1} \circ \phi_l \circ \eta$ segment as $\eta$. Now use which follows from expanding $\check{\lambda}$, the unit equation for the braiding (which follows from the hexagon equation), and the $\eta$–$\epsilon$ adjunction equation. □

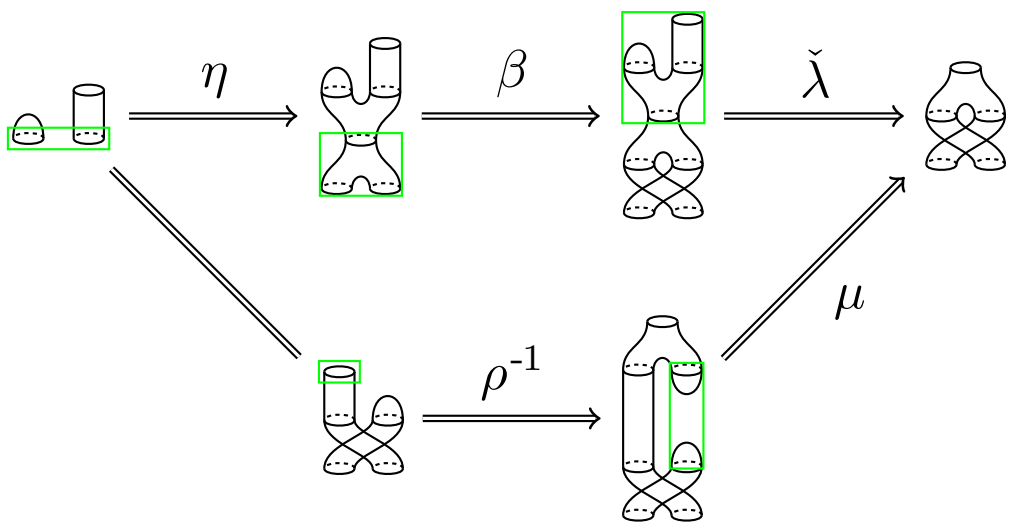

We can now finally complete the proof of Theorem 49.

**Lemma 56.** *Diagram* (55) *commutes.*

*Proof.* We need to show that (58) equals (52). Start with (58). Replace the (black) $\beta \circ \eta$ segment using Lemma 53, and then the (red) $\beta \circ \eta$ segment using the version of Lemma 53 rotated about the $z$ axis. Now use Lemma 54 and Lemma 55. The result, after some rearrangement, is (52). □

## 6. From the 2-Morse presentation to the modular presentation

In this section we prove there is a symmetric monoidal functor $K\colon \mathbf{F}(\mathcal{F}) \to \mathbf{F}(\mathcal{O})$. The 2-Morse presentation $\mathcal{F}$ is recorded in Appendix A.

**Defining the functor.** The anomaly-free modular presentation $\mathcal{O}$ has the same generating objects and generating 1-morphisms as $\mathcal{F}$, but only a subset of the generating 2-morphisms of $\mathcal{F}$. The functor $K$ will send those generators in $\mathcal{F}$ that are absent in $\mathcal{O}$ to composites of generators in $\mathbf{F}(\mathcal{O})$.

**Definition 57.** We define $K\colon \mathbf{F}(\mathcal{F}) \to \mathbf{F}(\mathcal{O})$ in terms of the generating objects, 1-morphisms and 2-morphisms of $\mathcal{F}$ as follows:

- On the generating object:

$$K(\circ) = \circ.$$

- On the generating 1-morphisms:

$K($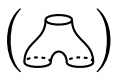$) =$ 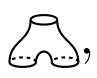, $K($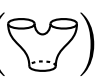$) =$ 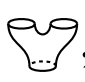, $K($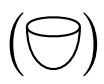$) =$ 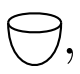, $K($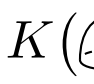$) =$ 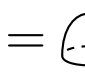.

- On the generating 2-morphisms:
  - (i) $K(\alpha) = \alpha$, and so on for all the generating 2-morphisms of $\mathcal{O}$ which are also generating 2-morphisms of $\mathcal{F}$;
  - (ii) $K(\check{\alpha})$, $K(\check{\rho})$, $K(\check{\lambda})$, and $K(\check{\beta})$ via the formulas (24)-(27) respectively;
  - (iii) $K(\phi_l)$ and $K(\phi_r)$ via the formulas (10) and (12) respectively;
  - (iv) $K(\mathrm{II})$, $K(\mathrm{III})$, and $K(\mathrm{III}')$ via the formulas (33), (35) and (36) respectively;
  - (v) for each of the 2-morphism generators $\omega$ in (i)–(iv), we define $K(\omega^{-1}) := K(\omega)^{\dagger}$.

In the remainder of this section we establish that $K$ is a well-defined functor by showing that it preserves the relations in $\mathcal{F}$.

6.1. **Proof of the invertibility relations.** We need to show that $K(\omega)^{\dagger} = K(\omega)^{-1}$ for each generator $\omega$ listed in Definition 57 (ii)–(iv).

We start with $\check{\alpha}$. If we take the dagger of the definition (24) of $K(\check{\alpha})$, we arrive at ${}^*\check{\alpha}$. Hence the relation $K(\check{\alpha})^{\dagger} = K(\check{\alpha})^{-1}$ is equivalent to the relation ${}^*\alpha = \alpha^*$, which follows from Lemma 19. Similarly, the invertibility relations for $K(\check{\rho})$, $K(\check{\lambda})$, $K(\check{\beta})$, $K(\phi_l)$, and $K(\phi_r)$ follow from Lemma 16, Lemma 20, and Lemma 18.

The invertibility relation for II follows from the left Rigidity relation, while the invertibility relations for III and III$'$ follow from Proposition 25.

6.2. **Proof of the Morse relations.** These are precisely the Biadjoint relations (15) and (16) together with their daggers and their rotated versions about the $x$ axis.

6.3. **Proof of the cusp flip relations.** The relation (CF) follows from naturality and the $\mu$–$\nu$ adjunction equation, after inserting $\rho^*$ for $\check{\rho}^{-1}$. The dagger of (CF) follows similarly except one uses ${}^*\rho$ for $\check{\rho}$. The six other versions of (CF) follow similarly.

6.4. **Proof of the swallowtail relations.**

*ST1.* The relation $(\mathrm{ST1})^x$ follows from coherence of the monoidal category axioms.[5] The relation (ST1) itself follows from taking duals of both sides. The other two versions follow similarly.

*ST2.* The relation (ST2) is precisely the unit equation (6). The relation $(\mathrm{ST2})^x$ follows from taking duals of both sides.

*ST3.* The relation (ST3) follows from inserting the decomposition of II as $\phi_l \circ \theta \circ \phi_l^{-1}$ and then moving the $\theta$ around the bend using the Ribbon relation (14). The rotated version $(\mathrm{ST3})^z$ follows similarly from inserting the decomposition of II as $\phi_r \circ \theta \circ \phi_r^{-1}$. The other two versions follow from taking duals.

6.5. **Proof of the beak relations.**

[5] This was actually one of the original axioms of a monoidal category, but Kelly showed [26] it is a consequence of the pentagon equation and (ST2).

*B1.* To prove $(\mathrm{B1})^x$, take duals of both sides and rearrange slightly. We must prove that the composite along the top row below is equal to $\rho$:

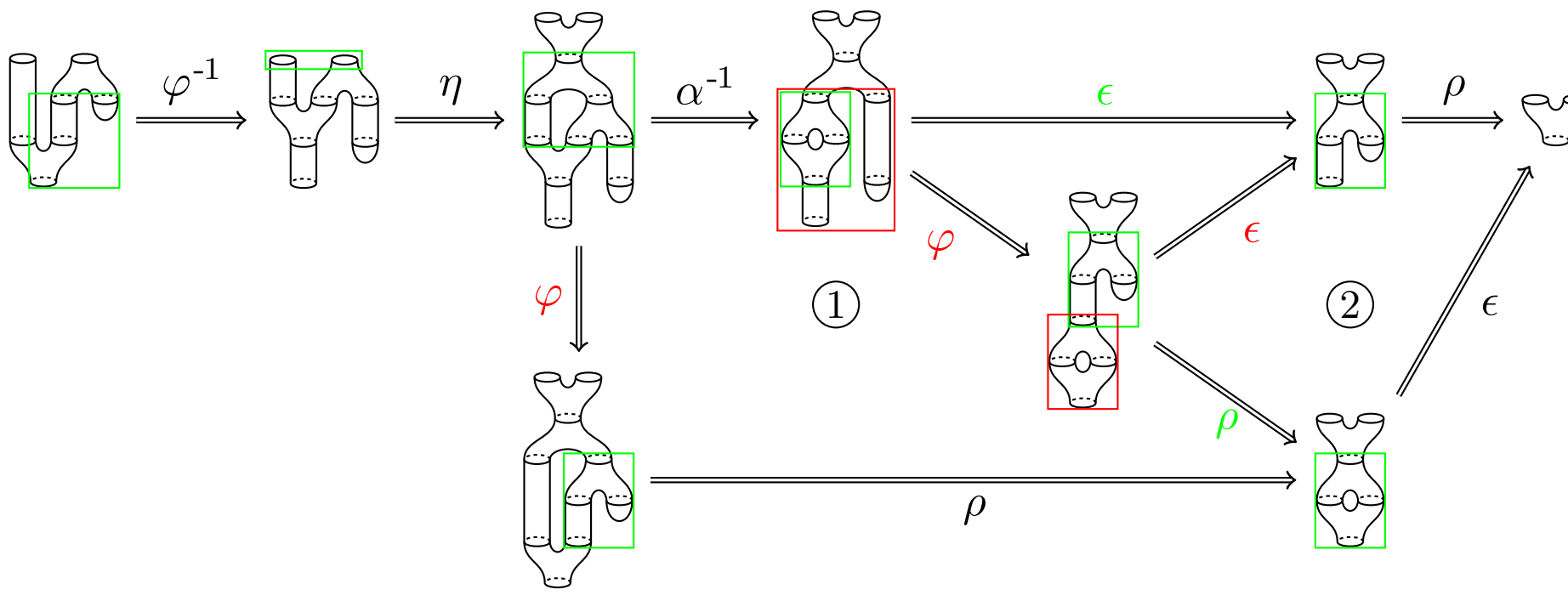

Diagram ① follows from monoidal coherence, and diagram ② follows from naturality. In the counterclockwise composite, move the $\eta$ past the $\rho$ and cancel with the $\epsilon$, and we are done. The relation (B1) follows from taking duals of both sides. The other two versions are similar.

*B2.* This version of (B2) commutes by monoidal coherence, and similarly for $(\mathrm{B2})^z$. The other two versions follow by taking duals of both sides.

6.6. **Proof of the Morse flip relations.**

*MF1.* The relation (MF1) follows from the definition of $\phi_l$ and the adjunction equations, and similarly $(\mathrm{MF1})^z$ follows from the definition of $\phi_r$ and the adjunction equations. Similarly for the other six versions.

*MF2.* We claim that the relation (MF2) is equivalent to the definition (35) of III. Indeed, insert the identity in the form of $\epsilon^\dagger$–$\eta^\dagger$ at the beginning of the counterclockwise composite of (MF2). We obtain

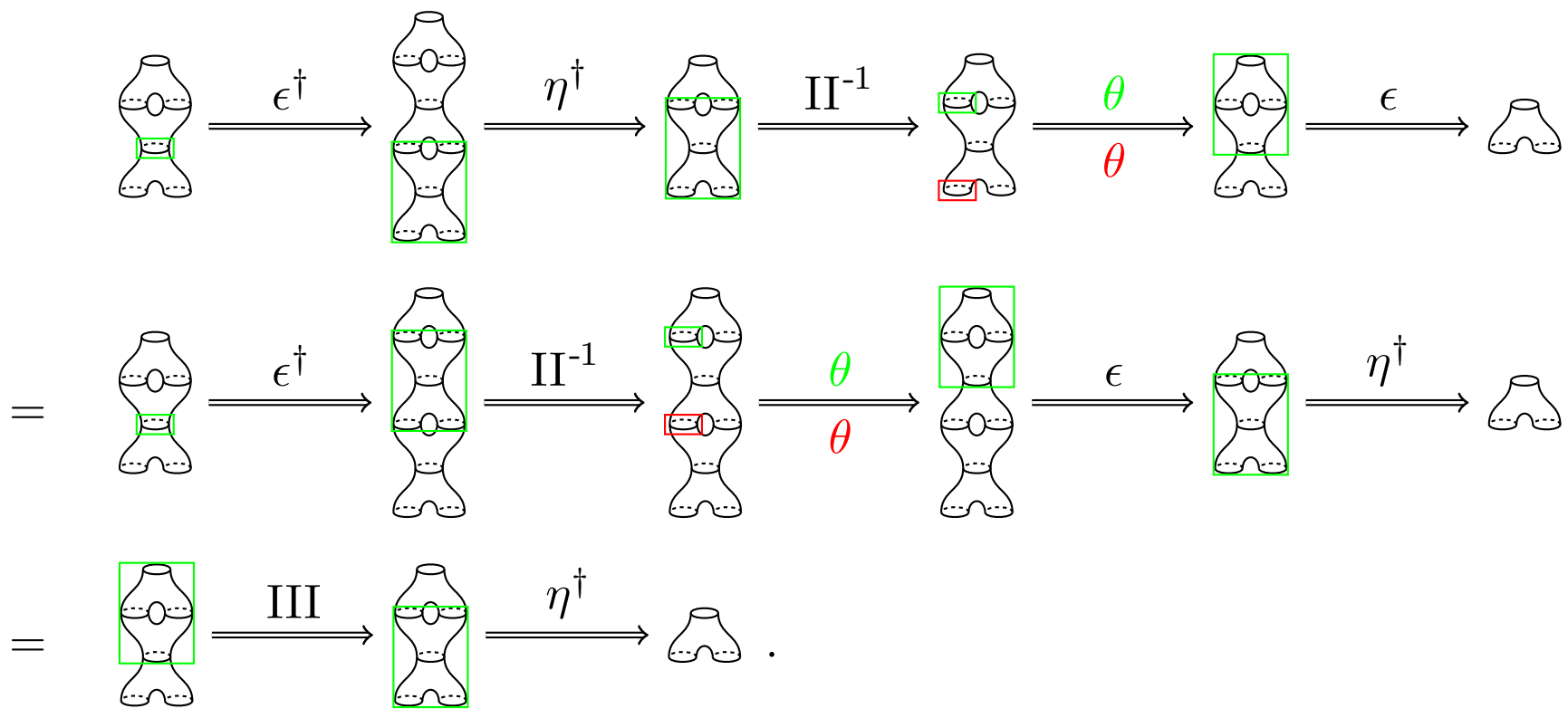

The rotated version $(\mathrm{MF2})^x$ follows from $\mathrm{III}^* = \mathrm{III}$, which was Corollary 30. The other two versions follow similarly.

6.7. **Proof of the triple point relations.**

*T1.* The relation (T1) is the Pentagon relation (5). Its rotated version $(\mathrm{T1})^x$ follows from taking duals of both sides.

*T2.* Insert the definition (35) of III in terms of $A$; the first and last $\theta$ terms immediately cancel and so we can replace each III in (T2) with an $A$. Then expand $A$ according to Lemma 22:

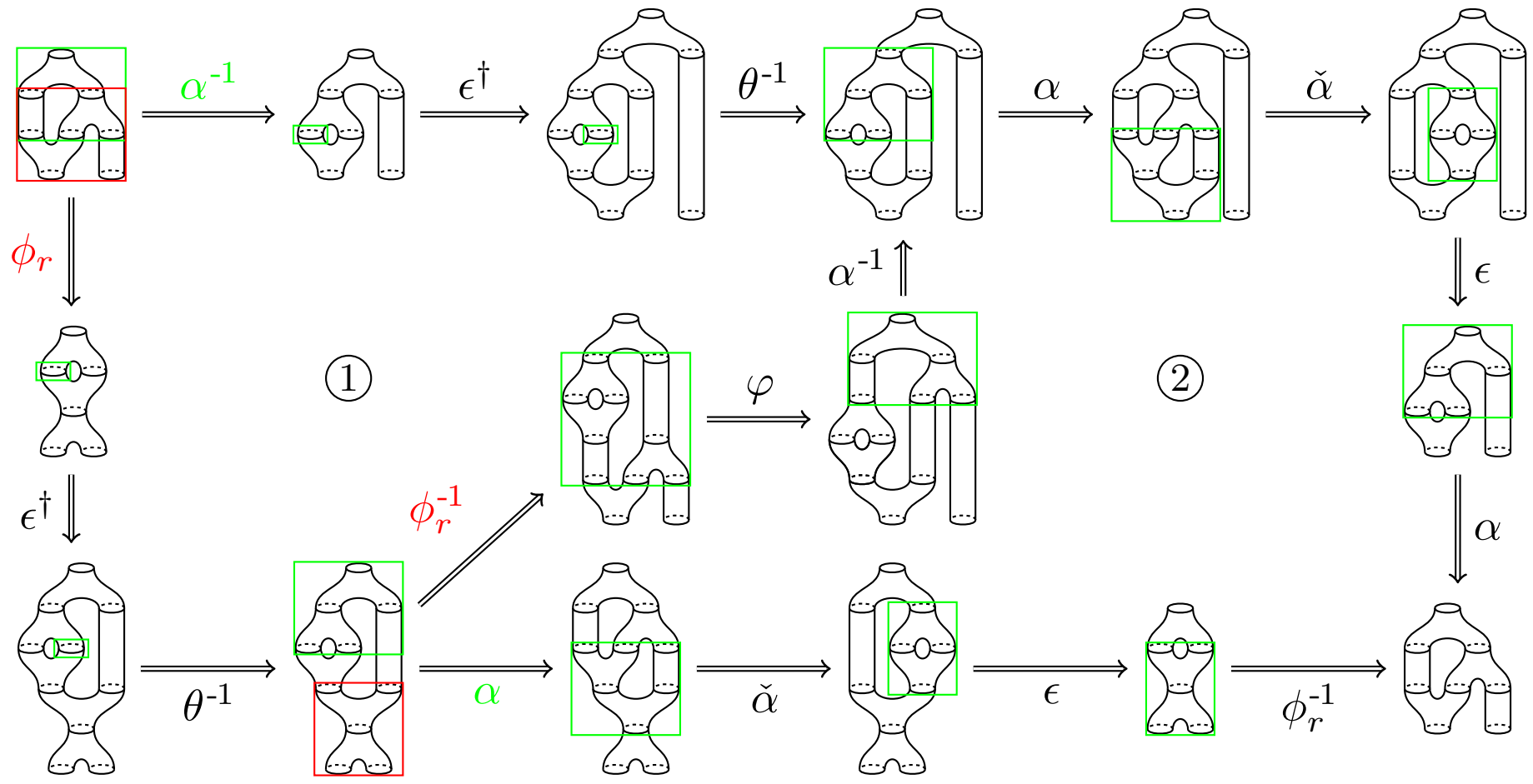

Diagram ① commutes by naturality and diagram ② commutes by expanding out $\breve{\alpha}$ as $(\alpha^{-1})^*$, and then using the pentagon equation and the $\eta$–$\epsilon$ adjunction equations. This proves (T2). The proof of $(\mathrm{T2})^z$ is essentially the same, except we use the expression for $A$ given in Corollary 29 instead. Finally, the two other rotated versions of (T2) follow by taking duals on both sides and using $\mathrm{III}^* = \mathrm{III}$ and $\mathrm{III}' = (\mathrm{III}')^*$.

*T3.* To prove (T3), write $\check{\Gamma} := \breve{\alpha}^{-1} \circ \phi_l^{-1}$. Consider the composite $\mathrm{III}^{-1} \circ \theta^{-1} \circ \theta \circ \mathrm{II} \circ \mathrm{III}$ at the bottom of (T3). We can eliminate the IIIs as follows:

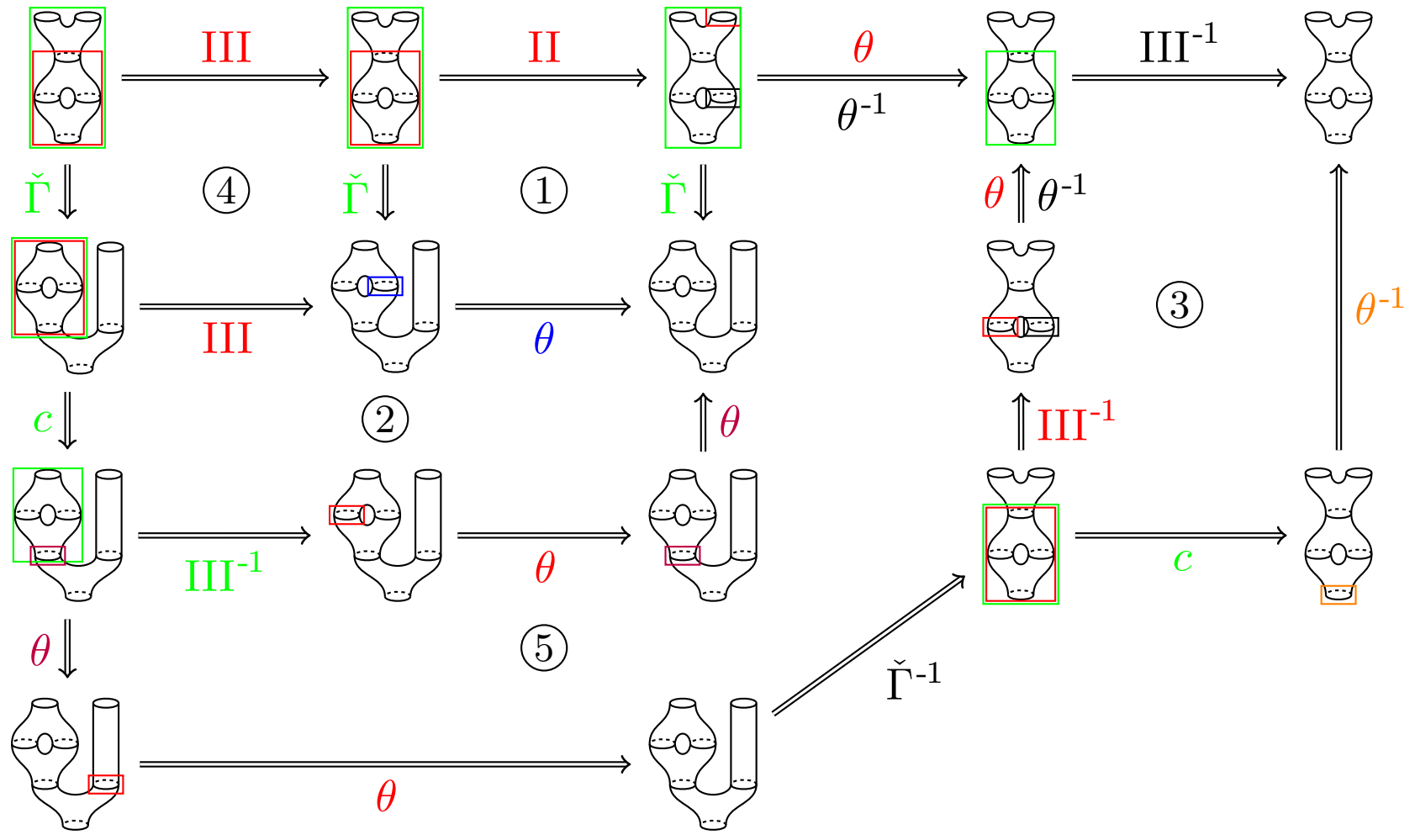

Diagram ① follows from expanding II. Diagrams ② and ③ are applications of Proposition 37, while diagrams ④ and ⑤ are (T2). Having eliminated the III 2-morphisms it is then not hard to show, using $(\text{FFG1})^z$ and the Balanced relations, that the counterclockwise composite above equals the remaining composite in (T3).

The rotated version $(\text{T3})^x$ can be proved in a similar way, using $(47)^x$. The other two rotated versions of (T3) follow from taking duals of both sides.

*T4.* Similarly to the proof of (T3), write $\Gamma := \alpha^{-1} \circ \phi_r^{-1}$. We have the following commutative diagrams:

(61)

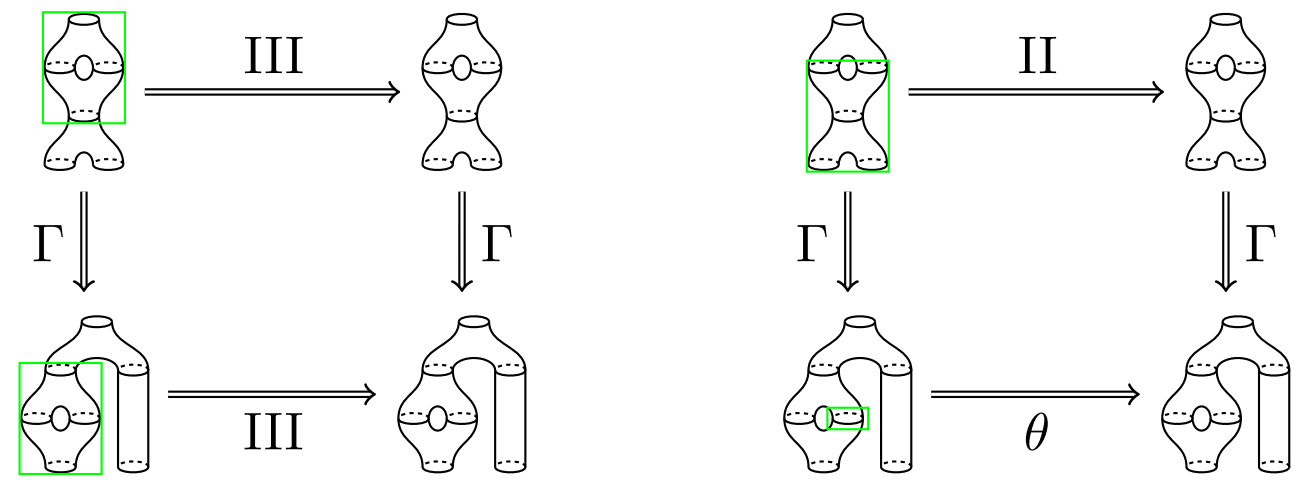

The first diagram is simply (54), and the second follows from expanding II as $\phi_r^{-1} \circ \theta \circ \phi_r$. Hence (T4) is equivalent to the following diagram on the twice-punctured torus (the one direction follows from adding a 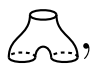, the other direction follows from capping off at the top and adding a 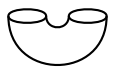 to the bottom):

(T4′)

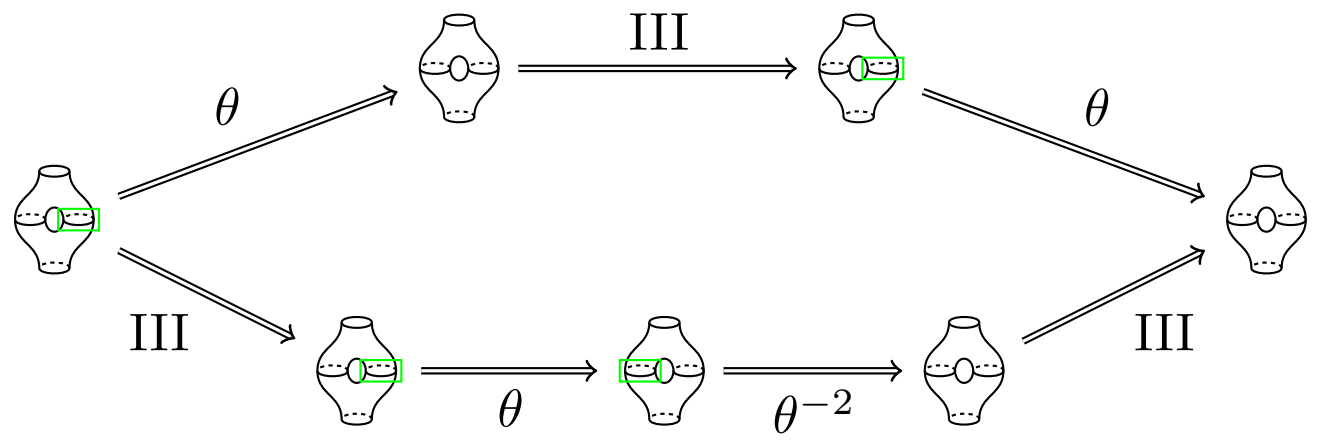

If we substitute in the definition of III in terms of $A$ from (35), equation (T4′) is precisely the 'anomaly-free' version of the anomalous relation (54), which was shown to follow from the other relations in Theorem 49. In the anomaly-free case, $x = \text{id}$ from Corollary 43 below, hence (T4′) follows. This establishes (T4). Since $(\text{T4})^x$ is also equivalent to (T4′) — by using $(61)^x$ — this establishes $(\text{T4})^x$ too. The other two rotated versions of (T4) are proved similarly, using the rotated version of Theorem 49, namely $(54)^z$.

*T5.* The relation (T5) follows from coherence for braided monoidal structures, since both sides produce the same underlying braid. Similarly for $(\text{T5})^z$. The other two rotated versions follow from taking duals.

*T6.* The relation (T6) follows from expanding the Frobeniusators, the pentagon equation and the adjunction equations between $\eta$ and $\epsilon$. Similarly for $(\text{T6})^z$. The other two rotated versions follow from taking duals.

*T7.* After expanding the Frobeniusators and using the adjunction equations, becomes equivalent to (T5). Similarly for the rotated versions.

*T8.* Follows from expanding the Frobeniusators, the pentagon equation and the adjunction equations for $\eta$ and $\epsilon$. Similarly for the rotated versions.

*T9.* The relation (T9) follows from expanding the Frobeniusators and II as $\phi_l \circ \theta \circ \phi_l^{-1}$. The rotated version $(\mathrm{T9})^z$ follows from expanding II as $\phi_l \circ \theta \circ \phi_l^{-1}$ using Proposition 31. The other two versions follow from taking duals on both sides.

*T10.* We can insert a 'bent arm' on the upper right hand circle of each node in order to express (T10) as an equivalent relation whose nodes are constructed solely from 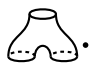. For instance, after attaching a bent arm the $\mathrm{II}^{-1}$ at the top left of (T10) can be rewritten as follows:

(62)

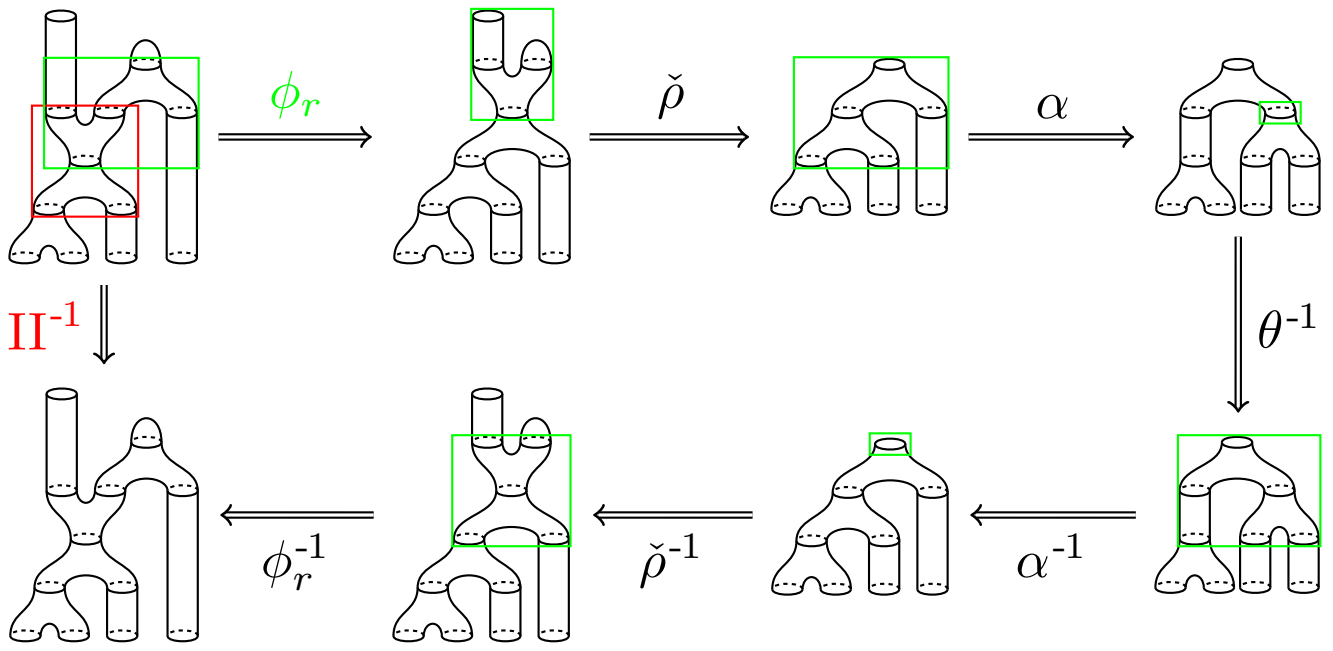

This diagram follows from expanding II in terms of $\phi_r$ using Proposition 31. Since the arms can be bent back due to invertibility of the Frobeniusator, this procedure indeed produces an equivalent relation. Applying this to all the 2-morphisms of (T10) we find that it is equivalent to the following diagram:

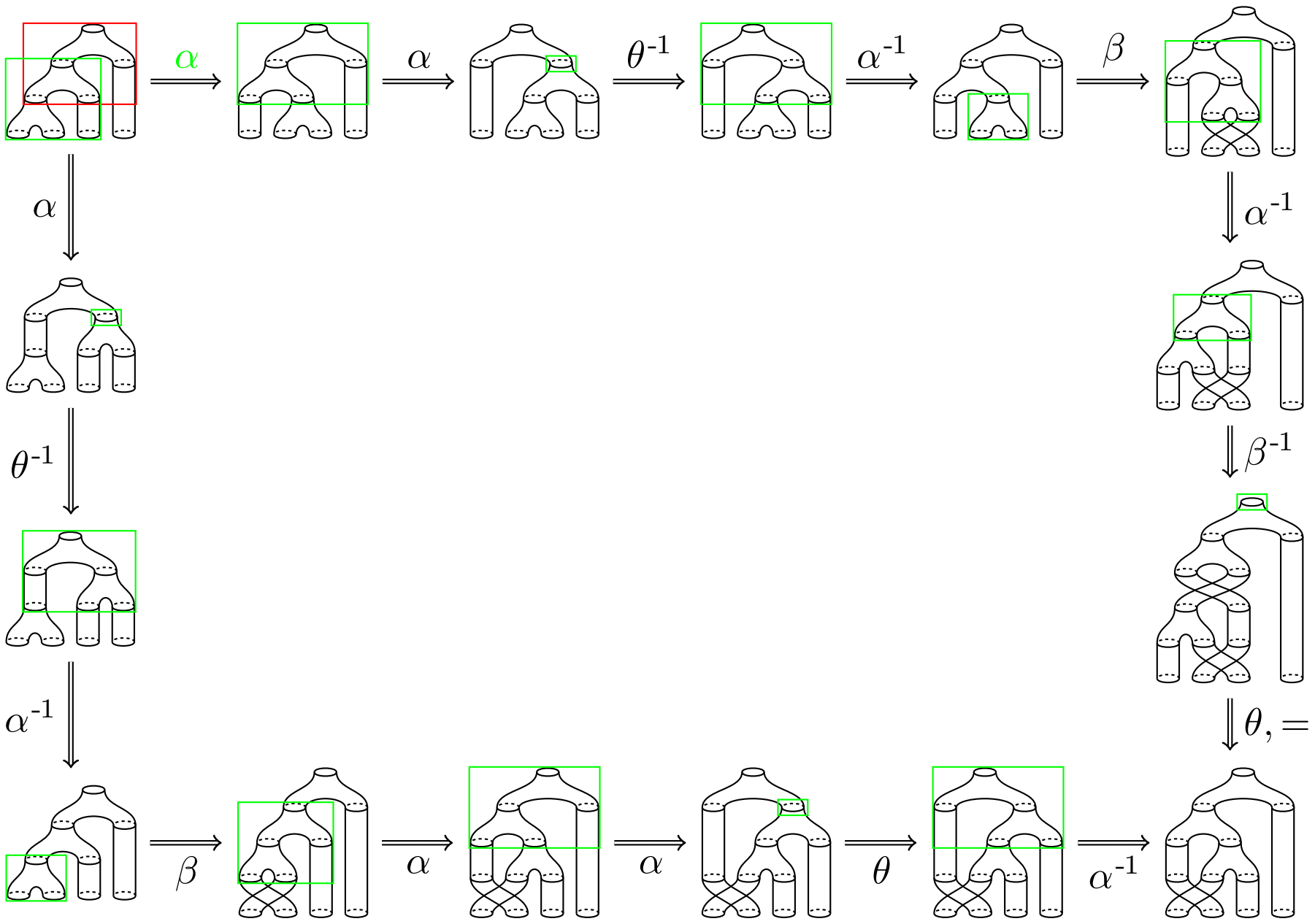

This diagram commutes by coherence for balanced tensor categories [24], since both sides produce the same underlying ribbon. The rotated version $(\mathrm{T10})^z$ follows similarly. The other two rotated versions follow from taking duals on both sides.

### 6.8. Proof of the fold-codim 2 gluing relations.

*FCT1, FCT2.* These are the Balanced relations (8) and (9). The rotated version follows from taking duals of both sides.

### 6.9. Proof of the fold-fold-glue relations.

*FFG1.* Start at the bottom right hand corner of (FFG1), and consider the clockwise composite. After inserting the definition of $\phi_r$, we obtain the counterclockwise composite below:

(63)

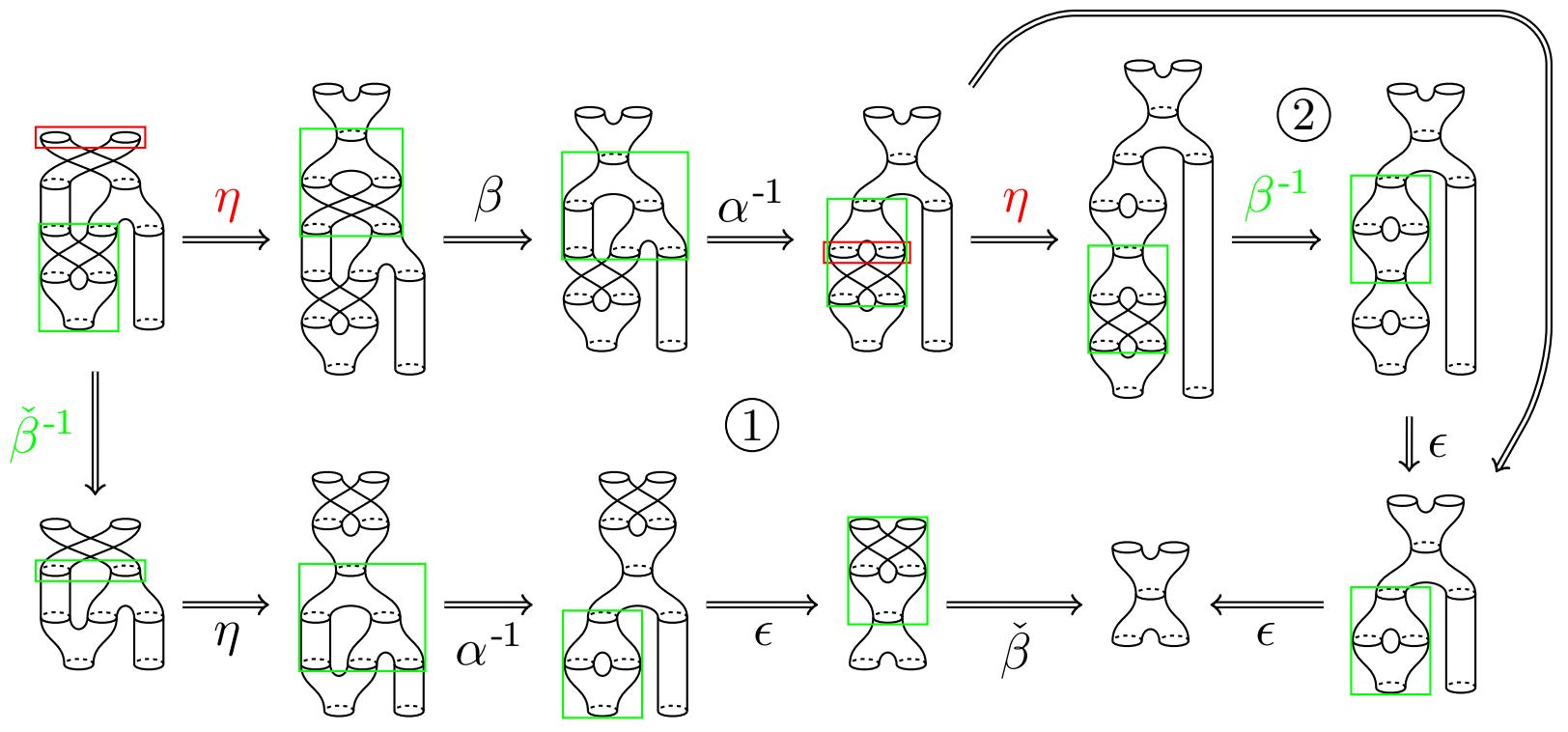

Diagram ① follows from naturality and the definition of $\check{\beta}$ as $\beta^*$. Diagram ② follows from naturality and the $\eta$–$\epsilon$ adjunction relation. On the other hand, after inserting the definition of $\phi_l$ and applying naturality, the counterclockwise composite of relation (FFG1) computes as follows:

(64)

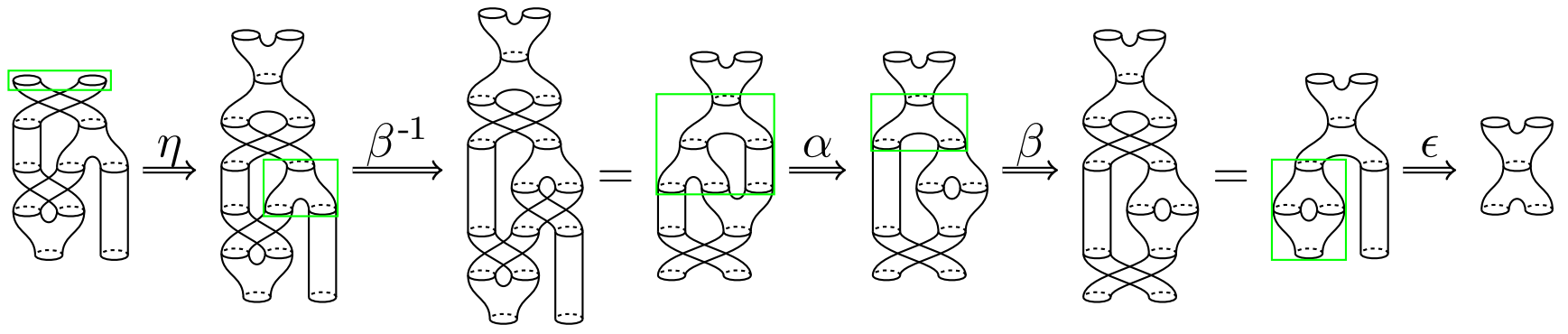

The $\beta^{-1}\alpha^{-1}\beta$ composite in (63) produces the same underlying braid as the $\beta\alpha\beta^{-1}$ composite in (64), and hence they are equal by coherence for braided tensor categories.

*FFG2.* The relation (FFG2) follows from (FFG1) by replacing the II in the clockwise composite with $\phi_r \circ \theta \circ \phi_l^{-1}$ and the II in the bottom composite with $\phi_l \circ \theta \circ \phi_l$.

*FFG3.* Write $\Gamma_l := \lambda \circ \phi_r \circ \alpha$ and $\Gamma_r := \rho \circ \phi_l \circ \alpha^{-1}$. Firstly, we will need the following diagram, which follows from (FFG1):

(65)

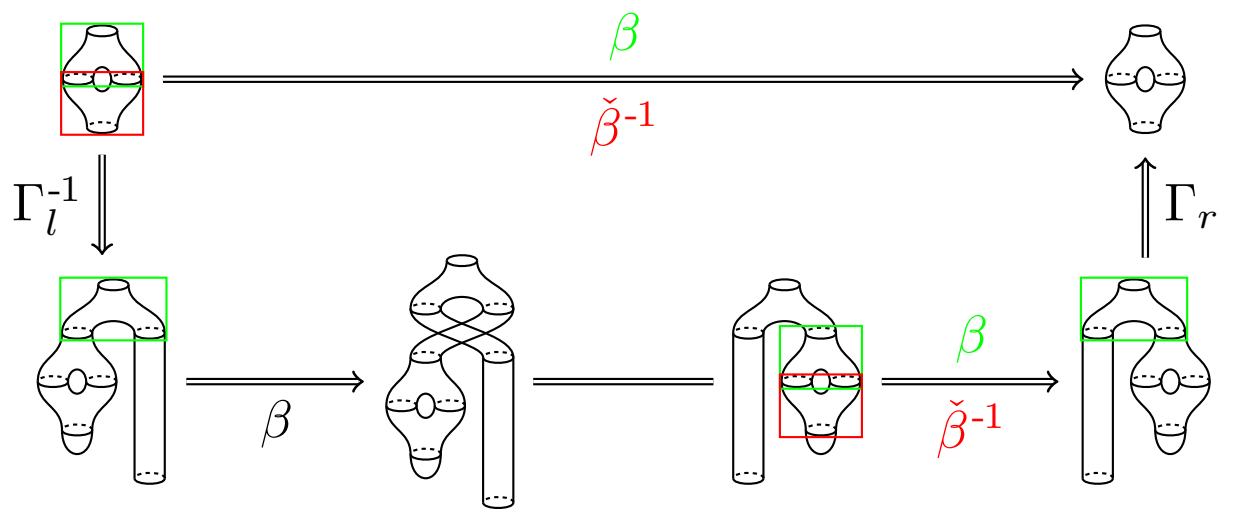

It suffices to prove (FFG3) on the once-punctured torus. To see this, write III and III′ in (FFG3) in terms of their actions on the once-punctured torus using (39) and $(39)^z$ (which can be proved in the same way as (39)). Also, recall from Remark 24 that III equals III′ on the once-punctured torus. Now apply (65) twice. We are left with establishing that III commutes with the hyperelliptic involution $c$ on the once-punctured torus:

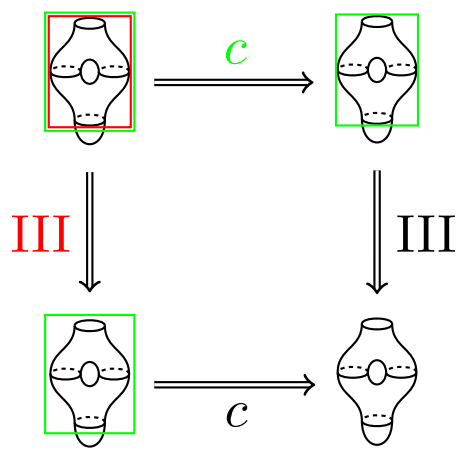

Corollary 38 implies both composites equal $\text{III}^3$, hence the diagram follows.

The proof for the version of (FFG3) rotated about the $z$-axis is the same except one uses the rotated version of (65).

*FFG4.* This is Charge Conjugation, Proposition 37.

### 6.10. **Proof of the cusp-glue relations.**

*CG.* This unit relation is known to follow from the hexagon equation (7) [24].

### 6.11. **Proof of the 3d Morse-glue relations.**

*3DMG.* Replace $\check{\beta}$ with ${}^*\beta$ using Lemma 20, then use naturality and the $\epsilon^\dagger$–$\eta^\dagger$-adjunction relation. Similarly for the dagger of this relation.

## 7. From the modular presentation to the 2-Morse presentation

In this section we establish a symmetric monoidal functor $L\colon \mathbf{F}(\mathcal{O}) \to \mathbf{F}(\mathcal{F})$, and prove that it provides an inverse for the functor $K\colon\ \mathbf{F}(\mathcal{F}) \to \mathbf{F}(\mathcal{O})$ established in Section 6.

The generators of $\mathcal{O}$ are a subset of the generators of $\mathcal{F}$. Hence the functor $L$ is defined in the natural way.

**Definition 58.** The symmetric monoidal functor $L\colon \mathbf{F}(\mathcal{O}) \to \mathbf{F}(\mathcal{F})$ sends each generator of $\mathcal{O}$ to the corresponding generator of $\mathcal{F}$.

To show that $L$ is well-defined, we need to check that $L$ respects the relations in $\mathcal{O}$. We will show this in Sections 7.2 to 7.9 respectively. We conclude the proof of Theorem 2 in Section 7.10.

7.1. **Decomposition of generators.** It is convenient to first perform a consistency check. Some of the 2-morphisms (such as III) are *defining generators* in $\mathcal{F}$, while they are defined as *composites* in $\mathbf{F}(\mathcal{O})$. We must check that these decompositions, which were *definitions* in $\mathbf{F}(\mathcal{O})$, are actually valid *equations* in $\mathbf{F}(\mathcal{F})$.

**Proposition 59.** *In $\mathbf{F}(\mathcal{F})$, the following decompositions of the generating 2-morphisms of $\mathcal{F}$ hold:*

(i) *the decomposition of $\phi_l$ and $\phi_r$ in terms of $\alpha$, $\eta$ and $\epsilon$ as expressed in* (10) *and* (12) *respectively, and the decompositions of $\phi_l^{-1}$ and $\phi_r^{-1}$ as expressed in* (11) *and* (13) *respectively;*
(ii) *the decomposition of $\check{\alpha}$, $\check{\rho}$, $\check{\lambda}$ and $\check{\beta}$ as the right duals of the corresponding 'unchecked' generators, expressed in* (24)–(27) *respectively;*
(iii) *the decomposition of* II *in terms of $\theta$ and $\phi_l$, expressed in* (33)*, and that* $\mathrm{II}^{-1}$ *equals the dagger of the right hand side of* (33)*;*
(iv) *the decomposition of* III *and* III′ *in terms of $\theta$, $\epsilon^\dagger$, $\epsilon$ and* II *as expressed in* (33)*,* (35) *and* (36) *respectively, and* $\mathrm{III}^{-1}$ *equals the dagger of the right hand side of* (35)*, and similarly for* III′*.*

*Proof.* (i) To prove the decomposition (10) of $\phi_l$, start with 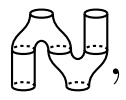, apply $\eta$ to the cylinders on the top two boundary circles, and then form an equation by postcomposing with either of the two composites from (MF1). Applying naturality and (MR2), the one composite becomes $\phi_l$ and the other becomes the right hand side of (10). A similar proof using $(\mathrm{MF1})^\dagger$ gives the decomposition of $\phi_l^{-1}$. Finally a similar argument, using $(\mathrm{MF1})^z$, establishes the corresponding result for $\phi_r$ and $\phi_r^{-1}$.

(ii) By using the $\eta$–$\epsilon$ adjunction equations, these decompositions are equivalent to the following relations in $\mathcal{F}$:

$$\begin{aligned} \check{\alpha} &= (\alpha^{-1})^* &&\Leftrightarrow &&(\mathrm{MF1})^z \\ \check{\rho} &= (\rho^{-1})^* &&\Leftrightarrow &&(\mathrm{CF}) \\ \check{\lambda} &= (\lambda^{-1})^* &&\Leftrightarrow &&(\mathrm{CF})^z \\ \check{\beta} &= \beta^* &&\Leftrightarrow &&(\mathrm{3DMG})^\dagger \end{aligned}$$

(iii) Cup off the middle input boundary circle in the relation (T9). Then apply relation $(\mathrm{ST3})^z$ to convert one of the II morphisms into a $\theta$. Apply the unit law (ST2) to convert this into a diagram whose source 1-morphism is 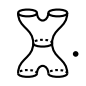. This is precisely the decomposition (33) of II in terms of $\phi_l$ and $\theta$. A similar argument, using $(\mathrm{T9})^\dagger$ and (ST3), establishes that $\mathrm{II}^{-1}$ equals the dagger of the right hand side of (33).

(iv) The decomposition (35) of III follows from (MF2) as we showed in Section 6.6. Similarly $(\mathrm{MF2})^\dagger$ implies that $\mathrm{III}^{-1}$ equals the dagger of the right hand side of (35). Similarly for III′. □

7.2. **Monoidal relations.** We now establish that $L$ preserves the relations in $\mathbf{F}(\mathcal{O})$, starting with the Monoidal relations. The pentagon relation (5) is (T1). The triangle relation (6) is (ST2).

7.3. **Balanced relations.** The hexagon relation (7) follows from (T5) and (CG) by 'cupping off' the second input circle. To see this, consider the first three 2-morphisms in the clockwise composite of (T5), namely $\beta \circ \alpha \circ \alpha$. If we cap off the second circle from the left, we obtain the following composite whose source 1-morphism is 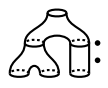:

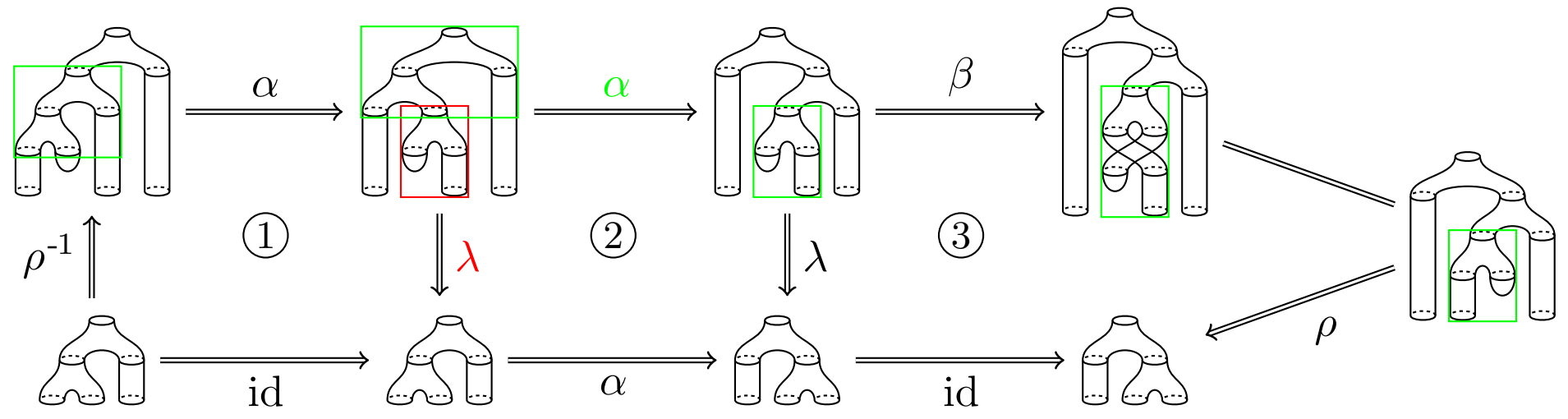

Diagram ① is (ST2), diagram ② is naturality, and diagram ③ is (CG). Carrying on in this manner produces precisely the hexagon relation (7).

The Balancing relations (8) and (9) are precisely (FCT2) and (FCT1) respectively.

7.4. **Rigidity relations.** Note that $\phi_l$ and $\phi_l^{-1}$ are *invertible generators* in $\mathcal{F}$—that is, there are relations $\phi_l \phi_l^{-1} = \text{id}$ and $\phi_l^{-1}\phi_l = \text{id}$. Hence checking that $L$ respects the Rigidity relations amounts to checking that the decompositions (10) and (11) hold in $\mathbf{F}(\mathcal{F})$. We showed this in Proposition 59 (i).

7.5. **Ribbon relation.** This follows from capping off $(\text{ST3})^x$, and then expanding II as $\phi_l \circ \theta \circ \phi_l^{-1}$ using Proposition 59 (iii), and also using (ST1) and (ST2).

7.6. **Biadjoint relations.** The relations (15) and (16) are precisely the Morse relations (MR1) and (MR2) respectively. Similarly for their rotated versions.

7.7. **Pivotality relations.** Recall from Proposition 59 that the relation (CF) is equivalent to the statement that $\check{\rho} = (\rho^{-1})^*$ in $\mathbf{F}(\mathcal{F})$. Similarly, $(\text{CF})^\dagger$ is equivalent to the statement that $\check{\rho} = {}^*(\rho^{-1})$. Hence $\rho^* = {}^*\rho$. The argument in the proof of Lemma 16 can then be run in reverse to prove that Pivotality on the cylinder (31) holds in $\mathbf{F}(\mathcal{F})$. Capping off the top and bottom then gives the Pivotality relation on the sphere (17). A similar argument, using $(\text{CF})^z$, establishes $(17)^z$.

7.8. **Modularity relations.** We know from Proposition 59(iv) that III and its inverse $\text{III}^{-1}$ can be decomposed in $\mathbf{F}(\mathcal{F})$ according to the right hand side of (35) and its dagger respectively. We can define $A$ in the 2-category $\mathbf{F}(\mathcal{F})$ in terms of III as $A := \theta^{-1}_{\text{left}} \circ \text{III} \circ \theta^{-1}_{\text{left}}$.

Now, III is an invertible generator in $\mathcal{F}$. Hence the relations $\text{III} \circ \text{III}^{-1} = \text{III}^{-1} \circ \text{III} = \text{id}$ in $\mathcal{F}$ imply that $A$ is invertible, with inverse $A^\dagger$. Moreover, we have established that the Ribbon and Biadjoint relations hold in $\mathbf{F}(\mathcal{F})$. Hence we can apply Proposition 25 to conclude that the Modularity relation holds in $\mathbf{F}(\mathcal{F})$.

7.9. **Anomaly-freeness.** At this point we have established that each of the relations for $\mathbf{F}(\mathcal{M})$ holds in $\mathbf{F}(\mathcal{F})$. Hence Theorem 49 can be applied, and so the equation (54) is true in $\mathbf{F}(\mathcal{F})$. On the other hand, in $\mathbf{F}(\mathcal{F})$, we have the relation (T4), which implies the relation (T4′) (which has no $x$ factor!) as we showed in Section 6.7. Since (54) and (T4′) both hold, we must have:

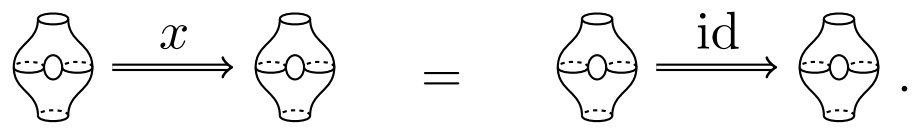

It is not necessary to specify exactly where $x$ is acting, due to $x$ being central in the sense of Corollary 45. We need to show that $x = \mathrm{id}$ on the cylinder, which is equivalent to $x' = \mathrm{id}$ from Lemma 42. Pre- and post-composing the above leads to:

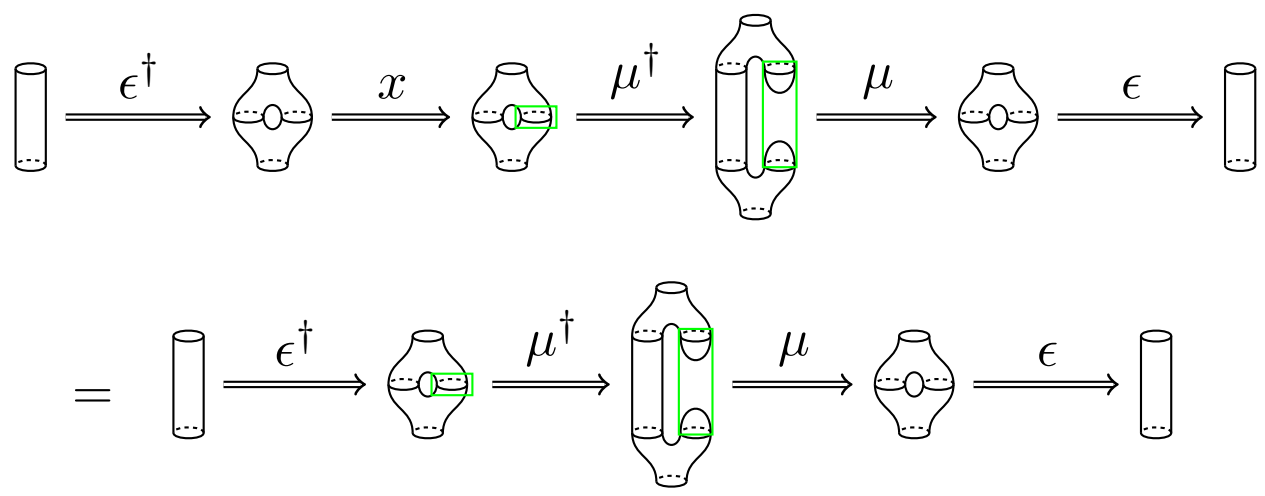

Now commute $x$ to the left, using the centrality of $x$ (Corollary 45). Apply the Pivotality relation (31) on the cylinder (which holds in $\mathbf{F}(\mathcal{F})$ as we showed above) to both sides. Hence $x = \mathrm{id}$ on the cylinder.

7.10. **Proof of Corollary 2.** We can now combine these results to establish the equivalence of the 2-Morse presentation $\mathcal{F}$ with the modular presentation $\mathcal{O}$. We have established symmetric monoidal functors $K\colon \mathbf{F}(\mathcal{F}) \to \mathbf{F}(\mathcal{O})$ and $L\colon \mathbf{F}(\mathcal{O}) \to \mathbf{F}(\mathcal{F})$. Clearly $K \circ L = \mathrm{id}_{\mathbf{F}(\mathcal{O})}$, and from Proposition 59 we see that $L \circ K = \mathrm{id}_{\mathbf{F}(\mathcal{F})}$. Hence Theorem 1, that there is an equivalence $\mathbf{F}(\mathcal{F}) \simeq \mathbf{F}(\mathcal{O})$, follows. Similarly Corollary 2, that there is an equivalence $\mathbf{Bord}^{\mathrm{or}}_{123} \simeq \mathbf{F}(\mathcal{O})$, follows from the theorem and the main result of [7], that there is a symmetric monoidal equivalence $\mathbf{F}(\mathcal{F}) \simeq \mathbf{Bord}^{\mathrm{or}}_{123}$.

## Appendix A. The 2-Morse presentation

In [7], we used 2-Cerf theory to show that the oriented 3-dimensional bordism bicategory $\mathbf{Bord}^{\mathrm{or}}_{123}$ is equivalent to the free quasistrict symmetric monoidal 2-category $\mathbf{F}(\mathcal{F})$ on a certain presentation $\mathcal{F}$. We call $\mathcal{F}$ the *2-Morse presentation.* In this Appendix we recall the definition of $\mathcal{F}$ from [7]. For our conventions on *daggers* and *rotations* of relations, see Appendix B.

**Definition 60.** The 2-Morse presentation $\mathcal{F}$ is defined as follows:

- Generating object: 
- Generating 1-morphisms:

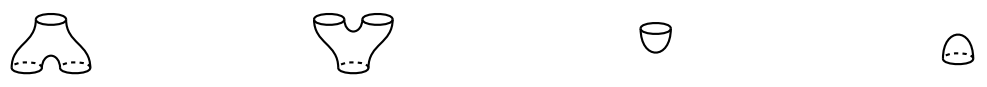

- Invertible generating 2-morphisms:

(66)

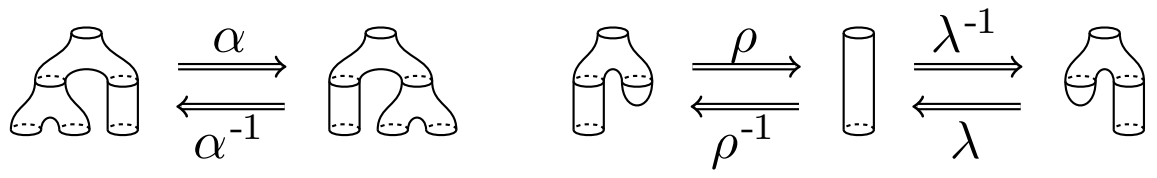

(67) 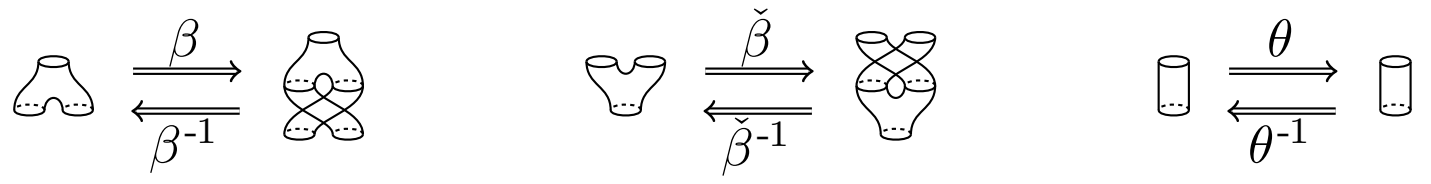

(68) 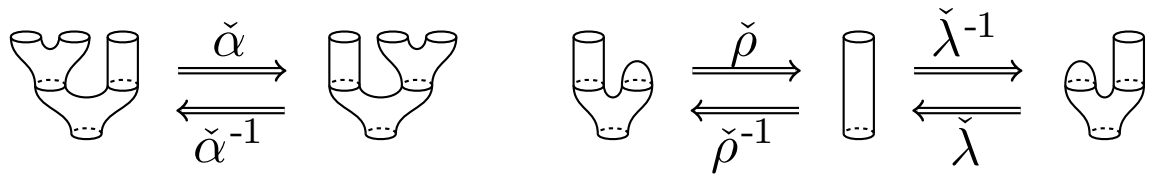

(69) 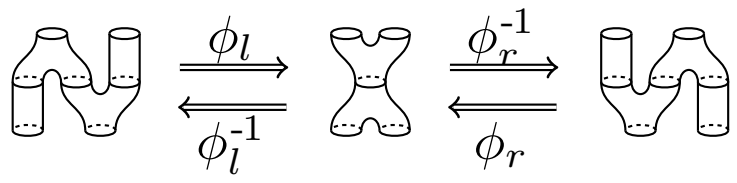

(70) 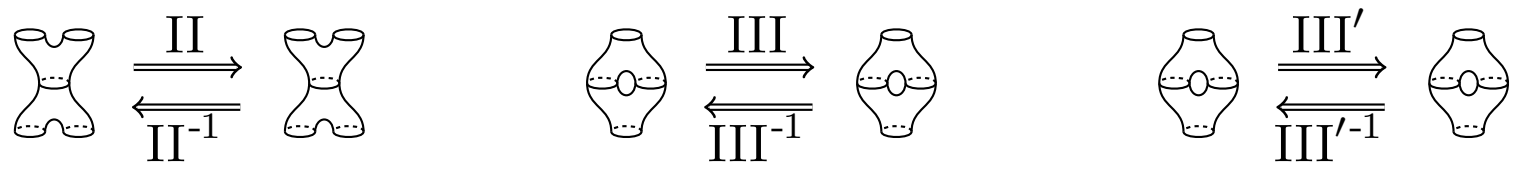

Noninvertible generating 2-morphisms:

(71) 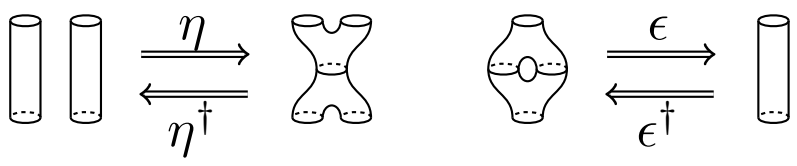

(72) 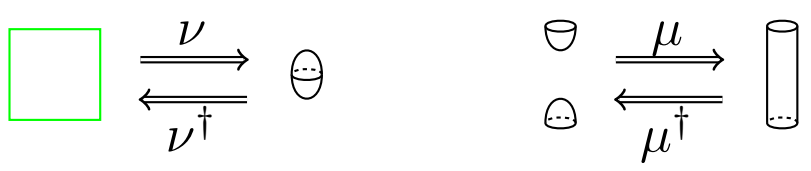

The relations are listed below in Sections A.1 to A.11.

We now list the relations. We only draw one version of each relation, with the understanding that all daggers and rotated versions of that relation about the $x$, $y$ and $z$ axes are also relations. Often applying such an operation reproduces the same relation, so for clarity we list the number of distinct equations each picture represents next to the equation label. For instance, '$x\dagger$' means there are four versions of the relation $\xi$, namely $\xi$, $\xi^\dagger$, $\xi^x$ and $(\xi^x)^\dagger$.

A.1. **Inverses.** Each of the invertible generating 2-morphisms $\omega$ satisfies $\omega\omega^{-1} = \mathrm{id}$ and $\omega^{-1}\omega = \mathrm{id}$.

A.2. **Morse.**

(MR1) $x\dagger$ 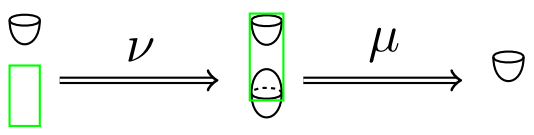 $=$ id

(MR2) $x\dagger$ 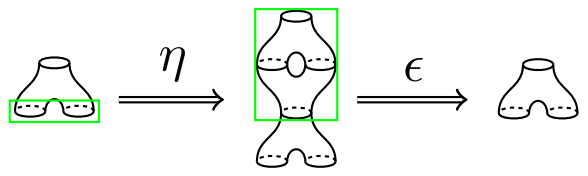 $=$ id

### A.3. **Cusp flip.**

(CF) $xz\dagger$

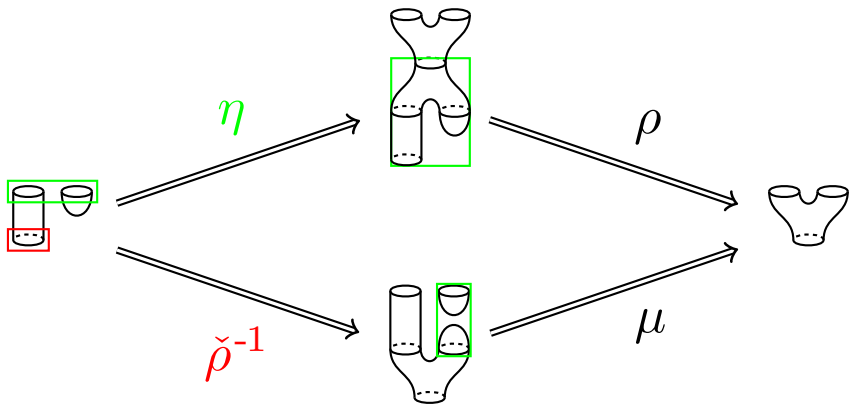

### A.4. **Swallowtail.**

(ST1) $xz$

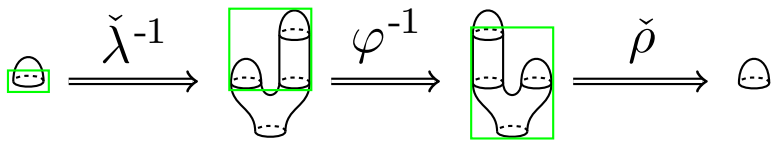

$=$ id

(ST2) $z$

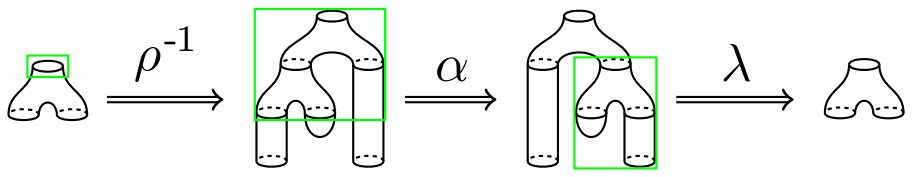

$=$ id

(ST3) $xz$

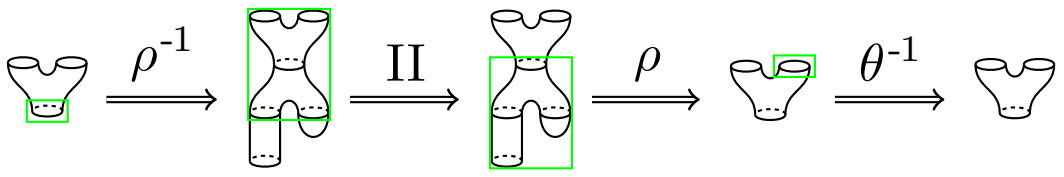

$=$ id

### A.5. **Beak.**

(B1) $xz$

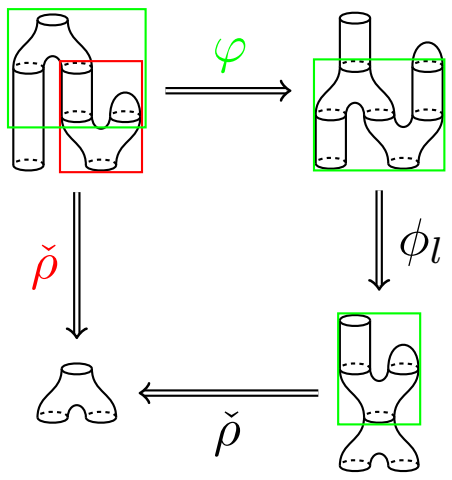

(B2) $xz$

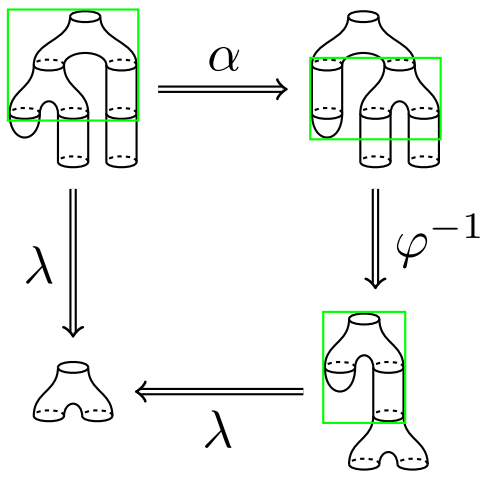

### A.6. **Morse flip.**

(MF1) $xz\dagger$

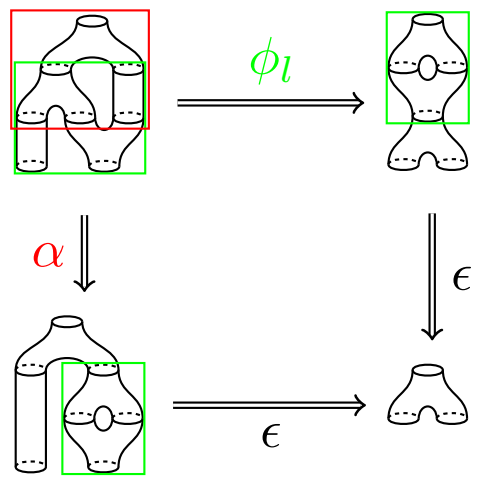

(MF) $xz\dagger$

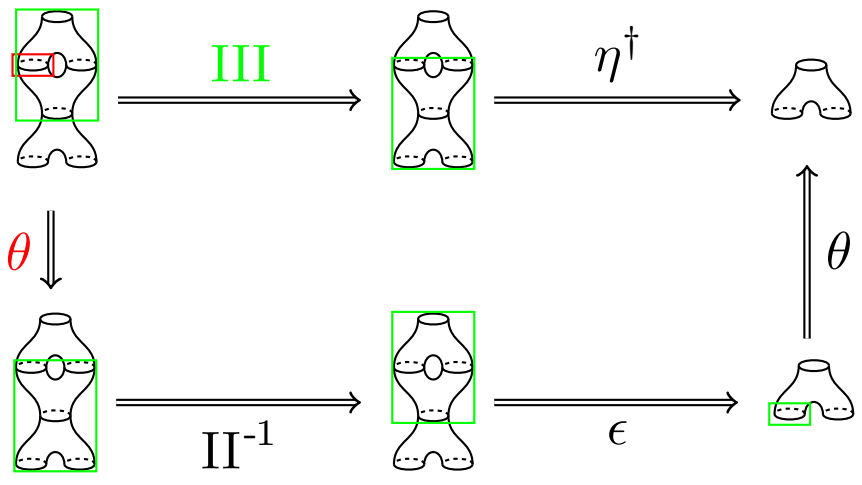

### A.7. **Triple point.**

(T1) $x$

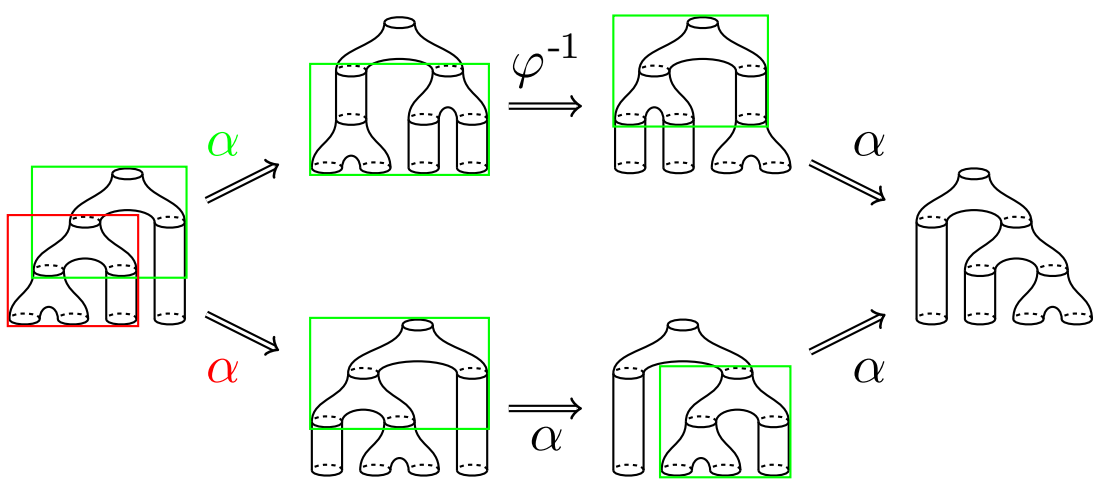

(T2) $xz$

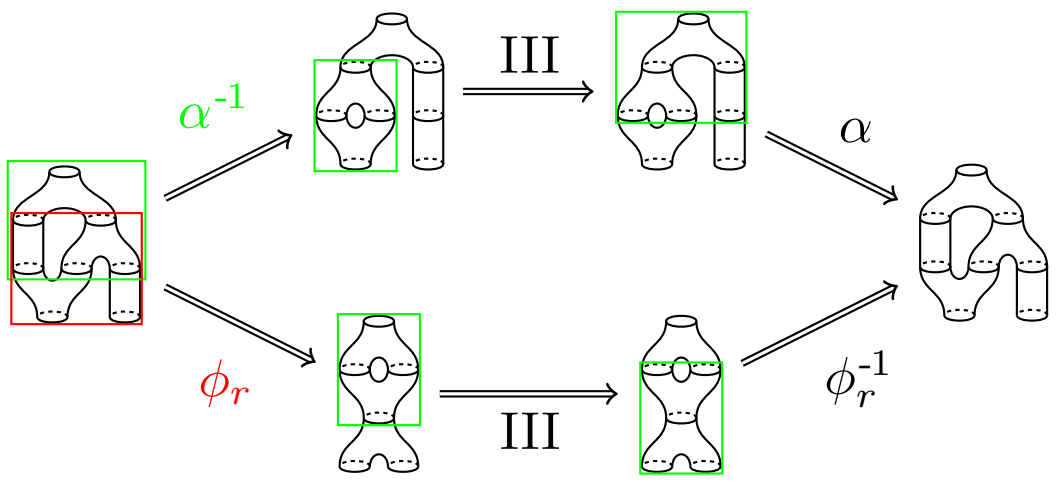

(T3) $xz$

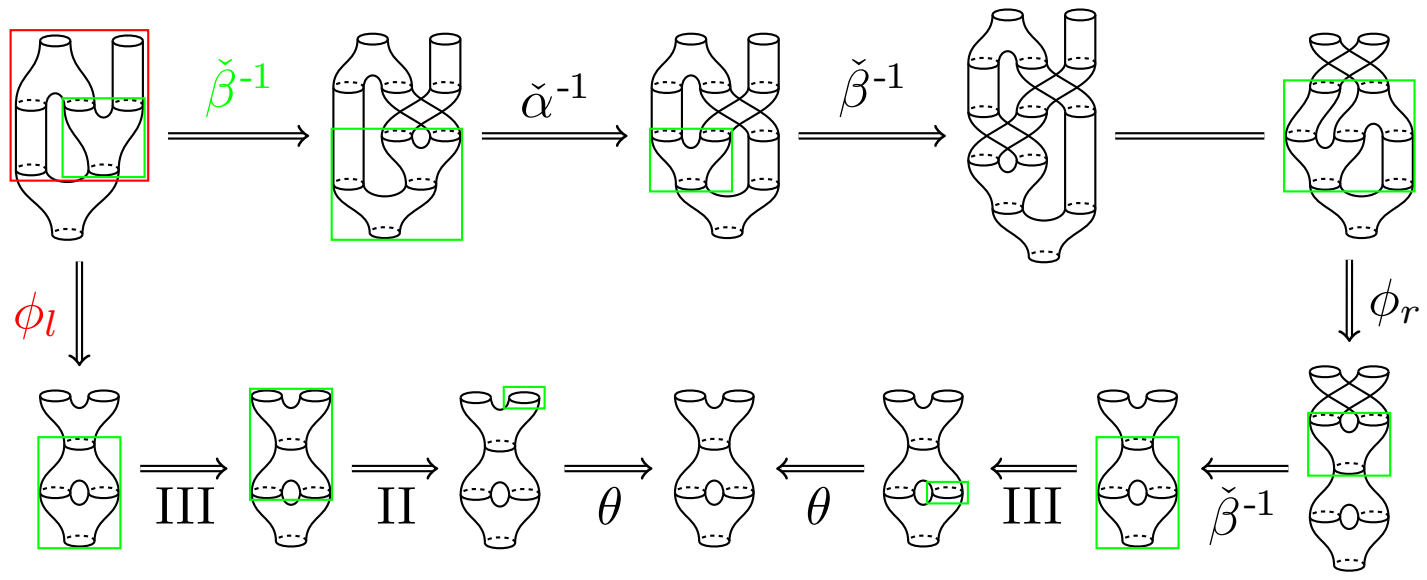

(T4) $xz$

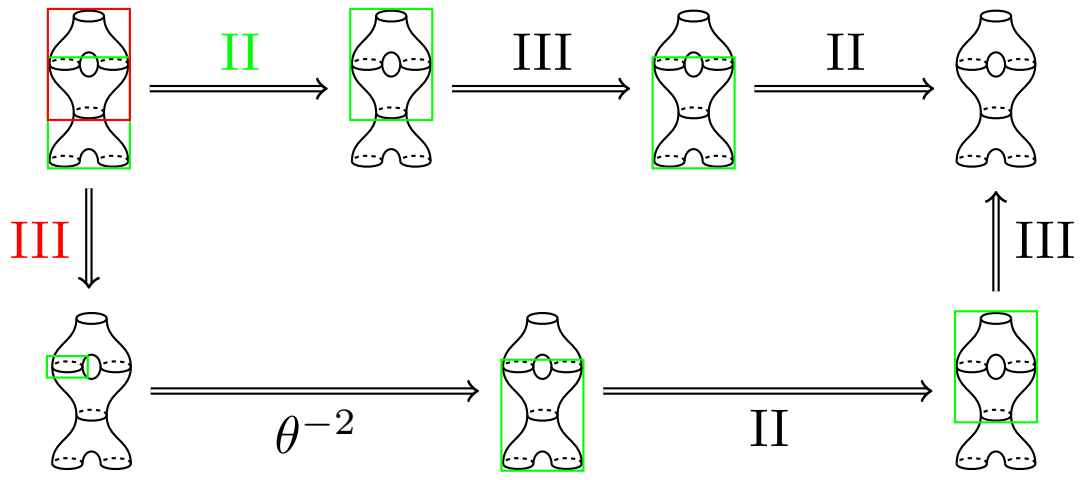

(T5) $xz$

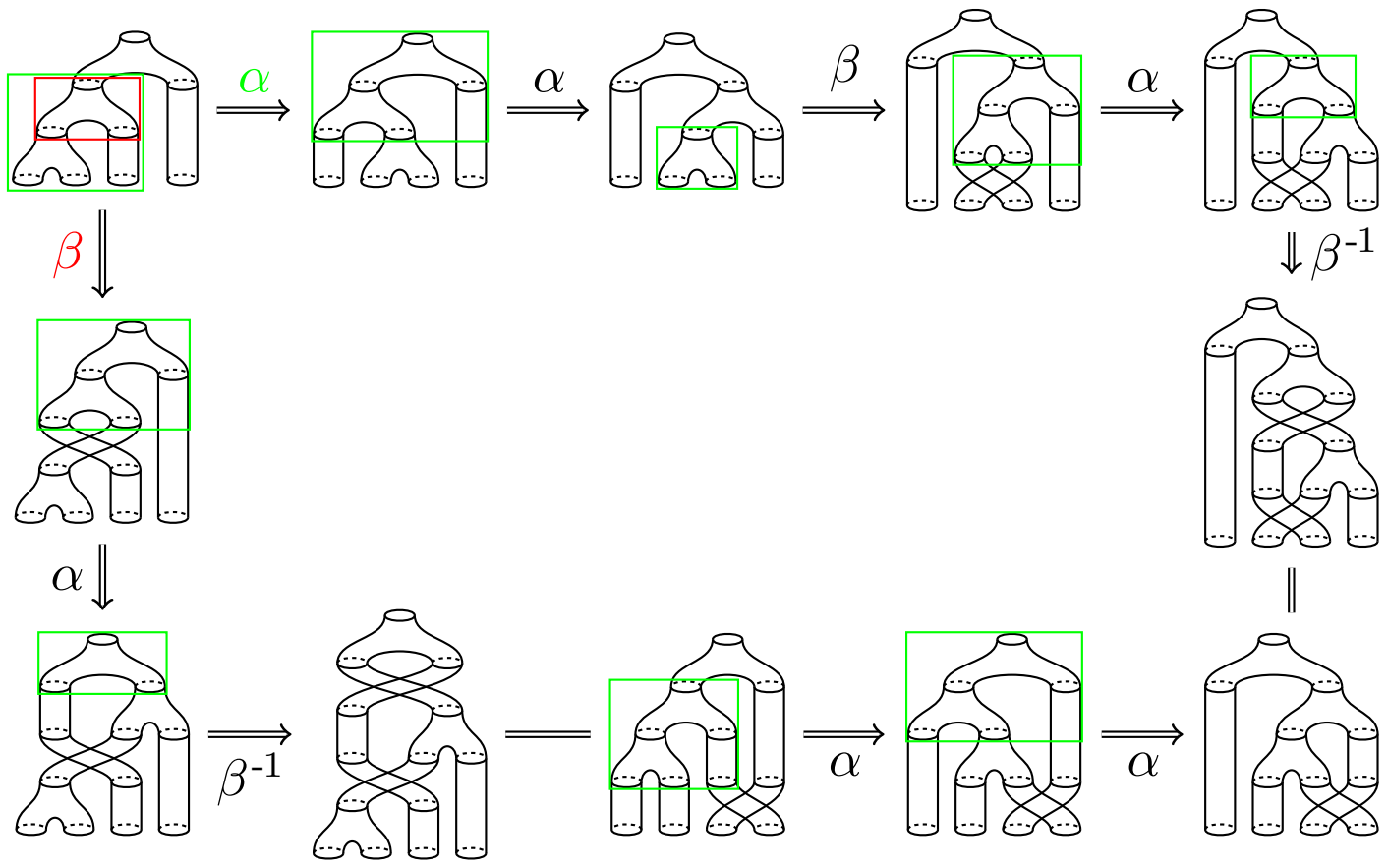

(T6) $xz$

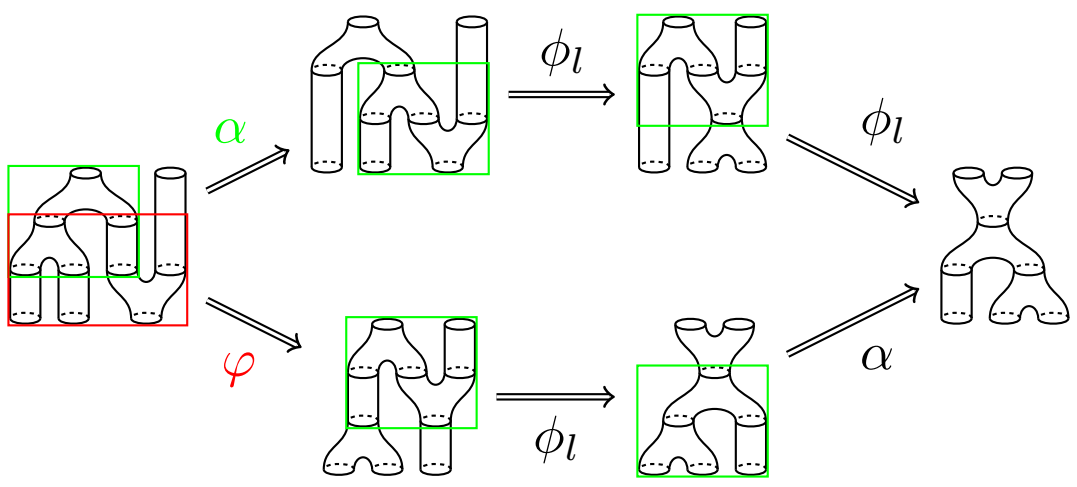

(T7) $xz$

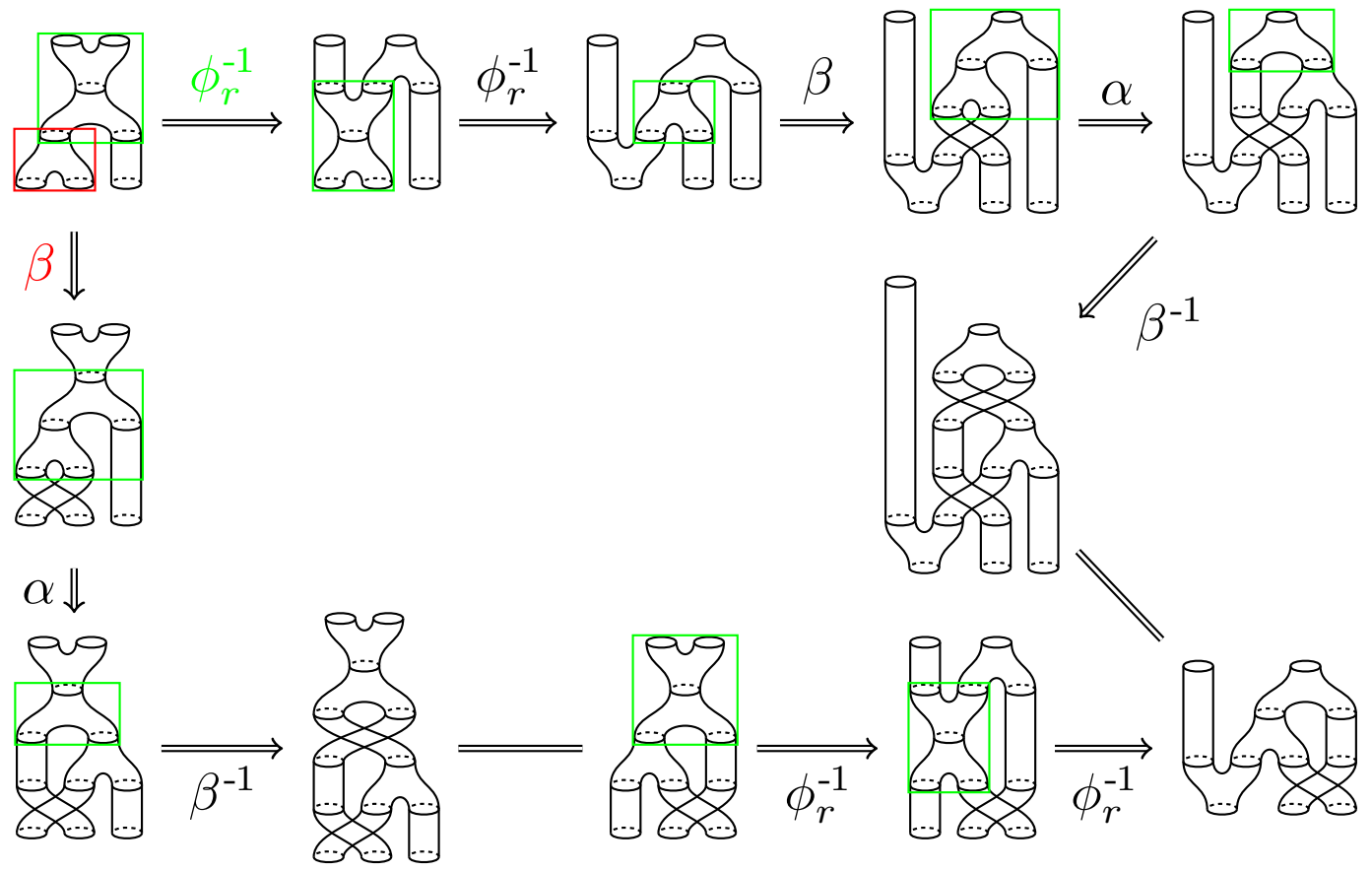

(T8) $xz$

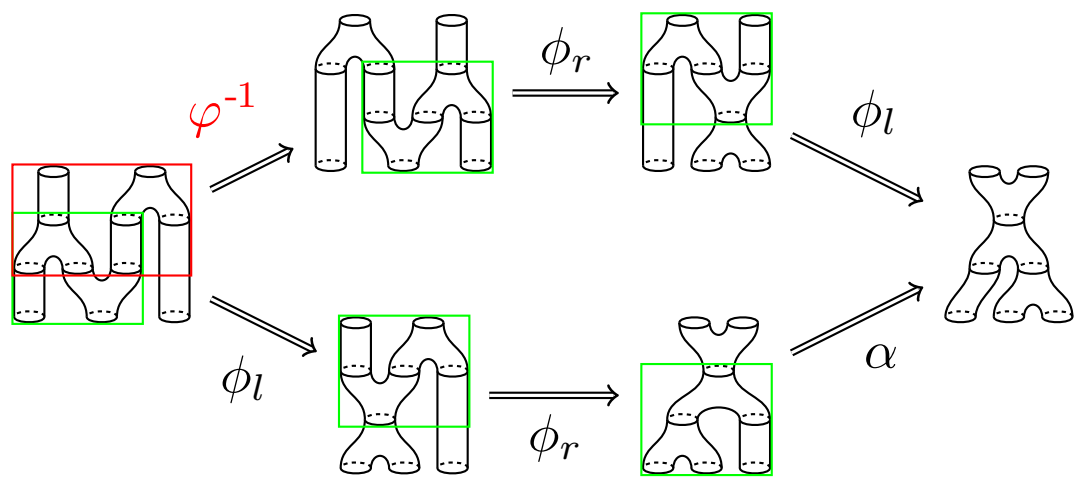

(T9) $xz$

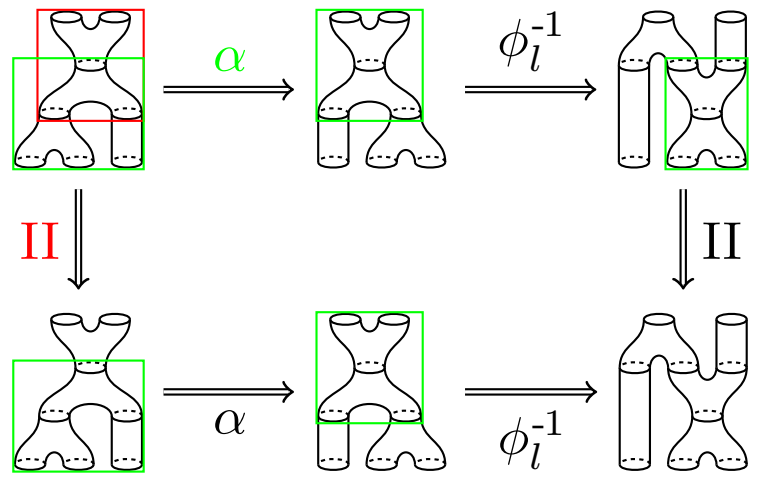

(T10) $xz$

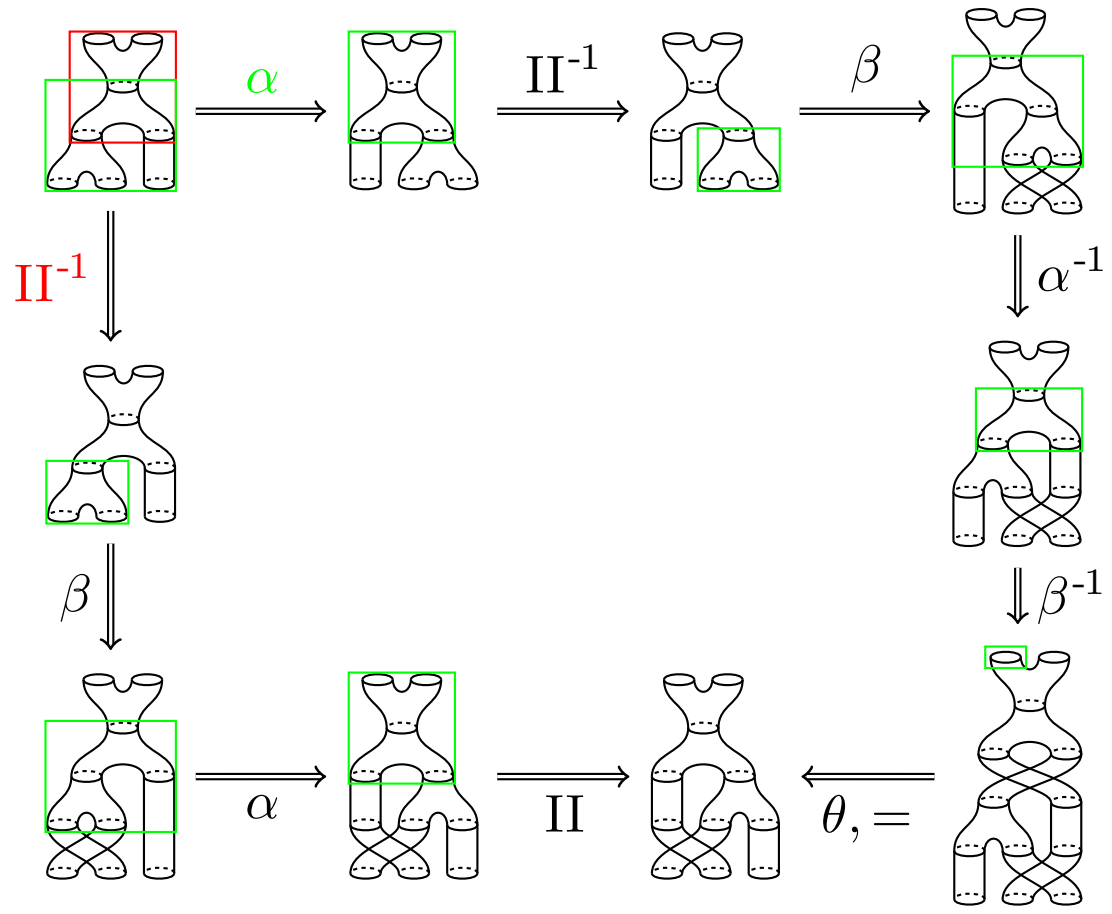

A.8. **Fold-codim 2 gluing.**

(FCT1) $z$

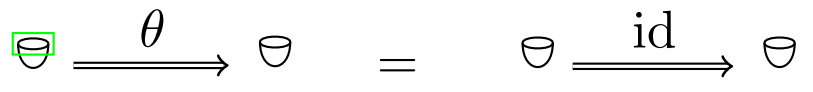

(FCT2) $z$

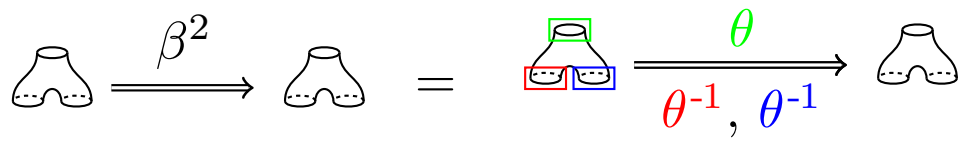

A.9. **Fold-fold-glue.**

(FFG1) $z$

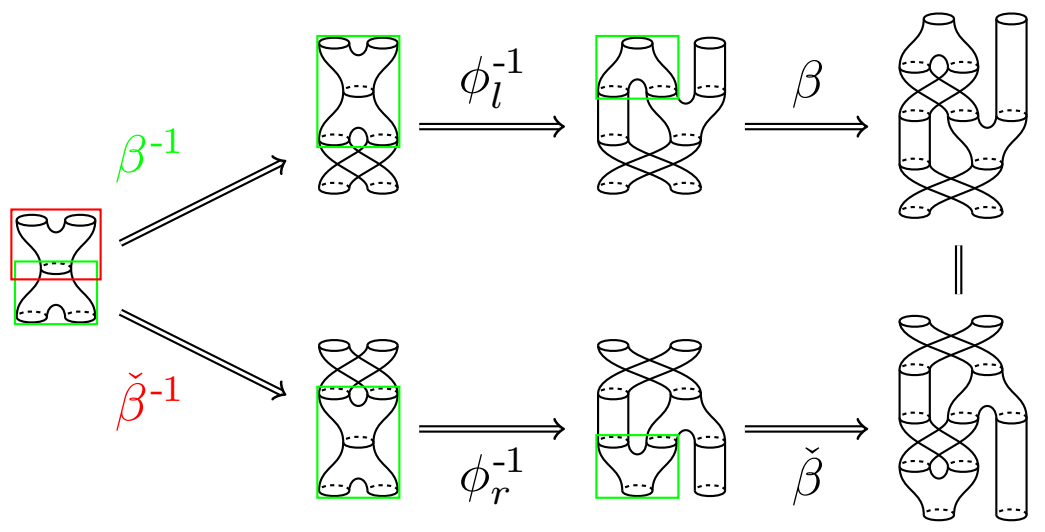

(FFG2)

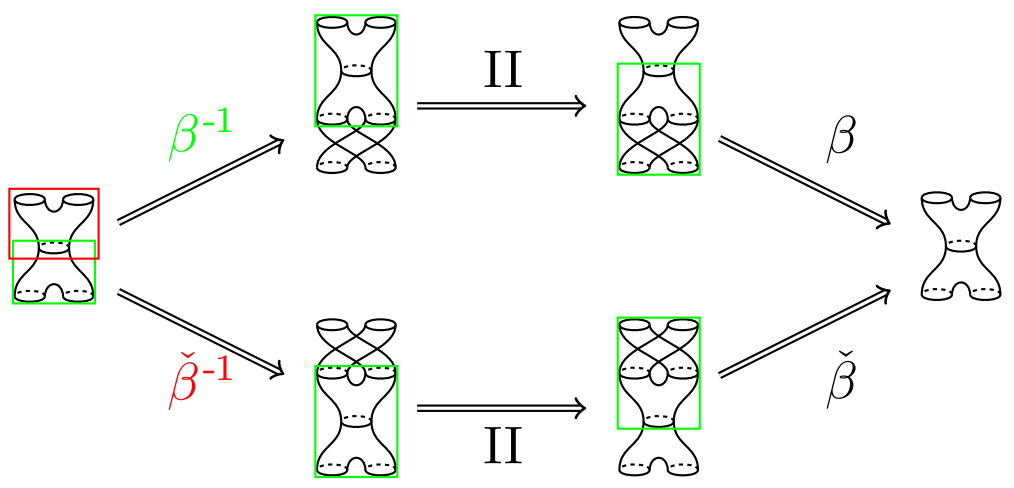

(FFG3)

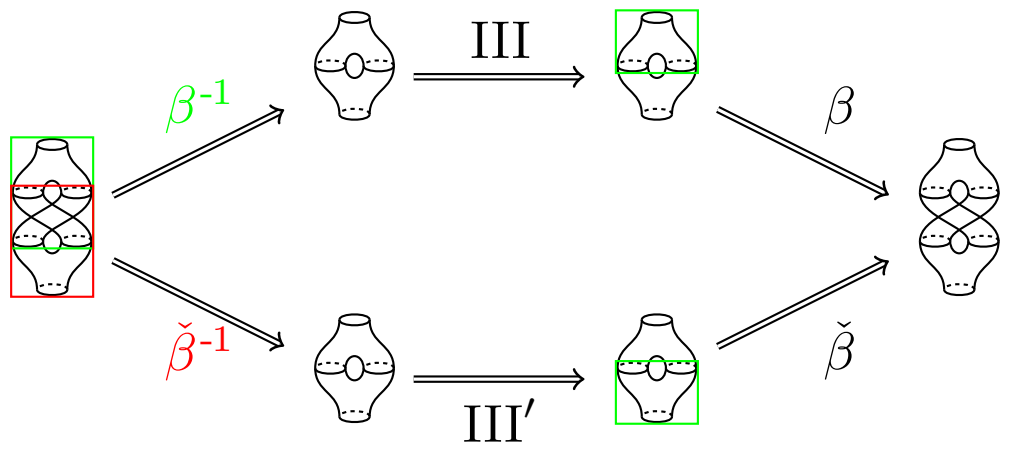

(FFG4) $z$

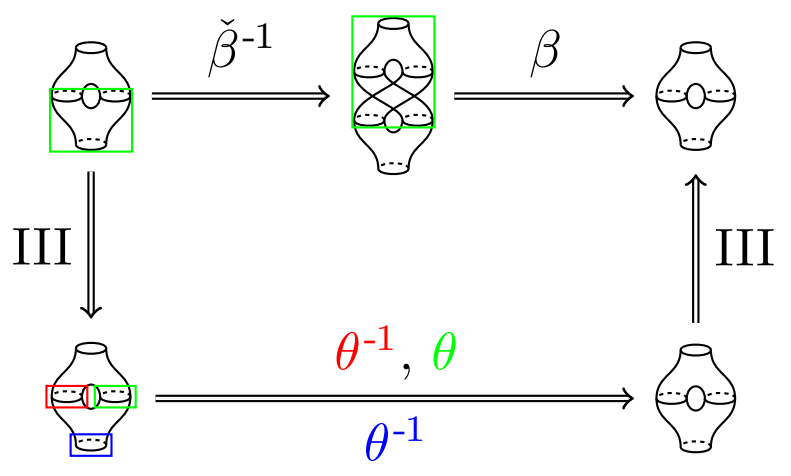

A.10. **Cusp-glue.**

(CG) $xz$

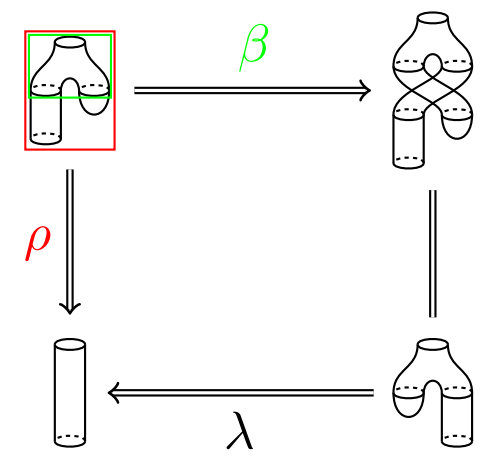

A.11. **3d Morse-glue.**

(3DMG) †

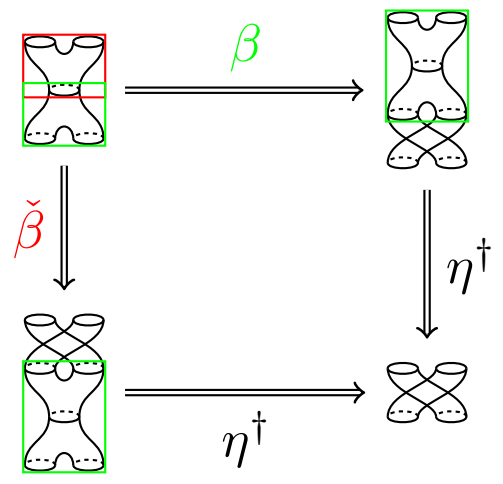

## Appendix B. Conventions regarding daggers and rotations

Let $\mathcal{H}$ refer to one of the presentations $\mathcal{F}$, $\mathcal{R}$, $\mathcal{M}$, or $\mathcal{O}$.

B.1. **Dagger of an equation.** Let $\omega$ be a generating 2-morphism of $\mathcal{H}$. We set

$$\text{dagger of } \omega := \begin{cases} \omega^{-1} & \text{if } \omega \text{ is invertible} \\ \omega^{\dagger} & \text{if } \omega = \eta,\ \epsilon,\ \mu,\ \nu \\ \omega & \text{if } \omega = \eta^{\dagger},\ \epsilon^{\dagger},\ \mu^{\dagger},\ \nu^{\dagger}. \end{cases}$$

and we write $\omega^{\dagger}$ for the dagger of $\omega$.

The *dagger* of an equation $\xi$ of the form $L = R$, where $L$ and $R$ are composites of generating 2-morphisms in $\mathcal{H}$, is a new equation $\xi^{\dagger}$ asserting that $L^{\dagger} = R^{\dagger}$. Here, $L^{\dagger}$ (resp. $R^{\dagger}$) is the new composite in $\mathbf{F}(\mathcal{H})$ obtained by taking the dagger of each 2-morphism $\omega$ in $L$ (resp. $R$) and then composing the result in the reverse order. For instance, if $L$ was the composite

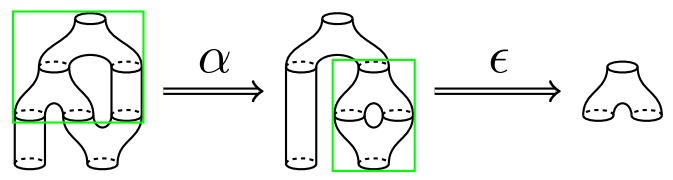

then $L^{\dagger}$ would be the composite

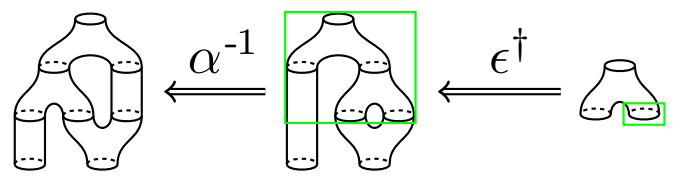

| $\omega$ | $\omega^x$ | $\omega^y$ | $\omega^z$ |
|---|---|---|---|
| $\alpha$ | $\breve{\alpha}$ | $\breve{\alpha}^{-1}$ | $\alpha^{-1}$ |
| $\rho$ | $\breve{\rho}$ | $\breve{\lambda}$ | $\lambda$ |
| $\beta$ | $\breve{\beta}$ | $\breve{\beta}$ | $\beta$ |
| $\theta$ | $\theta$ | $\theta$ | $\theta$ |
| $\breve{\alpha}$ | $\alpha$ | $\alpha^{-1}$ | $\breve{\alpha}^{-1}$ |
| $\breve{\rho}$ | $\rho$ | $\lambda$ | $\breve{\lambda}$ |
| $\breve{\lambda}$ | $\lambda$ | $\rho$ | $\breve{\rho}$ |
| $\breve{\beta}$ | $\beta$ | $\beta$ | $\breve{\beta}$ |
| $\phi_l$ | $\phi_r$ | $\phi_l$ | $\phi_r$ |
| $\phi_r$ | $\phi_l$ | $\phi_r$ | $\phi_l$ |
| II | II | II | II |
| III | III | III′ | III′ |

FIGURE 4. Definitions of rotated versions of invertible 2-morphism generators.

B.2. **Dagger 2-functor.** In each of the presentations $\mathcal{H}$, the following holds: if $\xi$ is a relation then $\xi^\dagger$ is also a relation. Hence, if we define the dagger of an interchangor 2-morphism (see Appendix C) to be its inverse, then the dagger operation extends to a strict symmetric monoidal 2-functor

$$\dagger\colon \mathbf{F}(\mathcal{H}) \to \mathbf{F}(\mathcal{H})$$

which is covariant on 1-morphisms and contravariant on 2-morphisms, and which squares to the identity.

B.3. **Rotation of an equation.** We use the 3-dimensional system of axes where $x$ runs from left-to-right, $y$ runs into the page, and $z$ runs upward:

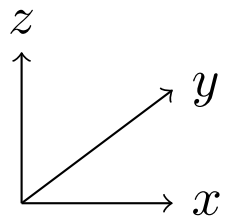

Let $u = x$, $y$, or $z$, and let $\xi$ be an equation of the form $L = R$ where $L$ and $R$ in $\mathbf{F}(\mathcal{H})$ are certain composites of 2-morphisms of $\mathcal{H}$. The *rotation of* $\xi$ *by* 180° *about the* $u$*-axis* (in brief, '$\xi$ rotated about the $u$-axis') is a new equation $\xi^u$ asserting that $L^u = R^u$. Here the composite $L^u$ (resp. $R^u$) is defined as follows. Firstly 'rotate' each 1-morphism occuring in $L$ in the intuitive manner. For instance,

rotated about the $z$-axis is 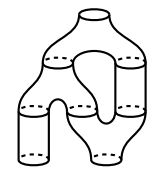.

Next, 'rotate' each generating 2-morphism $\omega$ in $L$ by replacing it with $\omega^u$, defined as follows. If $\omega$ is an invertible 2-morphism, consult Table 4, supplemented with the rule $(\omega^{-1})^u := (\omega^u)^{-1}$. If $\omega$ is $\eta$, $\epsilon$, $\mu$, $\nu$ or their daggers, we set $\omega^u = \omega$.

Table 4 needs to be appropriately interpreted in the context of the presentation $\mathcal{H}$. For instance, the equation (32) takes place in $\mathbf{F}(\mathcal{O})$. Its rotated version $(32)^z$ is

the following composite:

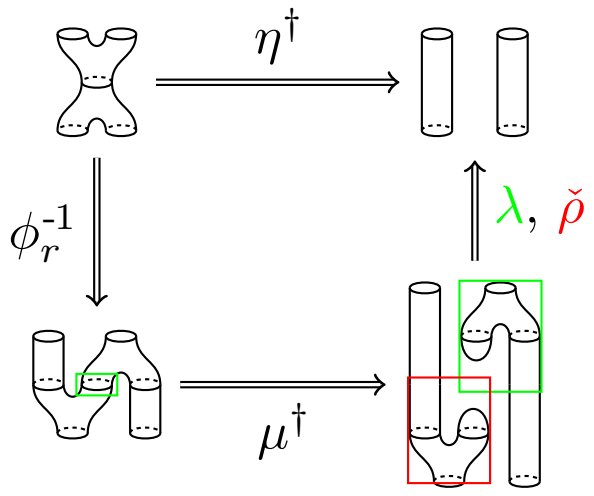

At the point where (32) was written down, the 2-morphisms $\phi_l^{-1}$, $\phi_r^{-1}$, $\check{\lambda}$ and $\check{\rho}$ had been defined as composites of the generators in $\mathbf{F}(\mathcal{O})$. Hence, the original equation (32) and its rotated version $(32)^z$ make sense. Needless to say we only ever talk about the 'rotated version of a relation' when it makes sense to do so.

As another example, consider equation (T3), which is a defining relation in $\mathcal{F}$. Its rotated version $(\mathrm{T3})^x$ is the following composite:

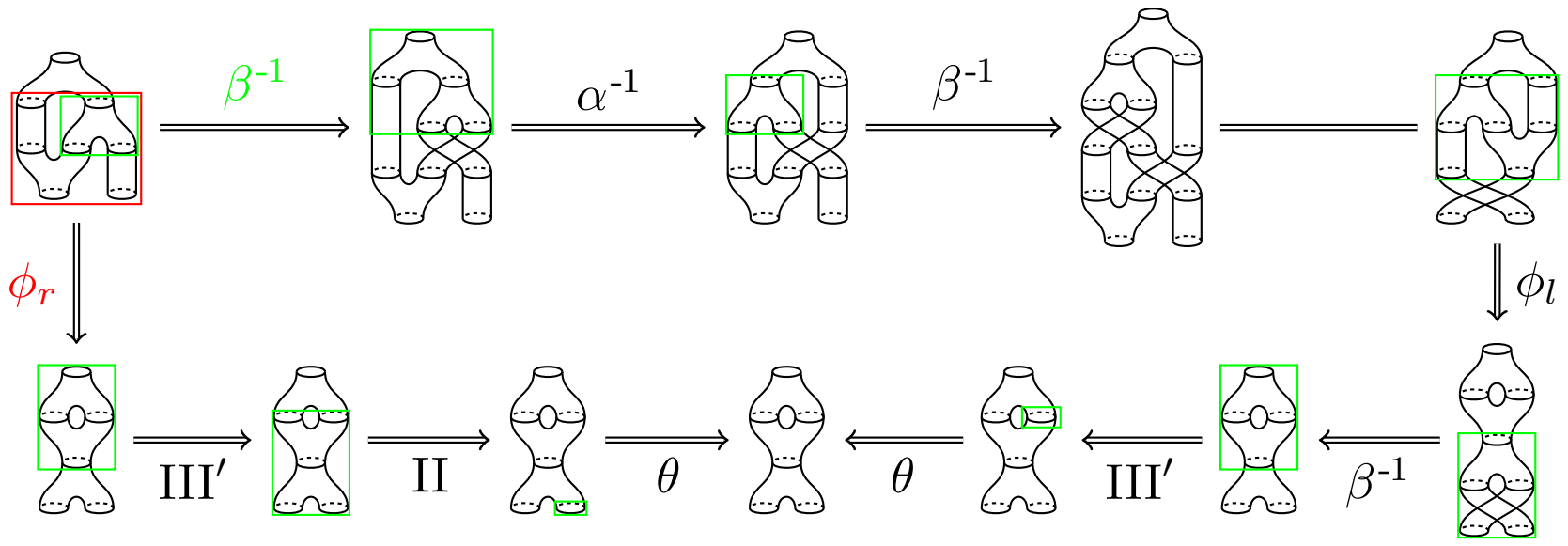

## Appendix C. Presentations of quasistrict symmetric monoidal 2-categories

Establishing the definition of a fully weak symmetric monoidal bicategory, and proving associated coherence theorems, has been an extensive effort by a number of authors; see the monographs [18, 40] for references and a pedagogical account. In the main body of this paper we will only need the notions of:

- *a quasistrict symmetric monoidal 2-category* [40];
- a *generators and relations presentation* for a quasistrict symmetric monoidal 2-category, in the formal sense of a *quasistrict symmetric monoidal 3-computad* [40] (referred to in this paper as a 'presentation' for simplicity); and
- the *wire diagrams* graphical calculus for quasistrict symmetric monoidal 2-categories [5].

We only need to deal with quasistrict symmetric monoidal 2-categories since in this paper we do not work with the 3-dimensional oriented bordism bicategory $\mathbf{Bord}_{123}^{\mathrm{or}}$ directly (which is a fully weak symmetric monoidal bicategory). Rather, our starting point is the result of [7] that $\mathbf{Bord}_{123}^{\mathrm{or}}$ is equivalent, as a a symmetric monoidal bicategory, to the quasistrict symmetric monoidal 2-category $\mathbf{F}(\mathcal{F})$. Our goal is merely to simplify this presentation, that is to establish a new computad $\mathcal{O}$ whose associated quasistrict symmetric monoidal 2-category $\mathbf{F}(\mathcal{O})$ is equivalent to $\mathbf{F}(\mathcal{F})$.

C.1. **Quasistrict symmetric monoidal 2-categories.** We recall some definitions. The fully weak definition of a *symmetric monoidal bicategory* can be found in [19, 40, 41]. We adopt here the notation and terminology of [5, 40]. Let $(\underline{\mathrm{Cat}}, \otimes_G)$ denote the category of strict 2-categories and strict 2-functors, equipped with the Gray tensor product $\otimes_G$ of 2-categories. A *semistrict monoidal 2-category* (equivalently a *Gray monoid*) is a monoid in $(\mathrm{Cat}, \otimes_G)$. Every monoidal bicategory is equivalent to a semistrict monoidal 2-category [15, 17]. A *Crans semistrict symmetric monoidal 2-category* is a symmetric monoidal bicategory whose underyling monoidal bicategory is a semistrict monoidal 2-category, with some further normalization conditions on the coherence data. Every symmetric monoidal bicategory is equivalent to a Crans semistrict symmetric monoidal bicategory [19].

There is an even stricter notion of a Crans semistrict symmetric monoidal bicategory, introduced in [40], a *quasistrict symmetric monoidal 2-category*. The definition can be formulated in such a way (the 'stringent' form of the definition from [5]) that the only nontrivial coherence data are the *interchangor* 2-isomorphisms coming from the monoidal 2-category structure. That is, the 2-isomorphisms

$$\varphi_{f,g}\colon (f \otimes \mathrm{id}) \circ (\mathrm{id} \otimes g) \Rightarrow (\mathrm{id} \otimes g) \circ (f \otimes \mathrm{id}) \tag{73}$$

where $f$ and $g$ are 1-morphisms. We will essentially spell out this 'stringent' form of the definition below, while introducing the graphical calculus. For a pedagogical account, see [5].

**Theorem 61** ([40])**.** *Every symmetric monoidal bicategory is equivalent to a quasistrict symmetric monoidal 2-category.*

We will not actually need this result in this paper but it is useful to bear in mind.

C.2. **Graphical notation.** We briefly review the *wire diagrams* graphical notation for quasistrict symmetric monoidal 2-categories from [5]. The basic idea is that the tensor product direction runs out of the page, composition of 1-morphisms and (horizontal) composition of 2-morphisms runs up the page, and (vertical) composition $\circ$ of 2-morphisms runs from left-to-right[6]. For instance:

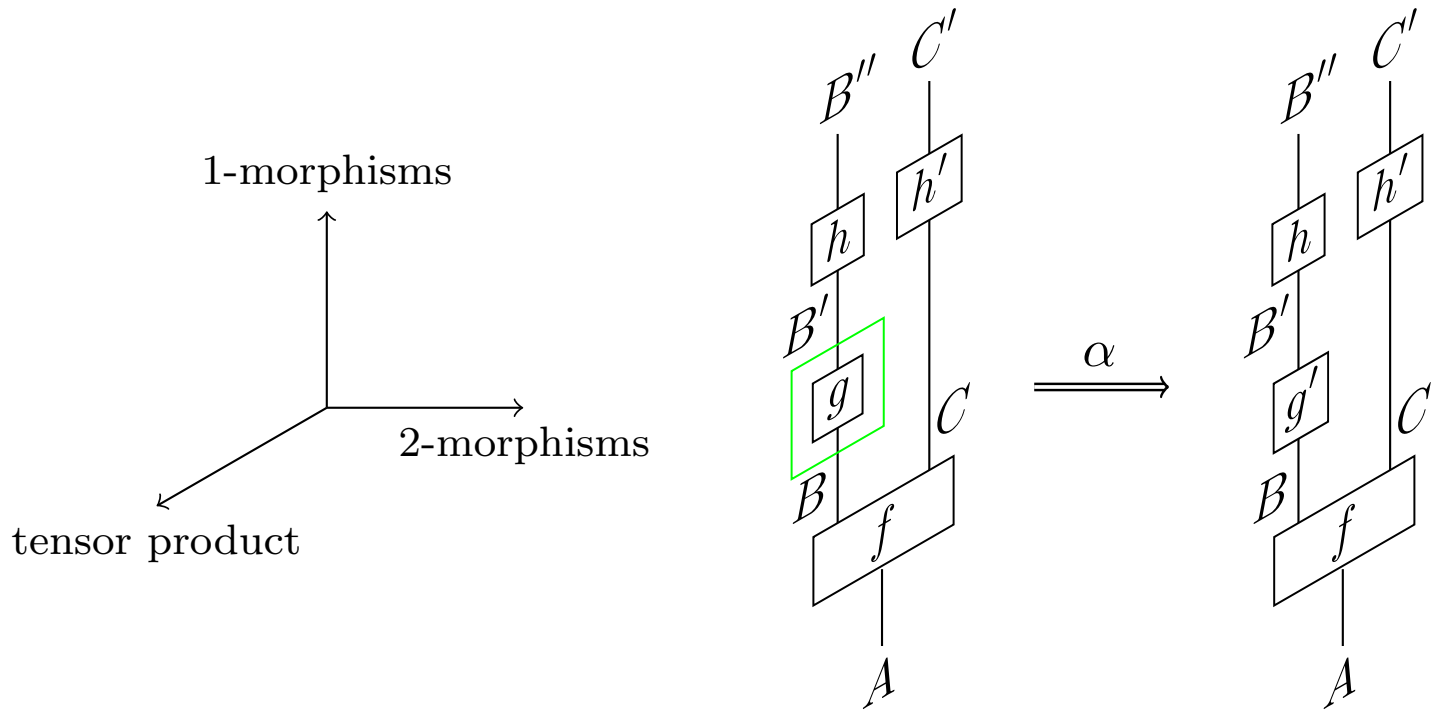

[6]Unfortunately, what is usually called *vertical* composition $\circ$ of 2-morphisms runs horizontally in wire diagrams, and what is usually called *horizontal* composition $*$ runs vertically!

To make the diagrams clearer, they will usually just be drawn flat in the page (but the 3-dimensional picture should be kept in mind):

(74)
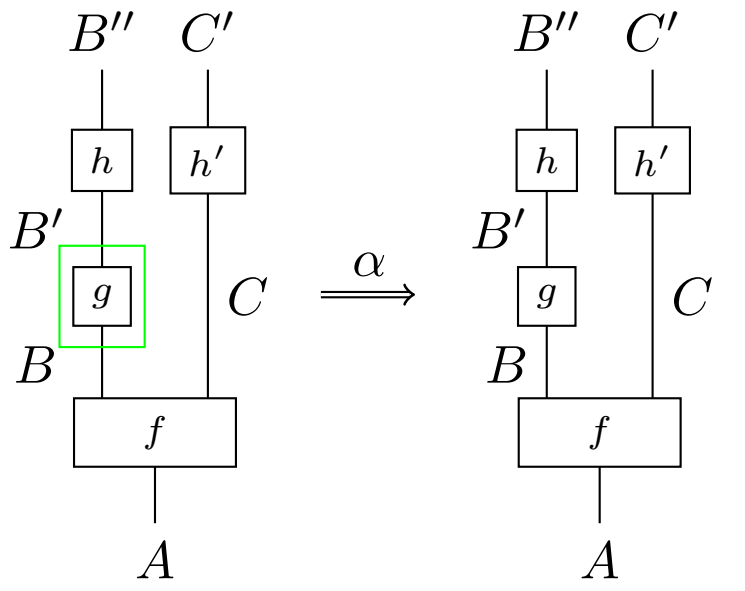

Let us explain this picture. It represents a 2-morphism between two 1-morphisms. Each 1-morphism is to be evaluated into a precise composite using $\circ$ and $\otimes$ according to the prescription *tensor first, then compose*. The source of the 2-morphism $\alpha$ is depicted with a coloured box. The full 2-morphism is then obtained by whiskering $\alpha$ with identity 2-morphisms corresponding to the regions where it is not acting. So, (74) evaluates to:

$$(h \otimes h') \circ g \circ f \xrightarrow{\mathrm{id}_{h\otimes h'} * (\alpha \otimes \mathrm{id}) * \mathrm{id}_f} (h \otimes h') \circ g' \circ f.$$

Note that this is a precise notation and not merely a mnemonic device. No further parentheses are needed in the above evaluation since the underlying 2-category is strict, and $A \otimes -$ and $- \otimes A$ are strict 2-functors, for each object $A$.

The interchangor 2-isomorphisms (73) are thus drawn as follows:

(75)
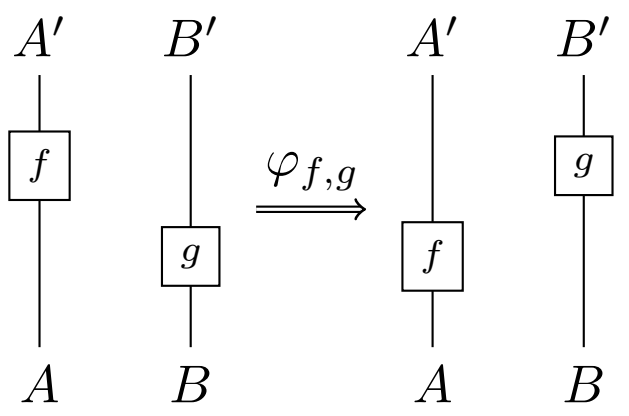

Importantly, though the interchangor is nontrivial, due to our conventions the following 'nudging' equation still holds:

(76)
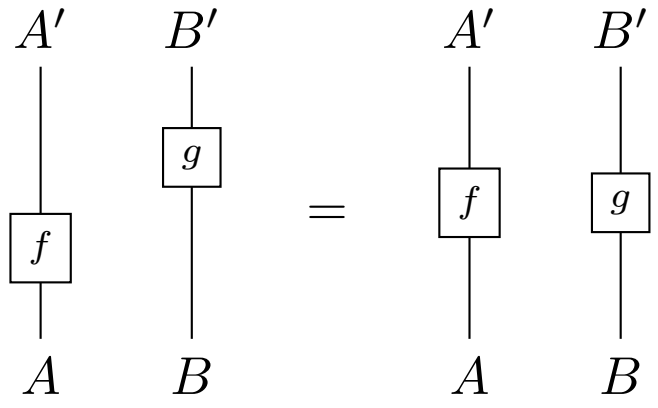

Note that by our 'tensor first, then compose' convention, the equation (76) looks as follows algebraically:

$$(\mathrm{id} \otimes g) \circ (f \otimes \mathrm{id}) = f \otimes g \tag{77}$$

If we had adopted the 'opcubical' convention to define a semistrict monoidal 2-category, as in [22], nudging would have worked in the opposite direction.

To avoid confusing the reader with this somewhat subtle manoeuvre, in this paper *generating 1-morphisms will never share the same height in the graphical calculus.* That is, 1-morphisms will exclusively be composites of the form $\mathrm{id} \otimes f \otimes \mathrm{id}$ where $f$ is a generating 1-morphism. Occasionally we may, purely to save some space, draw some at the same height, but then it is understood that their heights are meant to be *resolved* using the nudging equation (76), which can only ever be applied in a unique way.

The interchangor 2-isomorphisms satisfy coherence equations explicitly listed in [5]. The most important are the following diagrams (together with their 180° rotated versions), which must commute whenever they make sense:

(78)

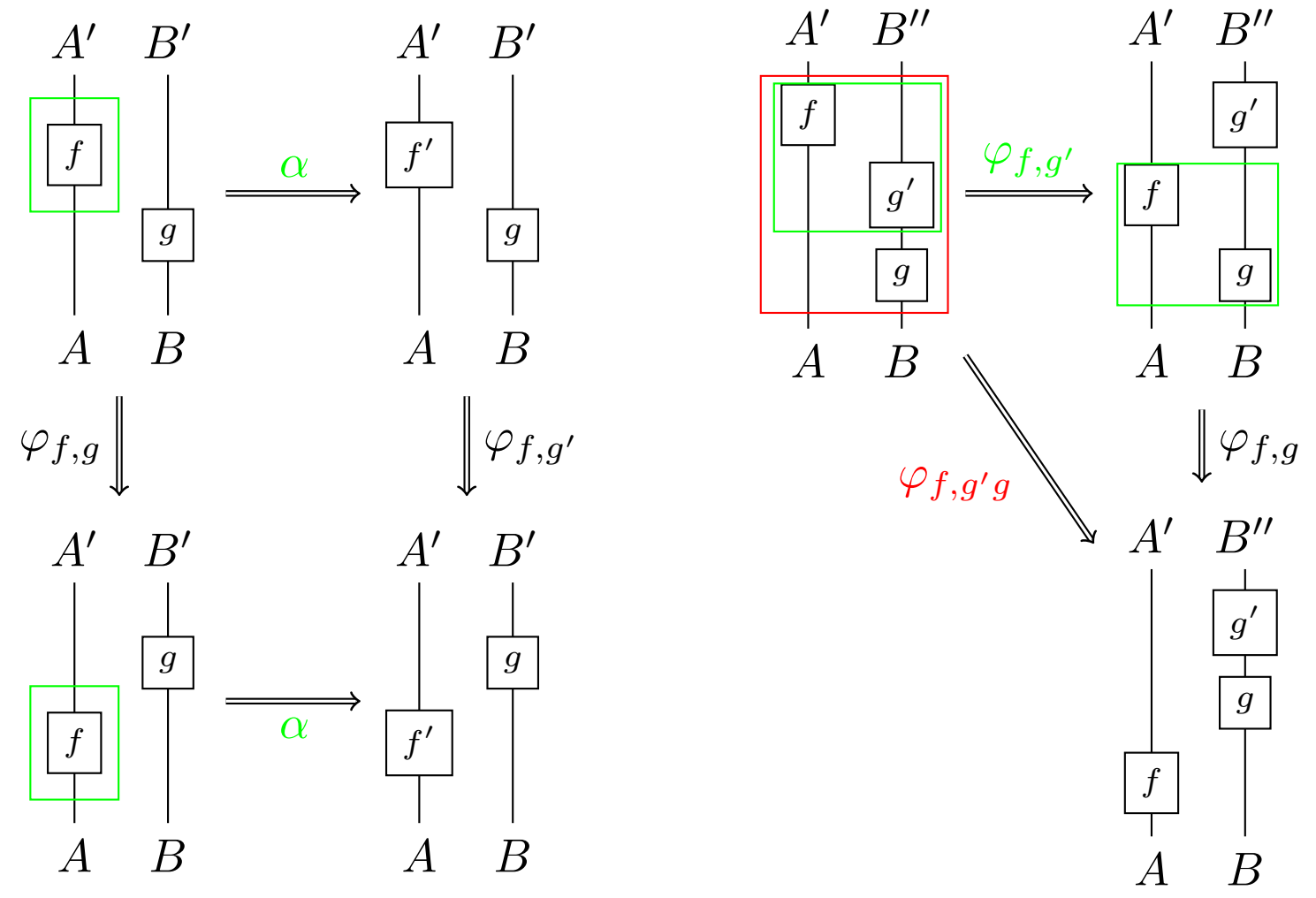

The braiding 1-isomorphisms $\gamma_{A,B} \colon A \otimes B \to B \otimes A$ in a quasistrict symmetric monoidal bicategory are drawn as:

(79)

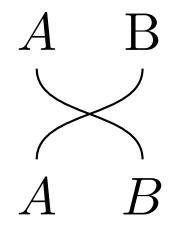

They are 1-isomorphisms *on-the-nose*, and we have:

(80)

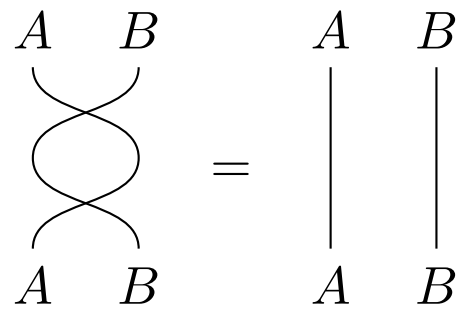

They must satisfy the following coherence equations *on-the-nose*, together with their 180° rotated versions[7]:

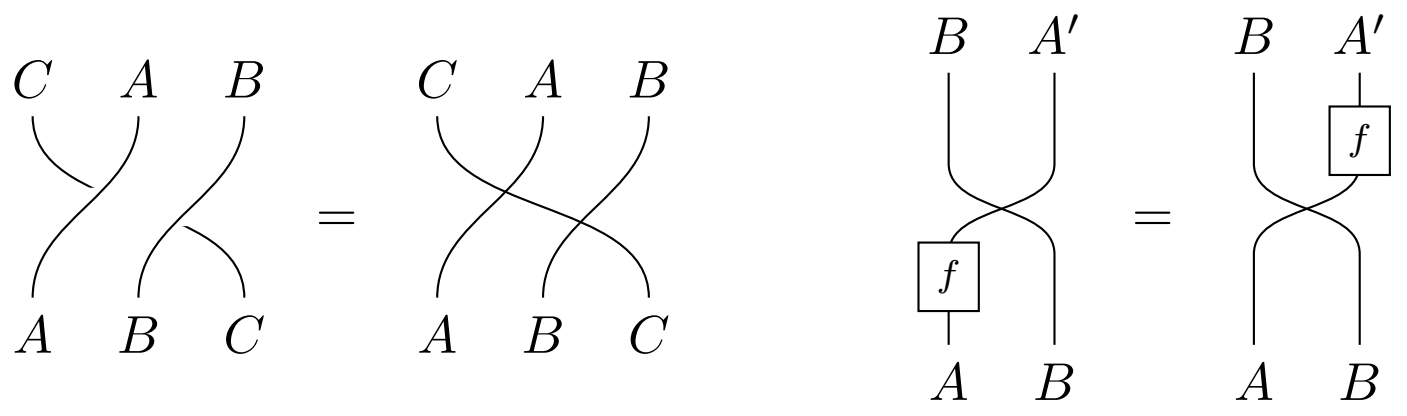

In addition, the following compatibility relation between the interchangor 2-isomorphisms $\phi_{f,g}$ and the braiding 1-isomorphisms $\gamma_{A,B}$ must hold, together with its rotated version (in particular the source and target of $\phi_{f,\beta_{B,C}}$ are required to coincide):

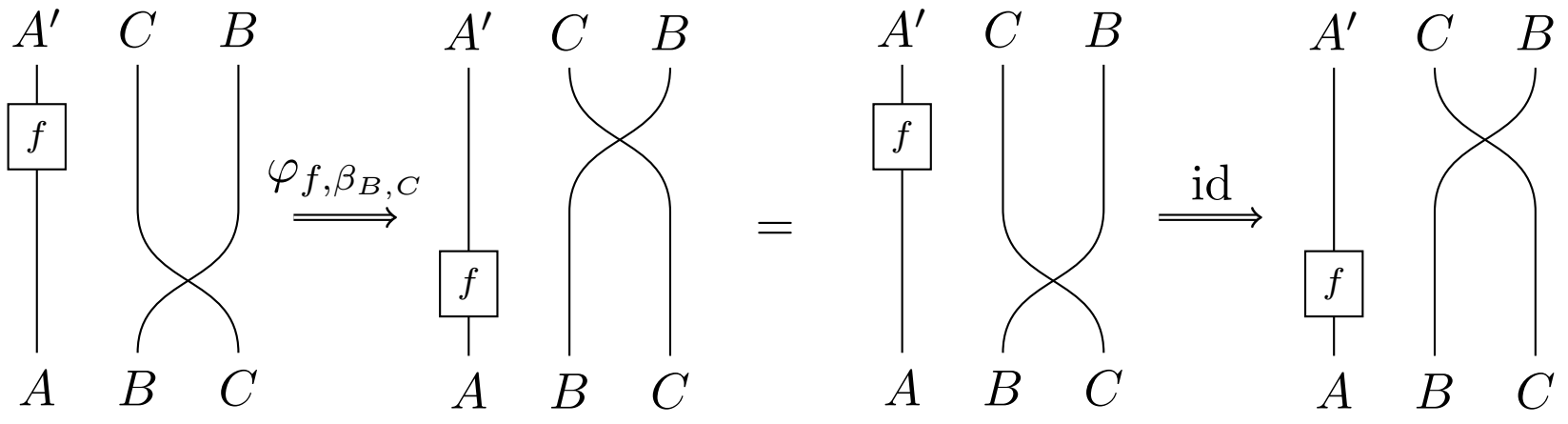

C.3. **Presentations of quasistrict symmetric monoidal 2-categories.** Let $\mathbb{G}_2$ be the category of two-dimensional globes, that is the category generated from the graph

$$0 \underset{t}{\overset{s}{\rightrightarrows}} 1 \underset{t}{\overset{s}{\rightrightarrows}} 2$$

subject to the relations $ss = st$ and $tt = ss$, and let $\mathsf{gSet}_2 := \mathrm{Fun}(\mathbb{G}_2^{\mathrm{op}}, \mathrm{Set})$ be the category of 2-globular sets. The axioms defining a quasistrict symmetric monoidal 2-category can be thought of as defining a (monotone, finitary) monad $T_{qs}$ on $\mathsf{gSet}_2$. To such a monad $T_{qs}$ there is a corresponding notion of a *3-computad*, which formalizes the notion of a generators-and-relations presentation for a quasistrict symmetric monoidal 2-category. We summarize the definition below; for a more precise account see [40].

**Definition 62.** A *presentation* of a quasistrict symmetric monoidal 2-category is a quasistrict symmetric monoidal 3-computad $\mathcal{H} = (G_0, G_1, G_2, R)$ in the sense of [40]. That is, it consists of:

- a set $G_0$ of *generating objects*;
- a set $G_1$ of *generating 1-morphisms*, whose source and target objects are words in the generating objects;
- a set $G_2$ of *generating 2-morphisms*, whose source and target 1-morphisms are sentences in the generating 1-morphisms (that is, composites of 1-morphisms generated by tensor product $\otimes$, composition $\circ$ of 1-morphisms and braiding generators $\gamma_{A,B}\colon A \otimes B \to B \otimes A$);

[7]The first equation states that $\gamma_{A\otimes B,C} = (\gamma_{A,C} \otimes 1_B) \circ (1_A \otimes \gamma_{B,C})$.

- a set $R$ of *relations* $(\alpha, \beta)$ between equivalence classes of paragraphs $\alpha$ and $\beta$ in the generating 2-morphisms. By a *paragraph in the generating 2-morphisms* we mean a composite of 2-morphisms from $G_2$ generated by tensor product $\otimes$, horizontal $*$ and vertical $\circ$ composition of 2-morphisms, and interchangor 2-isomorphisms $\varphi_{f,g}$. The equivalence relation is the finest such that the axioms of a (free) quasistrict symmetric monoidal 2-category are satisfied.

We write $\mathbf{F}(\mathcal{H})$ for the free quasistrict symmetric monoidal 2-category presented by $\mathcal{H}$.

So, an object of $\mathbf{F}(\mathcal{H})$ is a word (unparenthesized tensor product) in the generating objects. A 1-morphism of $\mathbf{F}(\mathcal{H})$ is a sentence (in the sense above) in the generating 1-morphisms. A 2-morphism is an equivalence class of paragraphs (in the sense above) in the generating 2-morphisms, where the relations are generated by $R$.

C.4. **Graphical calculus in the bordism setting.** This paper has to do with the presentation $\mathcal{F}$ defined in Appendix A, and the presentations $\mathcal{R}$, $\mathcal{M}$, and $\mathcal{O}$, defined in Section 2. Let $\mathcal{H}$ refer to one of these presentations. We will perform computations in $\mathbf{F}(\mathcal{H})$ using a version of the wire diagrams graphical calculus from [5] adapted to our context. They all have a single generating object $X$, which we will draw as a circle:

So, the objects of $\mathbf{F}(\mathcal{H})$ are arbitrary tensor products of $X$ — including the empty tensor product, which is the unit object 1 and is drawn with 'invisible ink' in the graphical calculus (which is fine since it is a *strict* unit). All the presentations have four generating 1-morphisms:

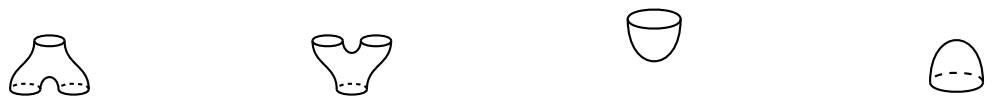

For instance, the 'pair-of-pants' 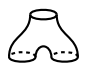 is a generating 1-morphism $M \colon X \otimes X \to X$:

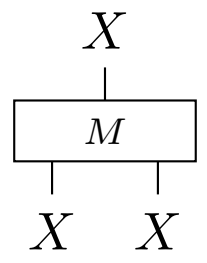

drawn as 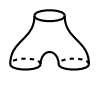.

Similarly, the 'cup' 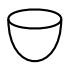 is a generating 1-morphism $1 \to X$:

drawn as

So, the 1-morphisms of $\mathbf{F}(\mathcal{H})$ are generated under $\circ$ and $\otimes$ from the above generators, together with the *braiding 1-isomorphisms* $\gamma_{X,X} \colon X \otimes X \to X \otimes X$

(81) 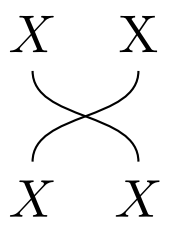

drawn as 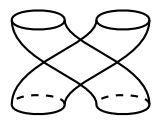.

The presentations differ in their generating 2-morphisms, but some generating 2-morphisms are common to all of them. For instance, there is the *associator* $\alpha\colon M \circ (M \otimes \mathrm{id}) \Rightarrow M \circ (\mathrm{id} \otimes M)$:

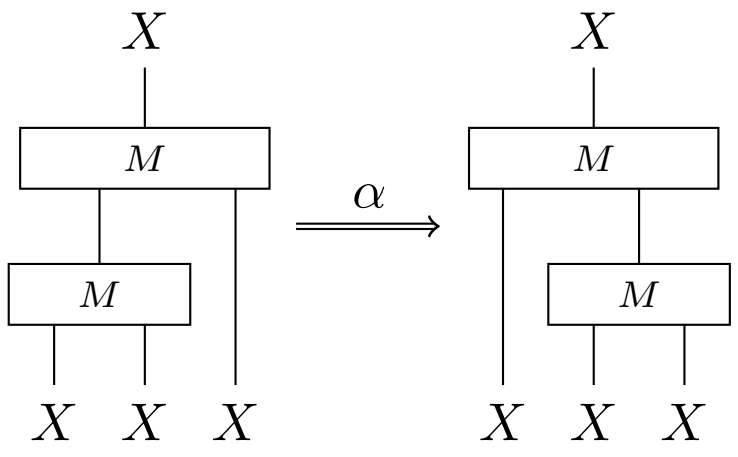

drawn as

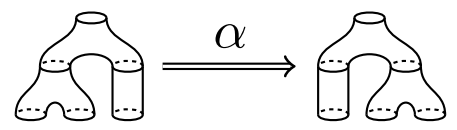

So, the 2-morphisms of $\mathbf{F}(\mathcal{H})$ are generated under $\circ$, whiskering, and tensor product by the 2-morphism generators, as well as the *interchangor* generators $\varphi_{f,f'}$ where $f$ and $f'$ are 1-morphism generators. For instance, in the pentagon relation (5), the inverse of the interchangor $\varphi_{M,M}$ appears:

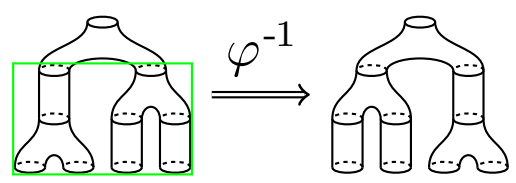

We will usually write the interchangors simply as $\varphi$, since the subscripts will usually be clear from the picture. Sometimes we will write $\varphi$ for a composite of interchangors; there is no ambiguity due to the coherence equations (78) as well as the further coherence equations listed in [5].

The *braiding generating 2-morphism* $\beta\colon M \Rightarrow M \circ \gamma_{X,X}$ will be another generating 2-morphism:

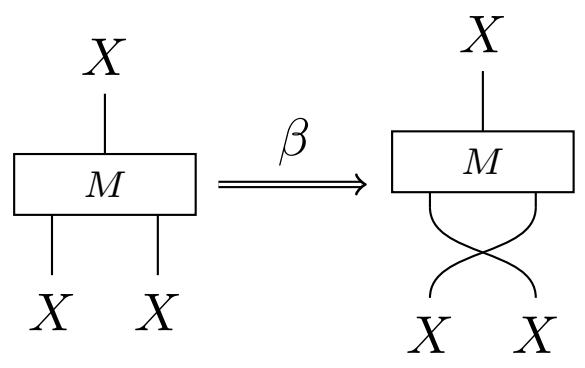

drawn as

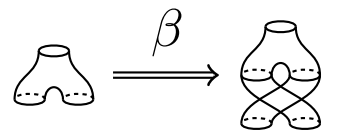

.

Note that $\beta$ is an explicit generating 2-morphism in the presentation and should not be confused with the braiding 1-morphism $\gamma_{X,X}$ from (81) which comes 'for free' in any 3-computad! One of the first relations will be that $\beta$ has an inverse

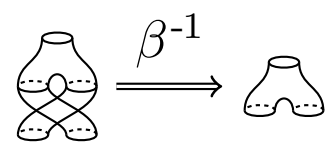

.

We will sometimes use the fact (80) that $\gamma_{X,X}^2 = \mathrm{id}$ to write $\beta^{-1}$ in the form

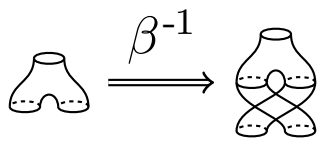

, i.e.

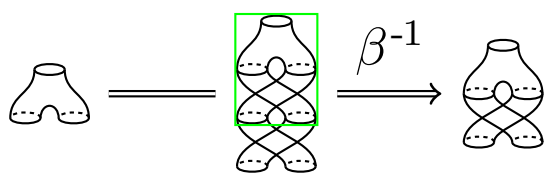

.

One property of the graphical calculus that we will use repeatedly is that if a pair of 2-morphisms have disjoint 1-morphisms as source, then they commute. For instance,

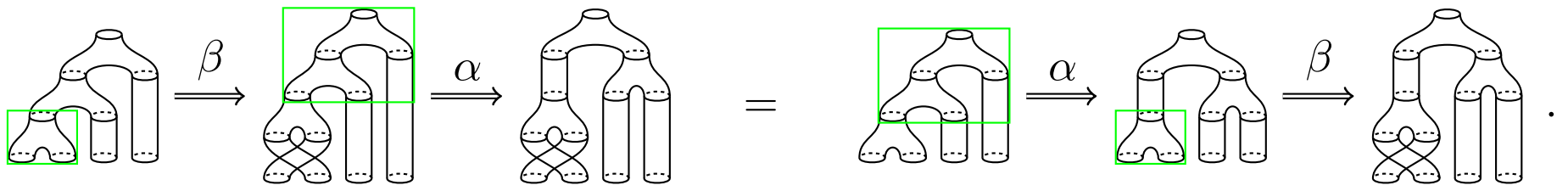

This follows since $\mathbf{F}(\mathcal{H})$ is a strict 2-category. Sometimes we will refer to this as *naturality*.

Mathematics Division, University of Stellenbosch
*Current address*: Mathematical Institute, University of Oxford
*E-mail address*: brucehbartlett@gmail.com

Mathematical Institute, University of Oxford
*E-mail address*: cdouglas@maths.ox.ac.uk

Max Planck Institute for Mathematics, Bonn
*E-mail address*: schommerpries.chris@gmail.com

Department of Computer Science, University of Oxford
*E-mail address*: jamie.vicary@cs.ox.ac.uk